\documentclass[12pt]{amsart}

\title[Mirror stabilizers]{Mirror stabilizers for lattice complex hyperbolic triangle groups}

\author{Martin Deraux}

\usepackage{url}
\usepackage{amssymb}
\usepackage{a4wide}
\usepackage{graphicx}

\newcommand{\Sc}{\mathcal{S}}
\newcommand{\Tc}{\mathcal{T}}
\newcommand{\Bc}{\mathcal{B}}
\newcommand{\RH}{H_{\mathbb{R}}}
\newcommand{\CH}{H_{\mathbb{C}}}
\newcommand{\CP}{P_{\mathbb{C}}}
\newcommand{\IP}{\mathbb{P}}
\newcommand{\IR}{\mathbb{R}}
\newcommand{\IC}{\mathbb{C}}
\newcommand{\IZ}{\mathbb{Z}}
\newcommand{\IN}{\mathbb{N}}
\newcommand{\IQ}{\mathbb{Q}}
\newcommand{\inn}[2]{\langle #1,#2\rangle}
\newcommand{\br}{\textrm{br}}
\newcommand{\tr}{\textrm{Tr}}

\newtheorem{thm}{Theorem}[section]
\newtheorem{prop}[thm]{Proposition}

\newtheorem{dfn}[thm]{Definition}

\theoremstyle{remark}
\newtheorem{rmk}[thm]{Remark}

\def\cqfd{\mbox{}\nolinebreak\hfill$\Box$\medbreak\par}
\newenvironment{pf}{\noindent\textbf{Proof:}}{\cqfd}

\address{Martin Deraux, Univ. Grenoble Alpes, CNRS, IF, 38000 Grenoble, France; Sorbonne
  Université and Université de Paris, CNRS, INRIA, IMJ-PRG, Ouragan,
  F-75005 Paris, France }
\email{martin.deraux@univ-grenoble-alpes.fr}

\date{January 18, 2023}

\begin{document}

\begin{abstract}
  For each lattice complex hyperbolic triangle group, we study the
  Fuchsian stabilizer of a reprentative of each group orbit of mirrors
  of complex reflections. We give explicit generators for the
  stabilizers, and compute their signature in the sense of Fuchsian
  groups. For some of the triangle groups, we also find explicit pairs
  of complex lines such that the union of their stabilizers generate
  the ambient lattice.
\end{abstract}

\maketitle

\section{Introduction}

Recall that the complex hyperbolic plane $\CH^2$ is a Hermitian
symmetric space of constant negative holomorphic sectional curvature
$-1$. Its real sectional curvatures satisfy
$-1\leq K_{sect}\leq -1/4$, where the lower and upper bounds are
realized by tangent planes of totally geodesic copies of $\CH^1$
(complex lines) and $\RH^2$ (Lagrangians), respectively.  Its group of
holomorphic isometries is isomorphic to $G=PU(2,1)$ (and this has
index two in the full isometry group, which includes antiholomorphic
transformations).

For each complex line $L$, there is an obvious embedding of $U(1,1)$
into $U(2,1)$ whose image preserves $L$; this gives a description of
$Stab_G(L)$ as a central extension of $PU(1,1)$ by the fixed point
stabilizer, which is the subgroup of complex reflections fixing $L$.

Each complex line is the fixed point set a 1-parameter family of
isometries, called complex reflections (beware these may have
arbitrary order), and we call a subgroup $\Gamma\subset G$ a
\emph{complex reflection group} if it is generated by complex
reflections. Complex reflection groups that can be generated by 3
complex reflections are called \emph{complex triangle groups}.

We call a subgroup $\Gamma$ of $G$ a lattice if it is discrete, and if
the quotient $\Gamma\backslash \CH^2$ has finite volume. In a sense,
the discreteness assumption says $\Gamma$ cannot be too large, whereas
the finite volume assumption says that $\Gamma$ cannot be too
small\dots which makes it tricky to find groups satisfying both
conditions at the same time. Lattices do exist though, due to very
general results of Borel~\cite{borel}; the basic idea is to use
``arithmetic methods'', which we think of as a souped up version of
the fact that $\IZ\subset\IR$ is discrete, and $\IR/\IZ$ has finite
Lebesgue measure.

Another important (but difficult) method to construct lattices is to
use (very wisely chosen) complex reflection groups. This was the
approach taken by Mostow~\cite{mostow-pacific} to construct the first
examples of non-arithmetic lattices (lattices that cannot be obtained
from the arithmetic construction). There are many other approaches,
using various ideas around uniformization
(see~\cite{bhh},~\cite{deligne-mostow},~\cite{chl} for instance), but
in this paper we focus only on the groups constructed
in~\cite{dpp2}; the latter paper can be thought of as one possible
generalization of the basic ideas developed by Mostow
in~\cite{mostow-pacific}.

The techniques that are used in~\cite{dpp2} in order to generalize
Mostow's construction are very intricate, and we will not go into
almost any of them. We freely use the results~\cite{dpp2} as a black
box, keeping technical details to a bare minimum.

For convenience of the exposition, we refer to the lattices
constructed in~\cite{dpp2} as \emph{the lattice complex hyperbolic
  triangle groups} (the terminology may be a bit misleading, because
there may actually be more lattice complex hyperbolic triangle groups,
but we expect our search in~\cite{dpp-census}
to be fairly exhaustive).

The basic point is that the lattice complex hyperbolic triangle
groups come in three families of groups
$$
\Sc(p,\tau),\quad \Tc(p,{\bf T}),\quad \Gamma(p,t)
$$
where $p\in\IN$, $p\geq 2$, $\tau\in\IC$, ${\bf T}\in\IC^3$,
$t\in\IQ$ are wisely chosen parameters (see the list in the appendix
of~\cite{dpp2}).

The lattices $\Gamma(p,t)$ are those that were studied by Mostow
in~\cite{mostow-pacific}, and we do not consider them in this paper
since they have been given many other descriptions
(see~\cite{deligne-mostow},~\cite{thurston-shapes} for instance).

For each lattice complex hyperbolic triangle group, the reader can
find explicit matrices generating the group, a description of a
fundamental domain for its action on $\CH^2$ and a presentation in
terms of generators in relations (see~\cite{dpp2},~\cite{spocheck};
see also~\cite{neat-subs} for more recent development).

In what follows, let $\Gamma$ denote a lattice complex hyperbolic
triangle group, and let $X=\Gamma\backslash\CH^2$ denote the
quotient. The main goal of the present paper is to study some detailed
properties of the (1-dimensional part of the) branch locus of the
quotient map $\CH^2\rightarrow X$. Indeed, from the results
in~\cite{dpp2}, one can extract a list $A_1,\dots,A_r$ of
representatives of conjugacy classes of complex reflections in the
group. For each $j$, $Stab_{\Gamma}(A_j)$ is then a central extension
of a lattice $F_j\subset PU(1,1)$. Our goal is to describe these
groups $F_j$; for each lattice complex hyperbolic triangle group, we
list the group orbits of mirror stabilizers, and for each of them, we
\begin{itemize}
\item find an explicit generating set and
\item determine the signature of the corresponding Fuchsian group, its
  trace field and its arithmeticity
\end{itemize}
From the results in the paper, one can write explicit presentations
for the mirror stabilizers (see equation~\eqref{eq:pres-stab-s1} for
one example of such a presentation).

We hope that our results give useful information for several
purposes. One purpose could be to try and understand the quotient
$X=\Gamma\backslash \CH^2$ and branching behavior of the quotient map
$\CH^2\rightarrow X$; in some cases the quotient map is completely
understood (see~\cite{deraux-klein} and~\cite{deraux-abelian}).
Another purpose is to study possible hybrid structures for complex
hyperbolic lattice triangle groups
(see~\cite{wells},~\cite{falbel-pasquinelli}); we will give some
explicit examples in section~\ref{sec:hybrids}.  Finally, this
information is useful in order to study the counting and
equidistribution properties of closed geodesic that are reflecting on
mirrors in $X= \Gamma\backslash \CH^2$, see \cite[\S 5.2]{epps}, which
was the initial motivation for writing this paper.

Part of the analysis of Fuchsian subgroups of lattice complex
hyperbolic triangle group was achieved by Sun in~\cite{lijie}. In
that paper, only some of the groups were treated, and only some of the
conjugacy classes of complex reflections, but we use the same
method. Information on the arithmeticity of the subgroups listed
by Sun can be found in a very recent paper by Jiang,~Wang
and~Yang~\cite{jiang-wang-yang}.

\section{Basic complex hyperbolic geometry}

Let $V$ be an (n+1)-dimensional complex vector space and
$h:V\times V\rightarrow\IC$ be a Hermitian form of signature
$(n,1)$. We write $\inn{X}{Y}$ for $h(X,Y)$ and $||X||^2$ for $h(X,X)$
(which is real, but not necessarily positive). Denote by $V^-$
(resp. $V^0$, $V^+$) the set vectors $X\in V$ with $||X||^2<0$
(resp. $=0$, $>0$).

The isometry group of the form is given by
$$
U(h)=\{ A\in GL(V) : \forall X,Y\in V, h(AX,AY)=h(X,Y) \},
$$
which we think of as a real Lie group. By Sylvester's law of inertia,
that Lie group is independent of the choice of the Hermitian form (up
to isomorphism of Lie groups).

Consider the projectivization map
$\pi:V\setminus\{0\}\rightarrow \IP(V)$, and define $\CH^2=p(V^-)$,
$\partial_\infty \CH^n=p(V^0\setminus\{0\})$.
It is a standard fact
that that $PU(h)$ acts transitively on $\CH^n$, and also on
$\partial\CH^n$ (in fact it also acts transitively on $p(V^+)$); these
homogeneous spaces are very different, since the stabilizer of a point
is compact only in the case of $\CH^n$.

It is a standard fact that $\CH^n$ carries a $PU(h)$-invariant
Riemannian metric (whereas $\partial_\infty\CH^n$ does not, nor does
$p(V^+)$), which makes it a Hermitian symmetric space. The interested
reader can find an expression for that metric
in~\cite{kobayashi-nomizu} for instance. The corresponding integrated
distance formula is easy to write down:
\begin{equation}\label{eq:distance}
  \cosh(\frac{1}{2}d({\bf X},{\bf Y}))=\frac{|\inn{X}{Y}|}{\sqrt{||X||^2||Y||^2}},
\end{equation}
where $X\in V$ (resp. $Y\in V$) is any lift of ${\bf X}\in\CH^n$
(resp. ${\bf Y}\in\CH^n$). The factor $\frac{1}{2}$ is included so
that the metric is scaled to have holomorphic sectional curvature
identically equal to $-1$. 

There are totally geodesic copies of $\CH^k$ inside $\CH^n$ for any
$k\leq n$, obtained by choosing a $(k+1)$-dimensional complex subspace
$W\subset V$ such that $h|_W$ has signature $(k,1)$. Similary, there
are totally geodesic copies of $\RH^k$, obtained by taking a totally
real (k+1)-dimensional subspace $R\subset V$, i.e. the $\IR$-span of
$k+1$ vectors $v_0,\dots,v_k$ with $\inn{v_j}{v_k}\in\IR$ for all
$j,k$. It is a well known fact that every complete totally geodesic
submanifold of $\CH^n$ is among the ones just described (see section
2.5 of~\cite{chen-greenberg}). 

Note that $\CH^n$ has non-constant curvature for all $n\geq 2$. More
precisely, we have the following (see~\cite{chen-greenberg}
or~\cite{goldman-book}).
\begin{prop}\label{prop:curvature}
  The real sectional curvatures of $\CH^n$ are contained in
  $[-1,-1/4]$. The real 2-planes with curvature $-1$ are the tangent
  planes to totally geodesic copies of $\CH^1$, whereas the ones with
  curvature $-1/4$ are the tangent planes to totally geodesic copies
  of $\RH^2$.
\end{prop}

Since $h$ is non-degenerate, it defines a bijection between the
Grassmannians of complex $k$-planes and $(n-k)$-planes, where $W$ is
in correspondence with
$W^\perp=\{X\in V:\forall Y\in W,\inn{X}{Y}=0\}$. As an important
special case, when $k=n-1$, we get a parametrization of totally
geodesic copies of $\CH^{n-1}$ by positive lines, i.e. by $p(V^+)$.
For $X\in V^+$, we will often abuse notation and write $X^\perp$ for
$p(V^-\cap X^\perp)$.

\begin{prop}\label{eq:angle}
  Let $X,Y\in V^+$. The complex hyperplanes $X^\perp$ and $Y^\perp$
  intersect in $\CH^2$ if and only if
  $$
  c:=\frac{|\inn{X}{Y}|}{\sqrt{||X||^2||Y||^2}}<1.
  $$
  If they intersect in $\CH^2$, they do so at a constant angle
  $\alpha$, and we have $c=\cos(\alpha)$.
\end{prop}

Given a positive vector $U\in V$, it is easy to see that, for any
$\zeta\in\IC$ with $|\zeta|=1$, the formula
\begin{equation}\label{eq:cxrefl}
  R_{U,\zeta}(X)=X+(\zeta-1)\frac{\inn{X}{U}}{\inn{U}{U}}U,
\end{equation}
defines an element of $U(h)$. It fixes pointwise $U^\perp$, and
rotates by $\arg(\zeta)$ in the complex directions orthogonal to its
fixed point set. In terms of linear algebra, it is a diagonalizable
linear transformation with eigenvalues $(\zeta,1,\dots,1)$, such that
$h$ restricts to a form of signature $(n-1,1)$ on the 1-eigenspace.

We call the corresponding isometry $R_{U,\zeta}\in PU(h)$ a
\emph{complex reflection} with mirror $U^\perp$. Replacing $U$ by
$\lambda U$ for $\lambda\in\IC$, leaves $R_{U,\zeta}$ unchanged, so we
may assume $\inn{U}{U}=1$. The complex number $\zeta$ is called the
\emph{multiplier} of the complex reflection, $arg(\zeta)$ is called
the \emph{angle} of the rotation.

Formula~\eqref{eq:cxrefl} makes sense if $U$ is a negative vector, in
which case we call $R_{U,\zeta}$ a \emph{complex reflection
  in a point}; in that case $p(U)$ is its only fixed point in $\CH^n$
(in particular, for $\zeta=-1$, we get isometric involutions with a
single fixed point).

\begin{rmk}\label{refl-commute}
  It is a standard fact that two complex reflections
  $R_{U_1,\zeta_1}$, $R_{U_2,\zeta_2}$ commute if and only if
  $U_1^\perp=U_2^\perp$ or $\inn{U_1}{U_2}=0$.
\end{rmk}

We will use the usual classification of isometries of $PU(h)$ into
elliptic, parabolic and loxodromic elements (see section~6.2
of~\cite{goldman-book}). Among elliptic elements, we distinguish regular
ones whose matrix representatives have $(n+1)$ distinct
eigenvalues. Non-regular elliptic elements are sometimes called
special elliptic.
There are also fine classifications for parabolic and loxodromic
elements, but we will not need that in the present paper.

From now on, we restrict to the case $n=2$. In this case, every
elliptic isometry is either regular elliptic, a complex reflection, or
a complex reflection in a point.

Now let $A\in PU(h)$ be regular elliptic. Choose a basis
$\mathcal{E}=e_0,e_1,e_2$ of $V$ that diagonalizes $A$; write
$\zeta_0,\zeta_1,\zeta_2$ for the corresponding (distinct)
eigenvalues. Since $A$ is elliptic, we must have $|\zeta|=1$ for all
$j=1,2,3$. Moreover, since eigenvectors with distinct eigenvalues are
orthogonal, $\mathcal{E}$ is $h$-orthogonal, and because of the
signature of $h$, one and only one of the basis vectors is
$h$-negative. By reordering and rescaling, we may assume
$||e_0||^2=-1$ and $||e_1||^2=||e_2||^2=1$. Then ${\bf e_0}=p(e_0)$ is
the isolated fixed point of $A$ in $\CH^2$. The complex 2-planes
$L_1=\IC\{e_0,e_1\}$ and $L_2=\IC\{e_0,e_2\}$ project to complex lines
${\bf L_1}$ and ${\bf L_2}$, each being $A$-invariant (and these are
the only complex lines through ${\bf e_0}$ that are $A$-invariant). In
${\bf L_j}$, $A$ acts as a rotation by $arg(\zeta_j/\zeta_0)$.

Conversely, let ${\bf L_j}=p(V^-\cap U_j^\perp)$ for positive vectors
$U_j$. We assume ${\bf L_1}$ and ${\bf L_2}$ are orthogonal,
i.e. $\inn{U_1}{U_2}=0$ (cf. formula~\eqref{eq:angle}). Then for any
$\zeta_j$ with $|\zeta_j|=1$, the product
$$
  R_{U_1,\zeta_1}R_{U_2,\zeta_2}
$$
of commuting complex reflection is elliptic. It is regular elliptic if
and only if $\zeta_1\neq \zeta_2$.

Note that if $A^k$ is special elliptic for some $k\in\IN^*$, then
$A^k$ is either a complex reflection in a point, or a complex
reflection. In the latter case, the mirror of $A^k$ must be either
${\bf L_1}$ or ${\bf L_2}$.

\section{Lattice complex hyperbolic triangle groups}\label{sec:list}

We refer the reader to~\cite{dpp2} for a detailed description of the
lattices $\Sc(p,\tau),\quad \Tc(p,{\bf T})$. Alternatively, the reader
can download matrices in Magma format from~\cite{neat-subs}.

Here we only give explicit generators $R_1,R_2,R_3$ for
$\Tc(p,{\bf T})$, where $p\in\IN^*$ and
${\bf T}=(\rho,\sigma,\tau)\in\IC^3$ (these specialize to generators
for $\Sc(p,\tau)$ by taking $\rho=\sigma=\tau$).

\begin{dfn}\label{dfn:Tc}
  For ${\bf T}=(\rho,\sigma,\tau)\in\IC^3$, the group $\Tc(p,{\bf T})$
  is the subgroup of $U(H)$ generated by $R_1,R_2$ and $R_3$ (see
  equations~\eqref{eq:form} and~\eqref{eq:generators}), where
  $u=e^{2pi/3p}$, $\alpha=2-u^3-\bar u^3$, $\beta_1=(\bar u^2-u)\rho$,
  $\beta_2=(\bar u^2-u)\sigma$, $\beta_3=(\bar u^2-u)\tau$.
\end{dfn}

\begin{equation}\label{eq:form}
H=\left(\begin{matrix}
  \alpha & \beta_1 & \overline{\beta}_3\\
  \overline{\beta}_1&\alpha&\beta_2\\
  \beta_3&\overline{\beta}_2&\alpha
\end{matrix}\right)
\end{equation}

{\small
\begin{equation} \label{eq:generators}
  R_1=\left(\begin{matrix}
  u^2 & \rho & -u\overline{\tau}\\
  0 & \bar u & 0\\
  0 & 0 & \bar u
  \end{matrix}\right);\ 
  R_2=\left(\begin{matrix}
  \bar u & 0 & 0\\
  -u\bar \rho & u^2 & \sigma\\
  0 & 0 & \bar u
\end{matrix}\right);\ 
  R_3=\left(\begin{matrix}
  \bar u & 0 & 0\\
  0 & \bar u & 0\\
  \tau & -u\bar \sigma & u^2
\end{matrix}\right)
\end{equation}
}

\begin{rmk}
  \begin{enumerate}
  \item If $e_1,e_2,e_3$ denotes the standard basis for $\IC^3$, then
    in terms of the notation of formula~\eqref{eq:cxrefl}, we have
    $R_j=\bar u R_{e_j,u^3}$, so the matrices $R_j$ can easily be
    reconstructed from the data of $H$.
  \item The factor $\bar u$ is included to get matrices in $SU(H)$
    rather than just $U(H)$.
  \end{enumerate}
\end{rmk}

When $\rho=\sigma=\tau$, the map
$$
J=\left(\begin{matrix}
  0 & 0 & 1\\
  1 & 0 & 0 \\
  0 & 1 & 0
\end{matrix}\right)
$$
is a regular elliptic element of order 3 in $U(H)$, and we have
$$
JR_1J^{-1}=R_2,\quad JR_2J^{-1}=R_3.
$$

\begin{dfn}
  For $\tau\in\IC$, $\Sc(p,\tau)$ is the subgroup of $U(H)$ generated
  by $R_1$ and $J$, where $u=e^{2pi/3p}$, $\alpha=2-u^3-\bar u^3$,
  $\beta_1=\beta_2=\beta_3=(\bar u^2-u)\tau$.
\end{dfn}

\begin{rmk}
  For an arbitrary choice of the parameters $p$ and ${\bf T}$, the
  group $\Tc(p,{\bf T})$ is usually not discrete; moreover, the
  Hermitian form on $\IC^3$ given by the Hermitian matrix $H$ is not
  even of signature $(2,1)$ in general.
\end{rmk}

We will only consider the groups with parameters as in
Table~\ref{tab:list}, which gives representatives of all the
lattices that are known to be complex hyperbolic triangle groups,
apart from Mostow groups. Note that some more parameters are known to
give lattices (in other words, some lattices can have several triangle
group descriptions, see section~7 of~\cite{dpp2}).

\begin{table}\label{tab:list}
\begin{tabular}{c|c}
  $\tau$ or ${\bf T}$                                               & $p$\\
  \hline
  $\sigma_1 = -1+i\sqrt{2}$                                         & 3,4,6\\
  $\bar\sigma_4=\frac{-1-i\sqrt{7}}{2}$                            & 3,4,5,6,8,12\\
  $\sigma_5=e^{-\pi i/9}\frac{\sqrt{5}+i\sqrt{3}}{2}$              & 2,3,4\\
  $\sigma_{10}=\frac{1+\sqrt{5}}{2}$                               & 3,4,5,10\\
  \hline
  ${\bf S_2}=(1+\frac{-1+i\sqrt{3}}{2}\frac{1+\sqrt{5}}{2},1,1)$   & 3,4,5\\
  ${\bf E_2}=(\sqrt{2},\frac{-1+i\sqrt{3}}{2},\sqrt{2})$           & 3,4,6\\
  ${\bf H_1}=(\frac{-1+i\sqrt{7}}{2},e^{-4\pi i/7},e^{-4\pi i/7})$ & 2\\
  ${\bf H_2}=(-1-e^{-2\pi i/5},e^{4\pi i/5},e^{4\pi i/5})$         & 2,3,5\\
\end{tabular}
\caption{Values of the parameters for the known lattice complex
  hyperbolic triangle groups}
\end{table}

\section{Braid relations between complex reflections}\label{sec:braiding}

The structure of each group (and of a fundamental domain for its
action on $\CH^2$) is conveniently encoded by \emph{braid relations}
between suitable conjugates of the generators. Recall that two
elements $a,b$ of a group satisfy a braid relation of length
$n\in\IN^*$ if
\begin{equation}\label{eq:braiding}
(ab)^{n/2}=(ba)^{n/2}.
\end{equation}
When $n$ is odd, by convention
$(ab)^{n/2}=a\cdot b\cdot a\cdot b\cdot a\cdots b\cdot a$ ($n$
factors). When equation~\eqref{eq:braiding} holds, we write
$\br_n(a,b)$; if $n$ is the smallest $n$ such that $\br_n(a,b)$, we
say the braid length of $a$ and $b$ is $n$, and we write $\br(a,b)=n$.

We will use braid relations only between complex reflections. An
important result is the following (see Proposition~2.5
of~\cite{dpp2}).
\begin{prop}\label{prop:center}
  Suppose $R_1$ and $R_2$ are complex reflections in $PU(2,1)$ with
  distinct mirrors, and multipliers $\zeta_1=e^{2\pi i/p_1}$ and
  $\zeta_2=e^{2\pi i/p_2}$ respectively ($p_j\in\IN$, $p_j\geq
  2$). Suppose $\br_n(R_1,R_2)$ for some $n\in\IN$, $n\geq 2$. Let
  $Z=(R_1R_2)^n$ if $n$ is odd, $Z=(R_1R_2)^{n/2}$ if $n$ is even.
  \begin{enumerate}
  \item $Z$ commutes with both $R_1$ and $R_2$.
  \item If $\br_n(R_1,R_2)$ with $n$ odd, then $R_1$ and $R_2$ are
    conjugate, in particular $p_1=p_2$.
  \end{enumerate}
\end{prop}
It is easy to see that $Z$ as defined in~\ref{prop:center} has a
repeated eigenvalue, hence it must be parabolic, a complex reflection
or a complex reflection in a point (recall that we work with $\CH^2$,
so our homogeneous coordinates live in a complex vector space of
dimension 3). The three cases are distinguished by looking at the
position with respect to $\CH^2$ of the intersection point in
projective space of the mirrors of $R_1$ and $R_2$. In fact $Z$ is a
complex reflection in a point (resp. parabolic, resp. a complex
reflection) if the mirrors intersect in $\CH^2$ (resp. in
$\partial_\infty\CH^2$, resp. outside $\overline{\CH}^2$).

Writing $R_j=R_{U_j,\zeta_j}$ (see formula~\eqref{eq:cxrefl}), the three
cases can be conveniently characterized computationally as follows. Consider
$X=U_1\boxtimes U_2$, which is the usual cross product of $U_1^*H$ and
$U_2^*H$; by construction $\inn{X}{U_1}=\inn{X}{U_2}=0$, so $\IC X$ is the
intersection of the mirrors in the projective plane. Then we have
\begin{prop}
  $Z$ is a complex reflection if $||X||^2>0$, parabolic if
  $||X||^2=0$, and a complex reflection in a point if $||X||^2<0$.
\end{prop}

When $Z$ is a complex reflection, part~(1) of
Proposition~\ref{prop:center} says that it commutes with both $R_1$
and $R_2$, so the mirror of $Z$ is the common perpendicular complex
line between the mirrors of $R_1$ and $R_2$.

The relevant braid relations for each group of section~\ref{sec:list}
appear in the Appendix of~\cite{dpp2}. One important point is that the
relations depend only on ${\bf T}$, not on $p$. In fact the
fundamental domains for lattices with a given value of ${\bf T}$ (but
different values of $p$) are almost the same, they differ only by a
simple truncation process (see the discussion in
section~\ref{sec:sides}).

\section{Sides of our fundamental domain}\label{sec:sides}

In this section, we briefly recall the construction of the sides of
the fundamental domains described in~\cite{dpp2}.

It is well known that there are no totally geodesic real hypersurfaces
in $\CH^2$ (this can be seen as a consequence of the fact that it has
non-constant real sectional curvature). In particular, the sides of
the fundamental domains constructed in~\cite{dpp2} are not totally
geodesic.

On the other hand, part of their skeleton is totally geodesic, namely
the 1-skeleton consists only of geodesic arcs, and some facets in the
2-skeleton are also totally geodesic.

Recall that each side of our domains is constructed inside a bisector,
i.e. the locus of points equidistant of two given points. The
geometric structure of bisectors (and even bisector intersections) is
by now well understood, see chapter~5 of~\cite{goldman-book}. Each
bisector $\Bc$ is topologically a 3-ball, and it has a foliation by
complex lines (called the complex slices of $\Bc$), and a singular
foliation by real planes (called the real slices) all intersecting in
a fixed real geodesic of $\CH^2$ (called the real spine of $\Bc)$. The
unique complex line containing that real geodesic is called the
complex spine of $\Bc$.

If $\Bc$ is a bisector and $p,q\in \Bc$, the real geodesic through $p$
and $q$ is contained in $\Bc$ if and only if $p,q$ are either in a
real or in a complex slice of $\Bc$. This is always true if one of the
two points is on the real spine. Our domains use the latter property,
namely each side is constructed as a (possibly truncated) geodesic
cone inside a bisector $\Bc$, the base of the cone being a polygon in
a complex slice of $\Bc$, and its apex is on the real spine of $\Bc$
(see Figure~4.2 of~\cite{dpp2} and also the pictures in the Appendix
of~\cite{dpp2} for examples of truncated cones).

The complex slices containing the base of the cones are given by well
chosen mirrors of complex reflections in the triangle group (in fact
suitable conjugates of the generating reflections $R_1$, $R_2$ and
$R_3$). The vertices of the base polygons are either intersections of
mirrors of reflections, or intersections with the common perpendicular
of two mirrors of reflections when some truncation is needed (see
section~\ref{sec:braiding}).

\begin{figure}[htbp]
  \centering
  \hfill\includegraphics{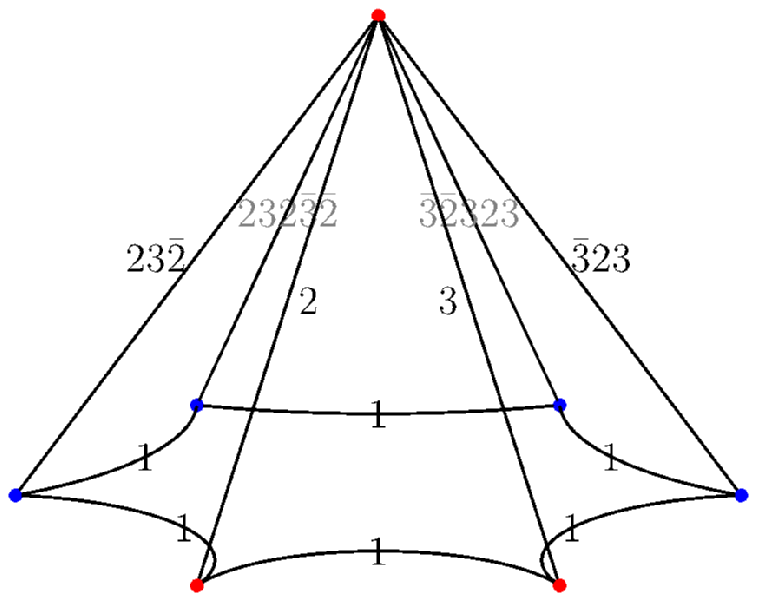}\hfill\includegraphics{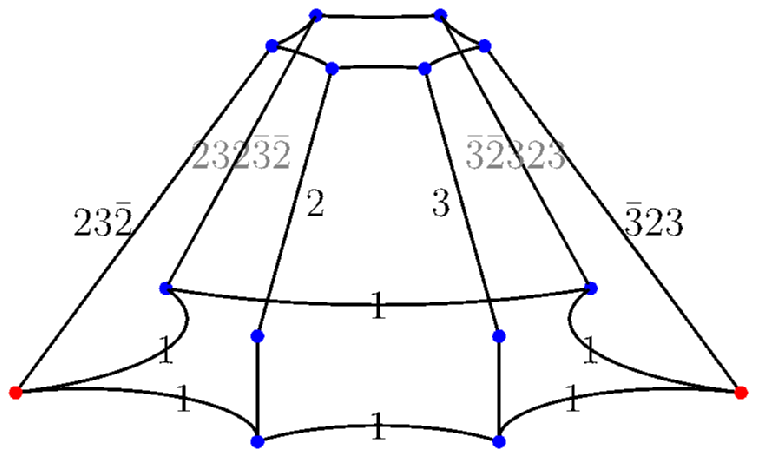}\hfill 
  \caption{Pictures of some sides for our fundamental domain for
    $\Sc(p,\sigma_1)$, $p=3$ (left) or $4$ (right)}\label{fig:s1}
\end{figure}

For example, in Figure~\ref{fig:s1}, we reproduce some figures from
the appendix of~\cite{dpp2}.
Let us focus first on the left half of the picture, which corresponds
to the case $p=3$, i.e. $R_j=R_{e_j,\zeta_3}$ with
$\zeta_3=e^{2\pi i/3}=\frac{-1+i\sqrt{3}}{2}$.  The labels indicate
that the base is a hexagon in the mirror of $R_1$ (which is
$e_1^\perp$), with vertices given by
$$
  e_1\boxtimes e_2,\quad
  e_1\boxtimes e_3,\quad
  e_1\boxtimes R_3^{-1}e_2,\quad
  e_1\boxtimes R_3^{-1}R_2^{-1}e_3,\quad
  e_1\boxtimes R_3^{-1}R_2^{-1}e_3,\quad
  e_1\boxtimes R_2e_3.
$$
Note for example that the mirror of $R_2R_3R_2^{-1}$ is the image
under $R_2$ of the mirror of $R_3=R_{e_3,\zeta_3}$, so it is
$R_2(e_3)^\perp$.  The apex of the cone is given by
$e_2\boxtimes e_3$.

Now consider the case $p=4$ (right half of the picture), which
correspond to generating reflections of the form $R_j=R_{e_j,i}$
($i^2=-1$). The corresponding side is similar (combinatorially it is
the same up to truncation). In that case, the base is still a hexagon
in $e_1^\perp$, but its vertices are now given by
$$
  e_1\boxtimes (e_1\boxtimes e_2),\quad
  e_1\boxtimes (e_1\boxtimes e_3),\quad
  e_1\boxtimes R_3^{-1}e_2,\quad
  e_1\boxtimes R_3^{-1}R_2^{-1}e_3,\quad
  e_1\boxtimes R_3^{-1}R_2^{-1}e_3,\quad
  e_1\boxtimes R_2e_3.
$$
The first vertex correspond to a truncation using the common
perpendicular complex line of $e_1^\perp$ and $e_2^\perp$, which is
$(e_1\boxtimes e_2)^\perp$. Note that due to truncation, there is also
an extra vertex, which is given by $e_2\boxtimes(e_1\boxtimes e_2)$.

The (right half of the) picture also indicates that the apex of the
cone $e_2\boxtimes e_3$ is not a vertex of that side, in fact the apex also
needs to be truncated; the truncation is done with the common
perpendicular complex line of $e_2^\perp$ and $e_3^\perp$, which is
$(e_2\boxtimes e_3)^\perp$. The truncation facet is then a hexagon on
that complex line, with vertices given by
$$
  m\boxtimes e_2,\quad
  m\boxtimes e_3,\quad
  m\boxtimes R_3^{-1}e_2,\quad
  m\boxtimes R_3^{-1}R_2^{-1}e_3,\quad
  m\boxtimes R_3^{-1}R_2^{-1}e_3,\quad
  m\boxtimes R_2e_3,
$$
where $m=e_2\boxtimes e_3$.

The 1-skeleton of the sides can be constructed from the 0-skeleton if
we know the combinatorics, since the 1-skeleton consists of geodesic
arcs. The 2-skeleton is more complicated to describe (see section 4.2
of~\cite{dpp2} for details); here we simply mention that the only
2-facets that lie in complex lines are the facets on the base of the
cone, as well as the truncation 2-facets that replace the apex of the
cone.

Note that in the group $\Sc(p,\sigma_1)$ corresponding to the pictures
in Figure~\ref{fig:s1}, we have $\br_6(R_2,R_3)$ (which is the reason
why we choose a hexagon in the base, see~\cite{dpp2}), and for $p=6$
$(R_2R_3)^3$ is a complex reflection fixing the top truncation facet.

In order to refer to sides, we will use the same notation as
in~\cite{dpp2} and use the notation
$$
[n]a;b,c
$$
where $a$, $b$ and $c$ are complex reflections such that $\br_n(b,c)$.
For the side written $[n]a;b,c$, $a$ stands for a reflection fixing
the base, and $b$, $c$ should be thought of as defining two of the
vertical edges of the cone. The full list of vertical edges is then
uniquely determined by $b$ and $c$; the number of elements in the list
is given by $n$, and the list of edges is determined by the following
list of reflections
$$
\dots,bcbc^{-1}b^{-1},bcb^{-1},b,c,c^{-1}bc,c^{-1}b^{-1}cbc,\dots
$$
given that this seemingly infinite list actually is $n$-periodic.

In order to describe a reflection in the group, we will use word
notation, writing $1$ for $R_1$, $2$ for $R_2$, $3$ for $R_3$, $\bar1$ for $R_1^{-1}$, $(123\bar2)^3$ for $(R_1R_2R_3R_2^{-1})^3$, etc.

For example the side described in Figure~\ref{fig:s1} is simply
written as $[6]1;2,3$.

The Appendix of~\cite{dpp2} contains a picture of (isometry types of)
all sides of a fundamental domains, for every triangle group with
parameter in~\ref{tab:list}. All pictures are drawn so that 2-facets
that are in complex lines are drawn in a horizontal plane (either as
the base of the side, or as the top side when there is truncation at
the top). The base sides are always fixed by some complex reflection
(in fact by a conjugate of one of the generating triple of
reflections). The top sides are usually also fixed by a reflection
(see Proposition~\ref{prop:center}), but not always; in some rare
cases, the ``reflection'' actually degenerates to the identity.

For convenience, we reproduce some of the information in the tables
from the appendix of~\cite{dpp2}, see Tables~\ref{tab:Sc}
and~\ref{tab:Tc}. For each family of groups in the list, we give a
description of side representatives, and mention values of the order
$p$ of the generating reflections where the corresponding side has top
truncation.

\begin{table}
\centering
  {\Large $\Sc(p,\sigma_1)$ groups}\\
  \begin{tabular}{|c|c|c|}
    \hline
    Triangle & \#($P$-orb) & Top trunc.\\
    \hline
    $[6]\ 1;\ 2,\ 3$ &  8  & $p=4,6$  \\                   
    $[4]\ 2;\ 1,\ 23\bar2$ &  8  & $p=6$\\
    $[3]\ 23\bar2;\ 1,\ 232\bar3\bar2$ & 8 &     \\
    $[3]\ 232\bar3\bar2;\ 1,\ \bar3\bar2323$ & 8 &  \\
    \hline
  \end{tabular}
\ \\[0.3cm]
%% s4c
  {\Large $\Sc(p,\bar \sigma_4)$ groups}\\
\begin{tabular}{|c|c|c|}
\hline
                Triangle & \#($P$-orb) & Top trunc.\\
\hline
        $[4]\ 1;\ 2,\ 3$ &     7       & $p=5,6,8,12$\\
  $[3]\ 2;\ 1,\ 23\bar2$ &     7       & $p=8,12$\\
\hline
\end{tabular}
\ \\[0.3cm]
{\Large $\Sc(p,\sigma_5)$ groups}\\
\begin{tabular}{|c|c|c|}
\hline
  Triangle                                                   & \#($P$-orb) & Top trunc.\\
\hline
  $[4]\ 1;\ 2,\ 3$                                           &    30       &        \\              
  $[5]\ 2;\ 1,\ 23\bar2$                                     &    30       & $p=4$  \\
  $[6]\ 2\bar3\bar2123\bar2;\ 2,\ \bar3\bar2\bar123\bar2123$ &     5       & $p=4$  \\
\hline
\end{tabular}
% s10
\ \\[0.3cm]
{\Large $\Sc(p,\sigma_{10})$ groups}\\
\begin{tabular}{|c|c|c|}
\hline
  Triangle & \#($P$-orb) & Top trunc.                                      \\
  $[5]\ 1;\ 2,\ 3$ &  5  & $p=4,5,10$                                     \\              
  $[3]\ 2;\ 1,\ 23\bar2$ &  5 & $p=10$                                    \\
\hline
\end{tabular}
\ \\[0.3cm]
\caption{Side representatives for $\Sc(p,\tau)$ groups}\label{tab:Sc}
\end{table}
\begin{table}
% S2
\ \\
{\Large $\Tc(p,{\bf S_2})$ groups}\\
\begin{tabular}{|c|c|c|}
\hline
  Triangle             & \#($P$-orb) & Top trunc.            \\
\hline
$[3]\ 1;\ 2,\ 3$       &  5  &                               \\                   
$[3]\ 23\bar2;\ 1,\ 3$ &  5  &                               \\                   
$[4]\ 3;\ 1,\ 2$       &  5  & $p=5$                         \\
$[5]\ 2;\ 1,\ 23\bar2$ &  5  & $p=4,5$                       \\
\hline
\end{tabular}
% E2
\ \\[0.3cm]
{\Large $\Tc(p,{\bf E_2})$ groups}\\
\begin{tabular}{|c|c|c|}
\hline
  Triangle & \#($P$-orb) & Top trunc.                                    \\
\hline
$[3]\ 1;\ 2,\ 3$ &  6  &                                          \\                   
$[4]\ 23\bar2;\ 1,\ 3$ &  6  & $p=6$                                  \\                   
$[4]\ 3;\ 1,\ 2$ &  6  & $p=6$                                        \\
$[4]\ 2;\ 1,\ 23\bar2$ &  6  & $p=6$                                  \\
$[6]\ \bar313;\ 12\bar1,\ 3$ &  3  & $p=4,6$                          \\
\hline
\end{tabular}
% H1
\ \\[0.3cm]
{\Large $\Tc(2,{\bf H_1})$ groups}\\
\begin{tabular}{|c|c|}
\hline
  Triangle & \#($P$-orb) \\
\hline
$[3]\ 1;\ 2,\ 3$ &  42 \\                   
$[3]\ 23\bar2;\ 1,\ 3$ &  42  \\                   
$[4]\ 3;\ 1,\ 2$ &  42  \\
$[7]\ 2;\ 1,\ 23\bar2$ &  42  \\
$[14]\ \bar2123\bar2\bar12;\ \bar212,\ \bar312\bar13$ &  3  \\
\hline
\end{tabular}
% H2
\ \\[0.3cm]
{\Large $\Tc(p,{\bf H_2})$ groups}\\
\begin{tabular}{|c|c|c|}
\hline
  Triangle                                                    & \#($P$-orb) & Top trunc.\\
\hline
$[3]\ 1;\ 2,\ 3$                                              & 15 &                                           \\                   
$[3]\ 23\bar2;\ 1,\ 3$                                        & 15 &                                     \\                   
$[5]\ 2;\ 1,\ 23\bar2$                                        & 15 & $p=5$                                  \\                   
$[5]\ 3;\ 1,\ 2$                                              & 15 & $p=5$                                        \\                   
$[10]\ (123)^2\bar212(\bar3\bar2\bar1)^2;\ 121\bar2\bar1,\ 3$ & 3 & $p=3,5$ \\                   
\hline
\end{tabular}
\ \\[0.3cm]
\caption{Side representatives for $\Tc(p,{\bf T})$ groups}\label{tab:Tc}
\end{table}

\section{Orbits of complex 2-facets}

In this section, we discuss how to list the conjugacy classes of
mirrors of complex reflections in all triangle groups with parameters
in Table~\ref{tab:list}.

The first observation is that the list can be read off any
``reasonable'' fundamental domain; indeed let $\Gamma$ be a complex
hyperbolic lattice, and let $\Pi$ be a closed, piecewise smooth,
finite-sided fundamental domain for $\Gamma$, in particular
\begin{enumerate}
  \item $\CH^2=\cup_{\gamma\in\Gamma}\gamma \Pi$ and
  \item for every $\gamma\in\Gamma$,
    $\gamma(\Pi^\circ)\cap\Pi^\circ = \emptyset$.
\end{enumerate}
where $\Pi^\circ$ denotes the interior of $\Pi$. We assume moreover
that $\Pi$ has a side pairing in the sense of the Poincar\'e
polyhedron theorem (see section~4.3 of~\cite{dpp2}, or section~3.2
of~\cite{dpp1}).

Now assume $R\in\Gamma$ is a complex reflection with mirror $m$. Pick
$x\in m$, and let $\gamma\in\Gamma$ be such that $\gamma(x)\in\Pi$
(see condition~(1)). By replacing $R$ by $\gamma R\gamma^{-1}$, we may
assume $m\cap\Pi\neq \emptyset$. Because of~(2), we must have
$m\cap \pi\subset \partial\Pi$. In particular, we have the following.
\begin{prop}\label{prop:refl-basic}
  Every reflection in $\Gamma$ is conjugate to one fixing a complex
  2-facet of $\Pi$.
\end{prop}

Moreover, we can easily check whether reflections fixing two different
complex 2-facets of $\Pi$ are conjugate in $\Gamma$, by tracking
cycles in the sense of the Poincar\'e; moreover, the Poincar\'e
polyhedron theorem allows us to determine the full stabilizer in
$\Gamma$ of every facet (in particular complex 2-facets). This is part
of what the computer program available at~\cite{spocheck} does
systematically, and the results are summarized in the computer output
files available with the source code.

A slight subtlety is that the domains we construct in~\cite{dpp2} are
not fundamental domains, but only fundamental domains for coset
decompositions, with respect to a finite cyclic group generated by a
regular elliptic element ($P=R_1J$ for $\Sc(p,\tau)$ groups,
$Q=R_1R_2R_3$ for $\Tc(p,{\bf T})$ groups). In the next few
paragraphs, we only discuss the case of $\Sc$ groups, everything goes
through for $\Tc$ groups, as long as we replace $P$ by $Q$).

The domains used in~\cite{dpp2} are of the form
$$
  \cup_{j=0}^n P^j(\Pi),
$$
where $\Pi$ is a true fundamental domain, $n$ is the order of $P$. For
future reference, we write $\Pi^{(P)}=\cup_{j=0}^n P^j(\Pi)$.

When all powers of $P$ are regular elliptic, this is inconsequential,
but when some power is a complex reflection, we need to add extra
conjugacy classes of complex reflections.
\begin{prop}\label{prop:refl-coset}
  Every reflection in a lattice $\Gamma=\Sc(p,\tau)$ is conjugate to
  one fixing a complex 2-facet of $\Pi^{(P)}$ or to one fixing the
  mirror of some complex reflection power of $P$ (and similarly for
  the groups $\Tc(p,{\bf T})$ with $P$ replaced by $Q$).
\end{prop}

It follows from the previous paragraphs (and the discussion of
section~\ref{sec:sides}) that the mirror of every complex reflection
in the group can be mapped via a suitable group element to one in the
following three families.
\begin{itemize}
\item The complex lines containing base polygons of various sides (all
  of these are mirrors of conjugates of the three generating
  reflections);
\item The complex lines containing the top 2-facets of the sides with
  truncated apex;
\item Mirrors of complex reflections obtained as
  powers of $P=R_1J$ (for groups in the $\Sc$-family) or $Q=R_1R_2R_3$
  (for groups in the $\Tc$-family) that are complex reflections.
\end{itemize}
We call the corresponding complex 2-facets in these three families
base 2-facets, top 2-facets and central 2-facets, respectively. The
list of base/top/central facets can easily be deduced from the tables
and pictures in the appendix of~\cite{dpp2}.

There are some obvious identifications between sides, coming from
\begin{itemize}
\item The side-pairing maps; these are all complex reflections fixing
  the base side.
\item The action of $P=R_1J$ (for $\Sc$ groups) or $Q=R_1R_2R_3$ (for
  $\Tc$ groups).
\end{itemize}
Accordingly, the lists of the appendix in~\cite{dpp2} only list one
representative for each $P$-orbit (resp. $Q$-orbit) of paired sides.

\begin{dfn}
  For future reference, we refer to the sides that are listed for a
  given $\Gamma=\Tc(p,{\bf T})$ in the appendix of~\cite{dpp2} as
  \emph{side representatives} for $\Gamma$.
\end{dfn}

Moreover, when two sides $s_1$ and $s_2$ appear in the appendix
of~\cite{dpp2}, the only identifications between $s_1$ and $s_2$ can
occur on their boundary (here by the boundary of a side $s_j$
contained in a bisector $\Bc$, mean boundary with respect to the
topology on $\Bc$ induced by the topology on $\CH^2$). The
identifications on the boundary can be read off by tracking orbits of
all side-pairing transformations.

The details are difficult to write down on paper because there are
many groups and each polytope has many facets, so we only summarize
the important points in Proposition~\ref{prop:ridge-injection} (which
follows by case by case analysis, tracking the cycles of 2-facets
systematically for all fundamental domains).
\begin{prop}\label{prop:ridge-injection}
  \begin{enumerate}
  \item Let $r_1$ and $r_2$ be distinct base or top complex 2-facets
    of (possibly different) side representatives, and let $x_j$ be in
    the open cell given by the interior of $r_j$. Then $x_1$ and $x_2$
    are not in the same $\Gamma$-orbit.
  \item Let $r$ be a base or top 2-facet of a side representative $S$.
    Assume $S$ is not invariant under any non-trivial power $P^j$ (for
    $\Sc$ groups) or $Q^j$ (for $\Tc$ groups), and $x_1,x_2$ are
    distinct points in the open cell given by the interior of
    $r$. Then $x_1$ and $x_2$ are not in the same $\Gamma$-orbit.
  \item Let $r$ be a base or top 2-facet of a side representative $S$
    contained in a bisector $\Bc$. Assume $S$ is invariant under
    $R=P^j$ or $Q^j$ (and let $j>0$ be minimal with that property);
    then $R$ is a complex reflection whose mirror is the complex spine
    of $\Bc$. Points in $r$ are $\Gamma$-equivalent if and only if they
    are invariant under the action of the cyclic group generated by $R$.
  \end{enumerate} 
\end{prop}

We think of these statements as giving injectivity for the quotient
map on various (unions of) complex 2-facets. Point~(2) says that the
interior of base and top 2-facets inject, unless the side is invariant
under a power of $P$ (or $Q$). Point~(3) says that if the side is
invariant under a power or $P$ (or $Q$), then the side is rotationally
symmetric, and one gets a sector (of angle given by the rotation angle
of $P^j$) that injects under the quotient map.

\begin{figure}[htbp]
  \centering
  \includegraphics{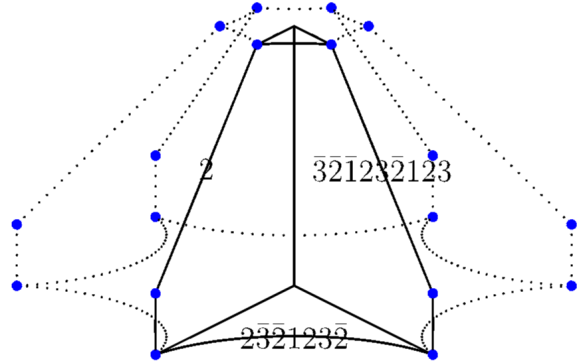}
  \caption{One of the sides for $\Sc(p,\sigma_5)$ is invariant under
    $P^5$, with rotational symmetry of angle $\pi/3$.}
\label{fig:s5}
\end{figure}

The exceptions to the injectivity statements as in~(3) do occur for
some of the groups we handle in~\cite{dpp2}, they can be spotted by
listing the sides whose $P$-orbit contains less than $\textrm{Ord}(P)$
elements (see the column headed $\#(P-orb)$) for $\Sc$-groups, and
similary with $P$ replaced by $Q$ for $\Tc$ groups (see
Tables~\ref{tab:Sc} and~\ref{tab:Tc}). For example, the side
$2\bar3\bar2123\bar2;2,Q^{-1}23\bar2Q$ is depicted in
Figure~\ref{fig:s5}; a fundamental domain for the action of the
rotation is a sector with angle $\pi/3$.

\section{Tiling mirrors with images of complex 2-facets} \label{sec:tiling}

In this section, we study the mirror stabilizers for a representative
of every orbit of a complex 2-facet of $\Pi$, and also of the mirror
of complex reflection powers of $P$ or $Q$ (if
any). Proposition~\ref{prop:refl-coset} says that this will give us
the list of all mirrors of complex reflections up to conjugation.

The basic observation is the following.
\begin{thm}\label{thm:fuchsian}
  If $m$ is the mirror of any complex reflection in $R$ in $\Gamma$,
  then $Stab_\Gamma(m)$ is a lattice in $Stab_{PU(H)}(m)$.
\end{thm}
\begin{pf}
  Every point on $m$ can be mapped to $\partial \Pi$ by a suitable
  element of $\Gamma$, so $m$ is tiled by $\Gamma$-images of complex
  2-facets in $\partial\Pi$. If $s$ is a complex 2-facet of $\pi$ and
  $\gamma_1(s)$ and $\gamma_2(s)$ are both contained in $m$, then
  $\gamma_1\gamma_2^{-1}\in Stab_\Gamma(m)$; $\Pi$ has finitely many
  complex 2-facets, and each has a finite area, so
  $Stab_\Gamma(m)\backslash m$ has finite area.
\end{pf}
Our goal is to make the corresponding lattices in $PU(1,1)$
explicit. The proof of the theorem suggests a method to make the
stabilizer explicit; we will list complex 2-facets of $\Pi$ that are
in the same group orbit, map them to lie in the same complex line $m$,
in such a way that their images in $m$ glue nicely to give a
fundamental domain for a lattice in $PU(1,1)$, with a well defined
side-pairing.

Note that the difficulty is to find an explicit side-pairing, the fact
that the images give a (possibly disconnected) fundamental domain for
the action of the stabilizer is automatic (see
Proposition~\ref{prop:ridge-injection} and the proof of
Theorem~\ref{thm:fuchsian}).

We now give a list of explicit maps that bring the complex 2-facets
back to a common complex line. For top and central complex 2-facets,
the stabilizers are all triangle groups, and it is very easy to obtain
a fundamental domain consisting of two adjacent triangles.

We treat the mirrors of the generating set of reflections $R_1$,
$R_2$, $R_3$ by a case by case analysis for all families of groups in
the list of Table~\ref{tab:list} (the maps are given by the same
formula for various values of $p$ within a fixed family
$\Tc(p,{\bf T})$ of groups, even though the facets themselves possibly
depend on $p$ via truncation).

For $\Sc(p,\tau)$ groups, we have $JR_jJ^{-1}=R_{j+1}$, so it is
enough to handle the mirror of $R_1$. For most groups
$\Tc(p,{\bf T})$, the mirrors of $R_1$, $R_2$ and $R_3$ are still
images of each other under a suitable group elements, because of the
fact that $\br_n(a,b)$ with $n$ odd implies that $a$ and $b$ are
conjugate, because
$$
   b = (ab)^{-\frac{n-1}{2}}a(ab)^{\frac{n-1}{2}},
$$
see also Proposition~\ref{prop:center}(2).

Indeed, the only family in the list of Table~\ref{tab:list} where all
braid lengths $\br(R_j,R_k)$ with $j\neq k$ are even is the family of
$\Tc(p,{\bf E_2})$ groups (see also Table~\ref{tab:Tc} for the braid
lengths between generating reflections). For that family, we will
handle both the mirrors of $R_1$ and $R_2$.

\begin{figure}
  \centering
  \includegraphics[width=0.7\textwidth]{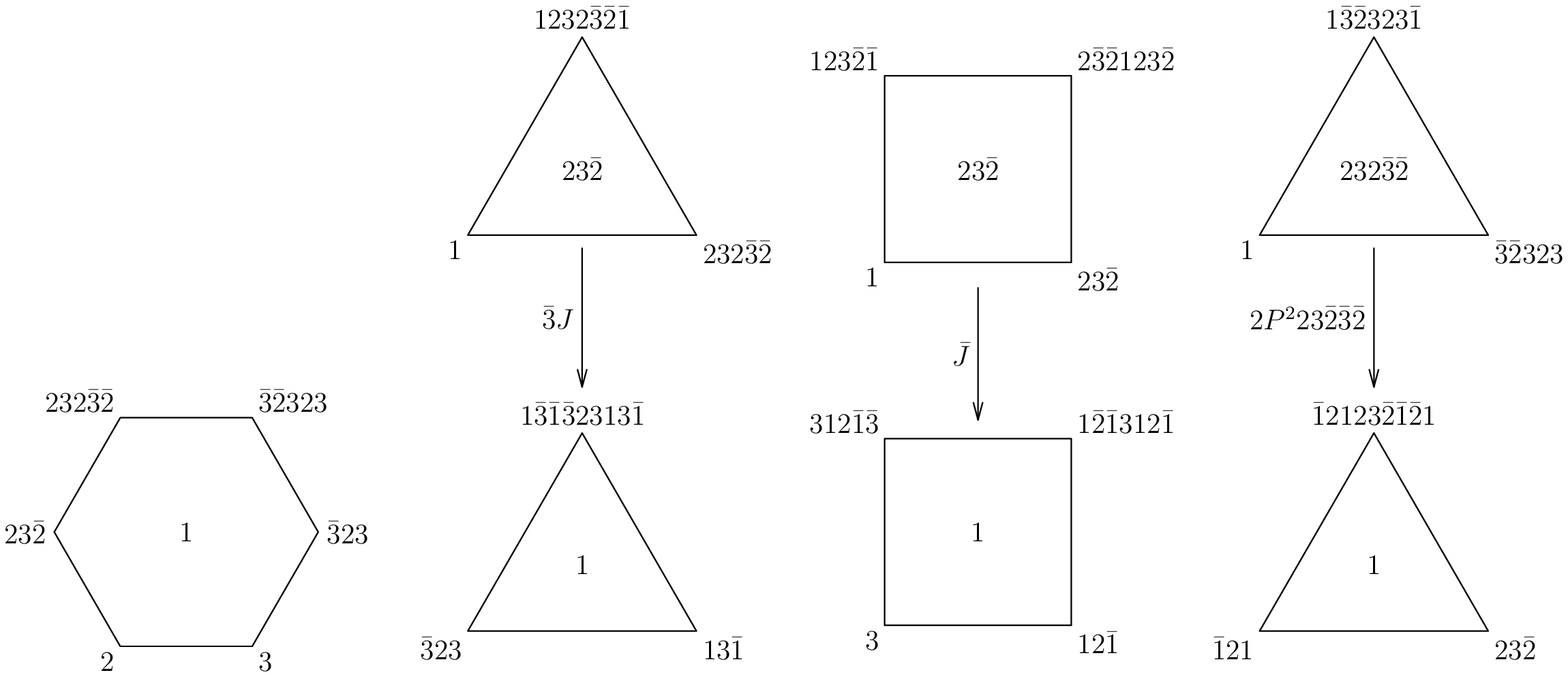}
  \caption{Maps to the mirror of $R_1$ for
    groups $\Sc(p,\sigma_1)$}\label{fig:maps-s1}
\end{figure}

\begin{figure}
  \centering
  \begin{tabular}{c}
    \includegraphics[width=0.25\textwidth]{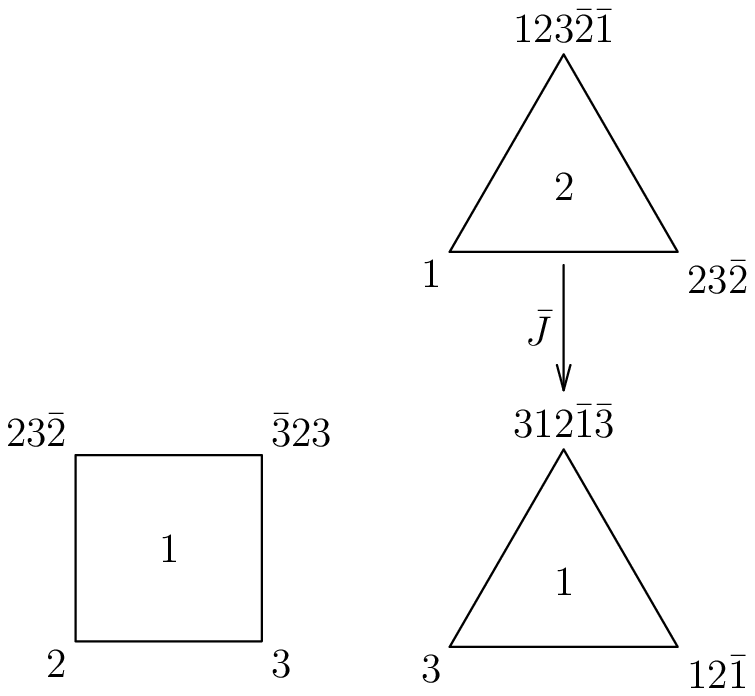}\\
    $\Sc(p,\bar\sigma_4)$
  \end{tabular}
  \hfill
  \begin{tabular}{c}
  \includegraphics[width=0.4\textwidth]{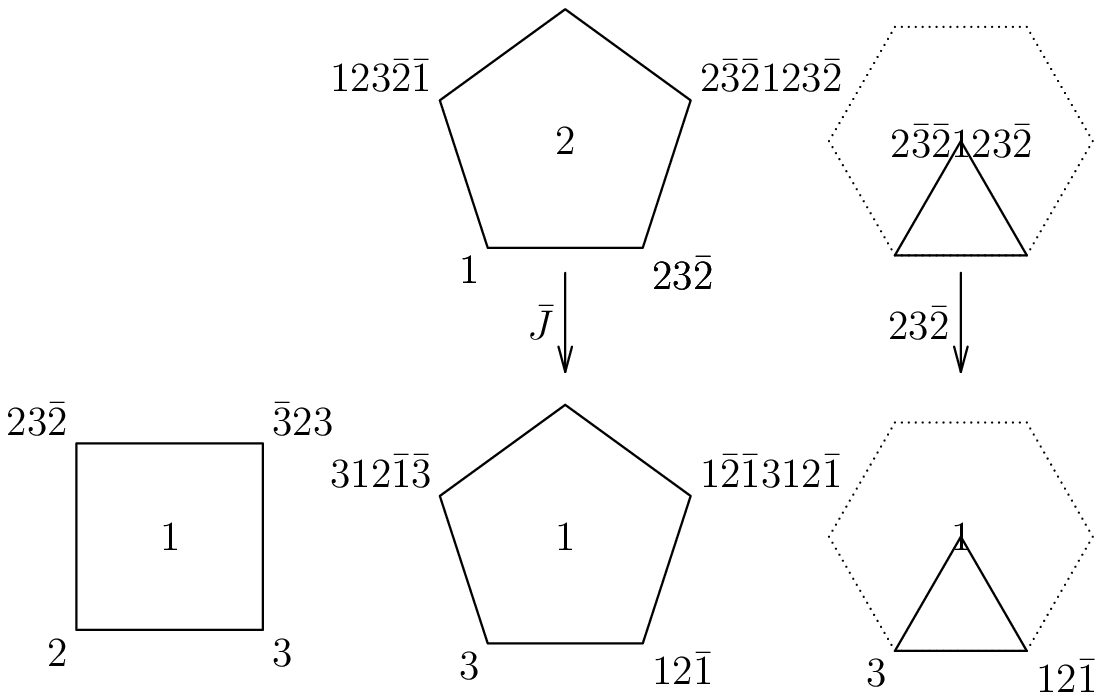}\\
    $\Sc(p,\sigma_5)$
  \end{tabular}
  \hfill
  \begin{tabular}{c}
  \includegraphics[width=0.25\textwidth]{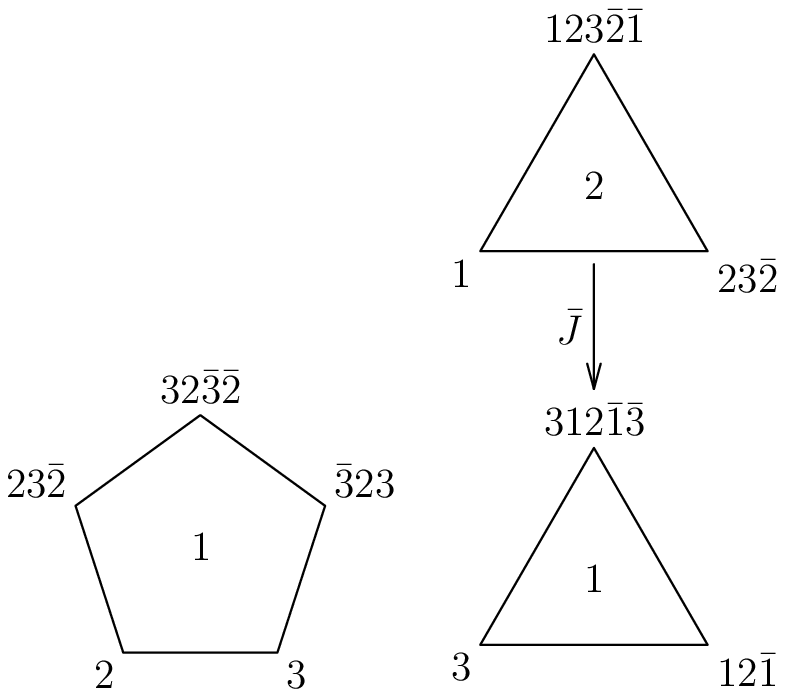}\\
    $\Sc(p,\sigma_{10})$
  \end{tabular}
  \caption{Maps to the mirror of $R_1$ for
    groups $\Sc(p,\bar\sigma_4)$, $\Sc(p,\sigma_5)$ and $\Sc(p,\sigma_{10})$}\label{fig:maps-s4c-s5}
\end{figure}

\begin{figure}
  \centering
  \includegraphics[width=0.7\textwidth]{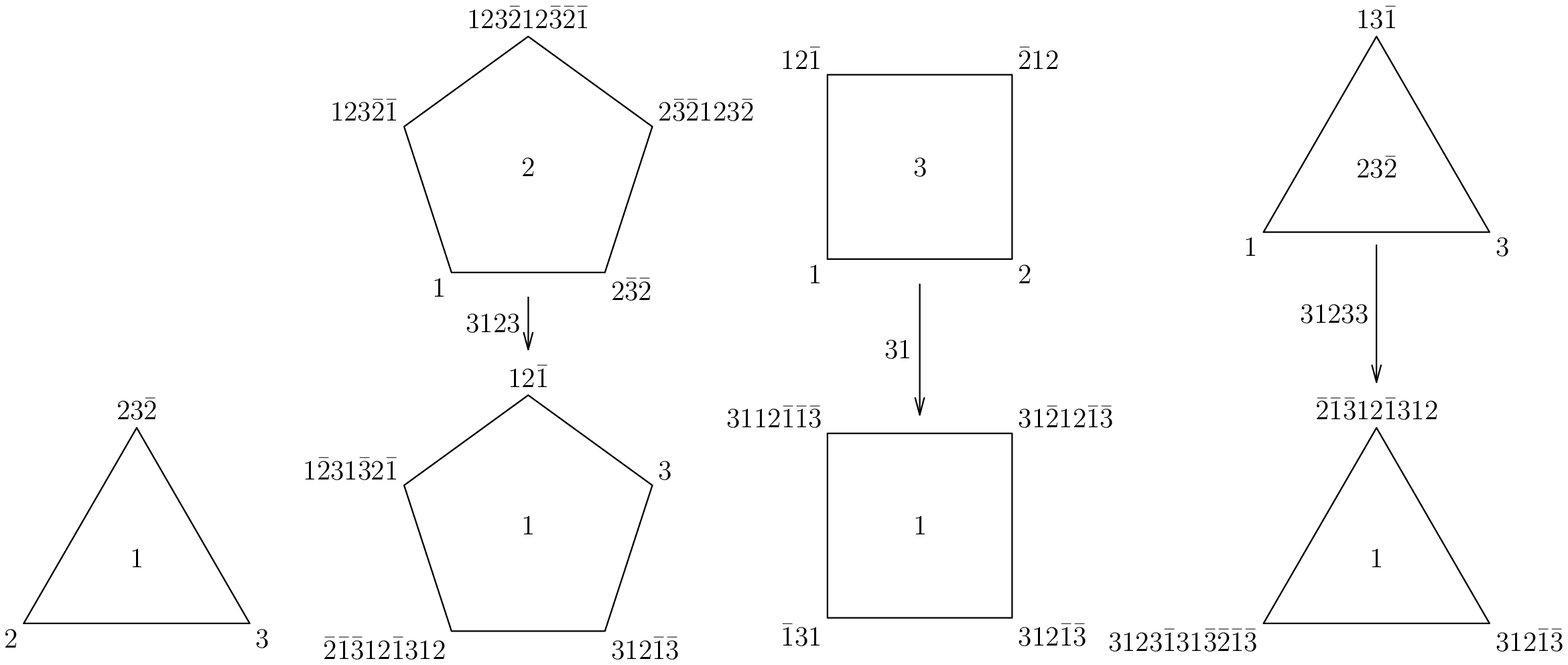}
  \caption{Maps to the mirror of $R_1$ for
    groups $\Tc(p,{\bf S_2})$}\label{fig:maps-S2}
\end{figure}

\begin{figure}
  \centering
  \includegraphics[width=0.7\textwidth]{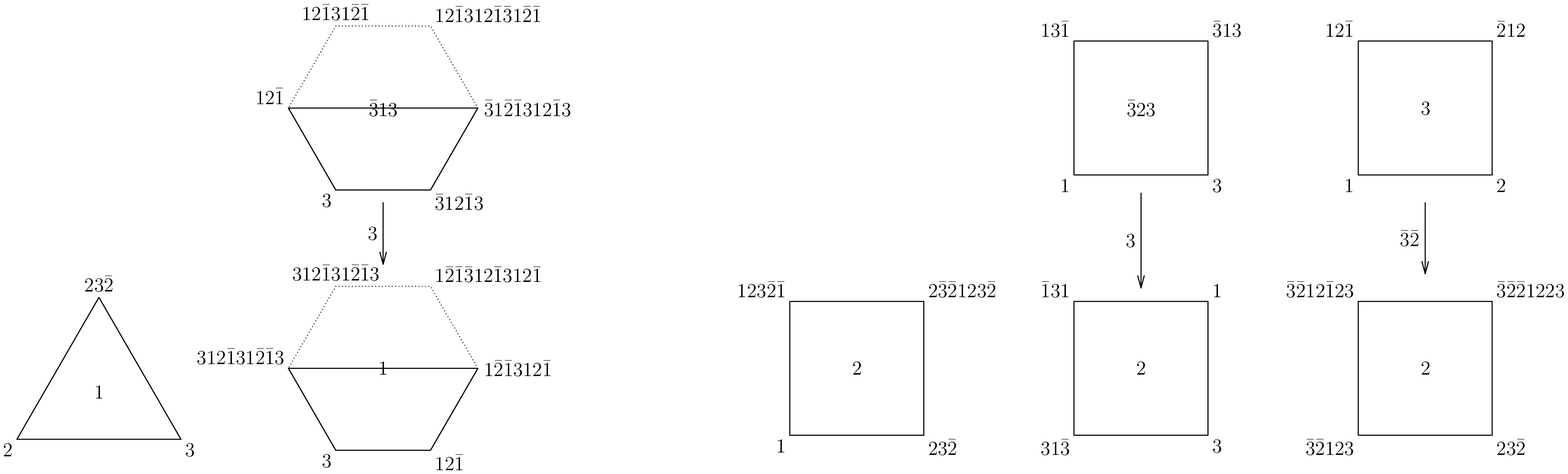}
  \caption{Maps to the mirror of $R_1$ for
    groups $\Tc(p,{\bf S_2})$}\label{fig:maps-E2}
\end{figure}

\begin{figure}
  \centering
  \includegraphics[width=0.7\textwidth]{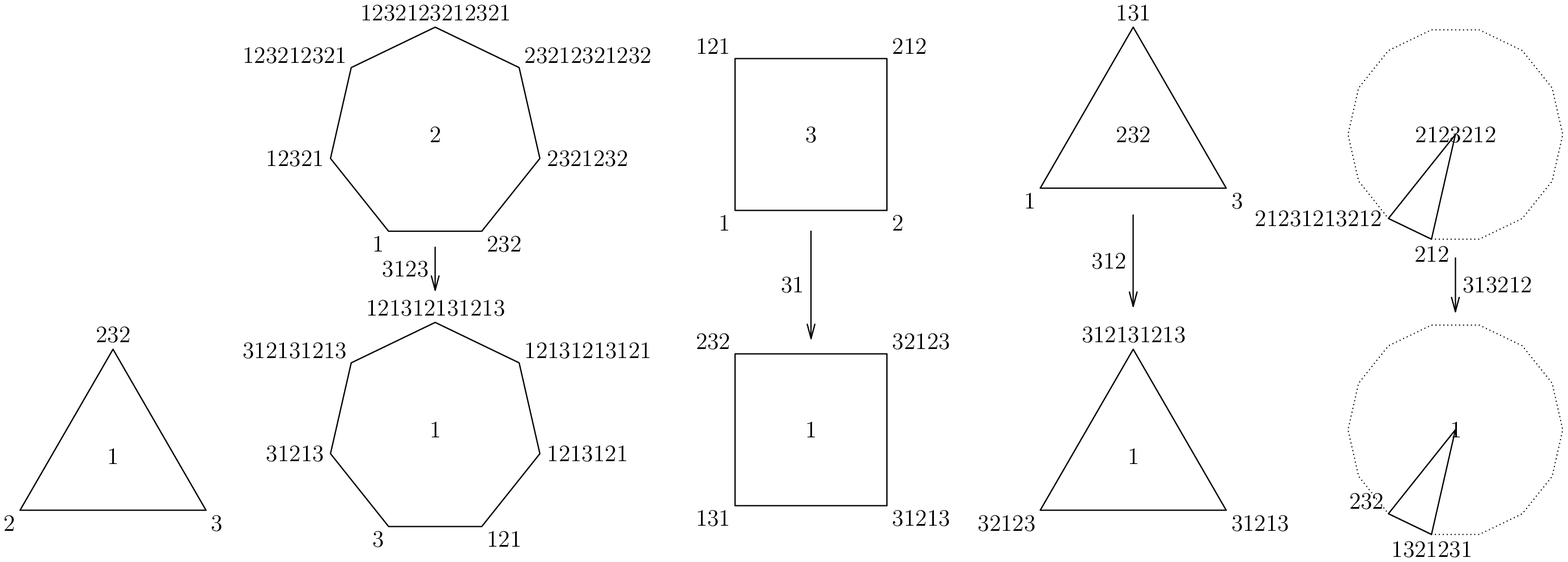}
  \caption{Maps to the mirror of $R_1$ for
    groups $\Tc(p,{\bf H_1})$}\label{fig:maps-H1}
\end{figure}

\begin{figure}
  \centering
  \includegraphics[width=0.7\textwidth]{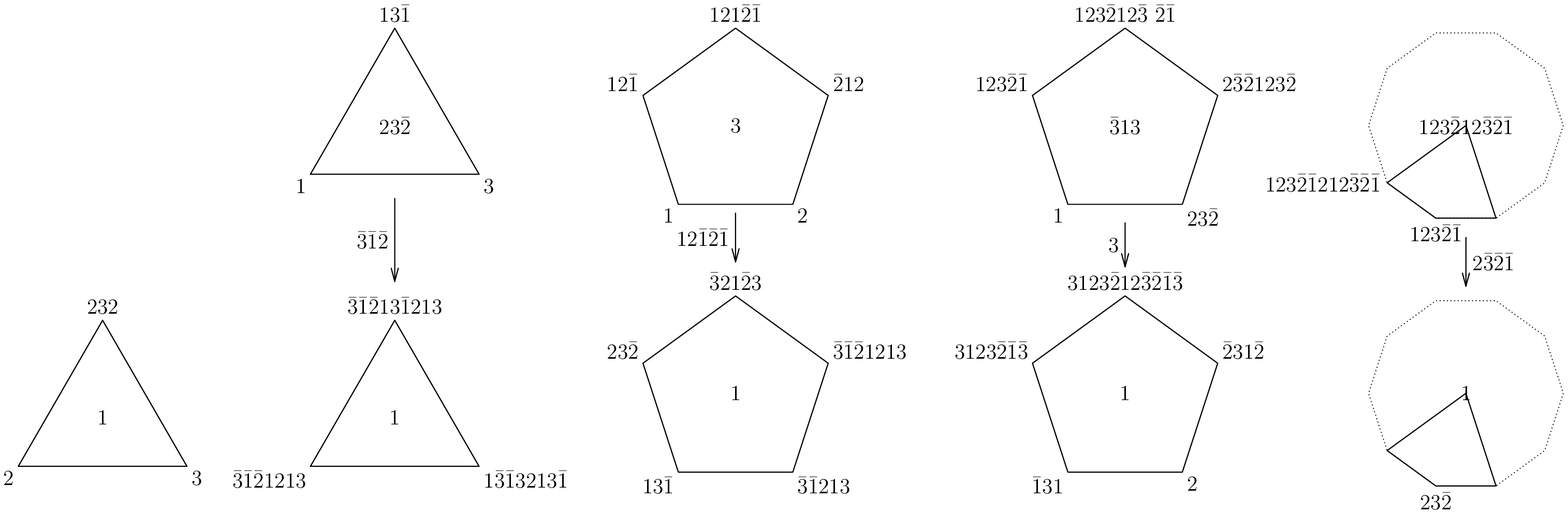}
  \caption{Maps to the mirror of $R_1$ for
    groups $\Tc(p,{\bf H_2})$}\label{fig:maps-H2}
\end{figure}

We list the results in the form of pictures in
Figures~\ref{fig:maps-s1}-\ref{fig:maps-H2}. Note that the pictures
are intended to be only combinatorial, not metric (but
coordinates/parametrizations for the actual vertices/sides can be
recovered from the labels on the pictures).

\section{Side-pairing and fundamental domains}

Our goal is now to check how the puzzle pieces described in
section~\ref{sec:tiling} fit together to form a fundamental domain for
the stabilizer, and to find an explicit side-pairing; from this, one
easily deduces a presentation (via application of the Poincar\'e
polyhedron theorem), as well as the signature of Fuchsian groups.

\begin{figure}[htbp]
  \centering
  \begin{tabular}{c}
    \includegraphics[width=0.3\textwidth]{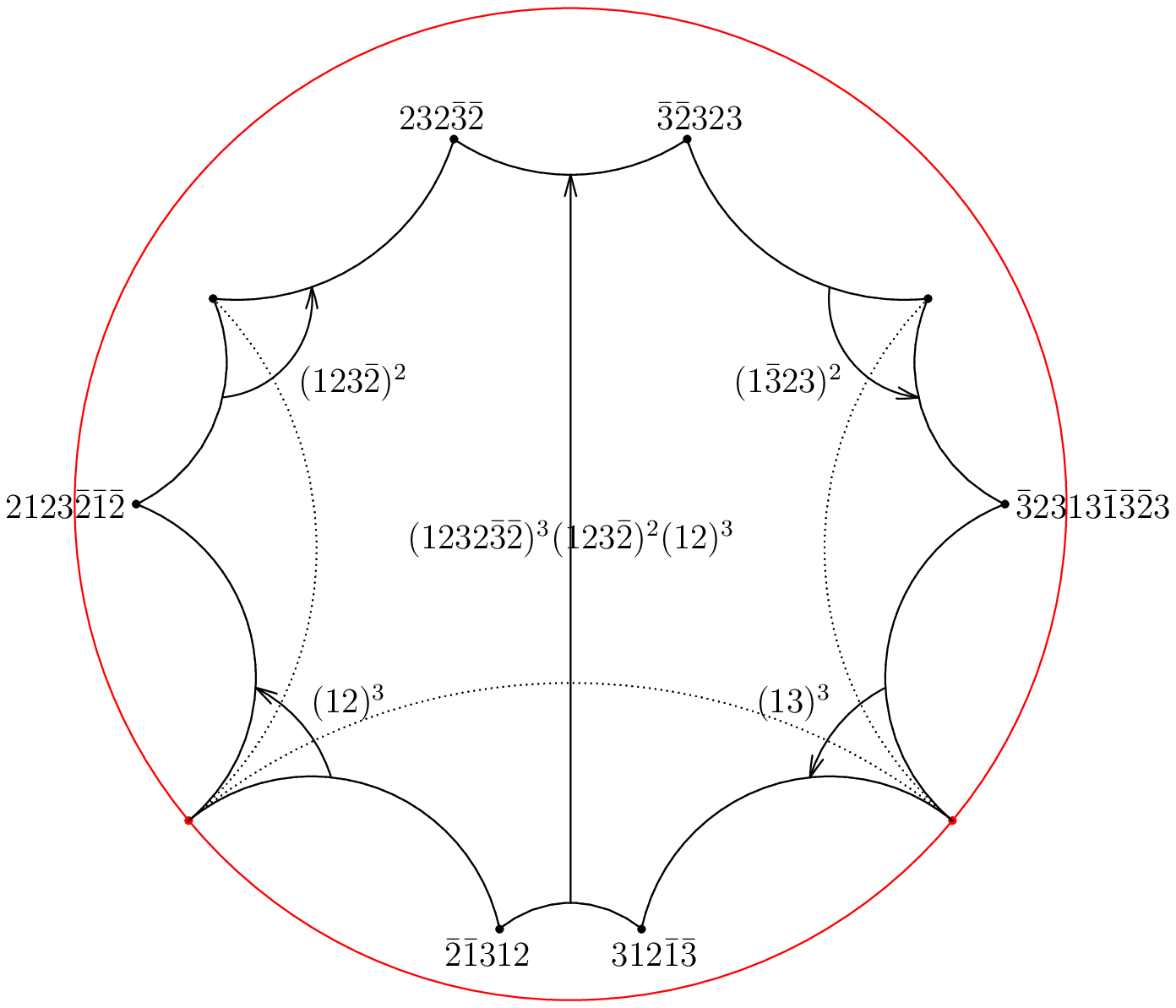}\\
    $p=3$
  \end{tabular}
  \begin{tabular}{c}
    \includegraphics[width=0.3\textwidth]{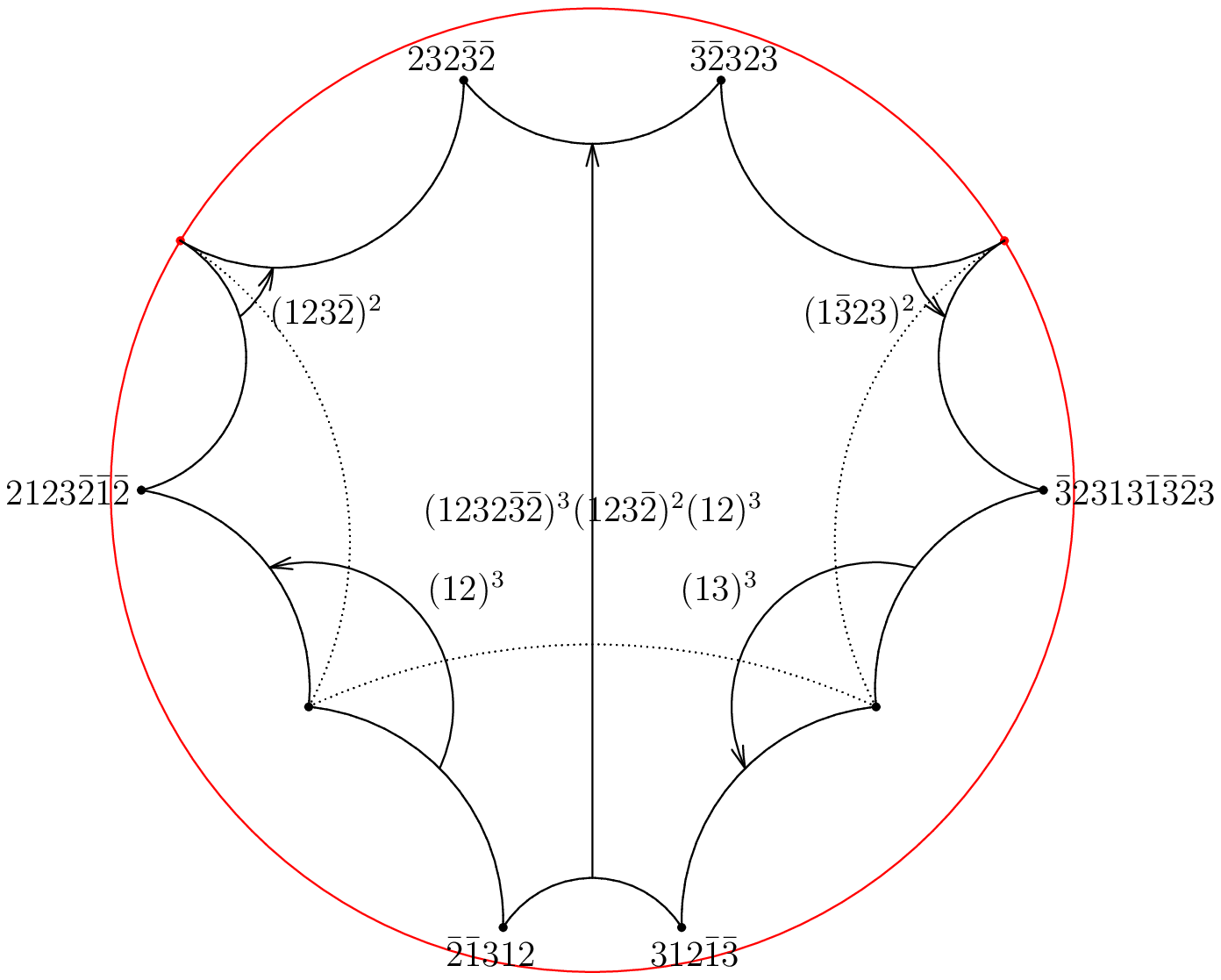}\\
    $p=4$
  \end{tabular}
  \begin{tabular}{c}
    \includegraphics[width=0.3\textwidth]{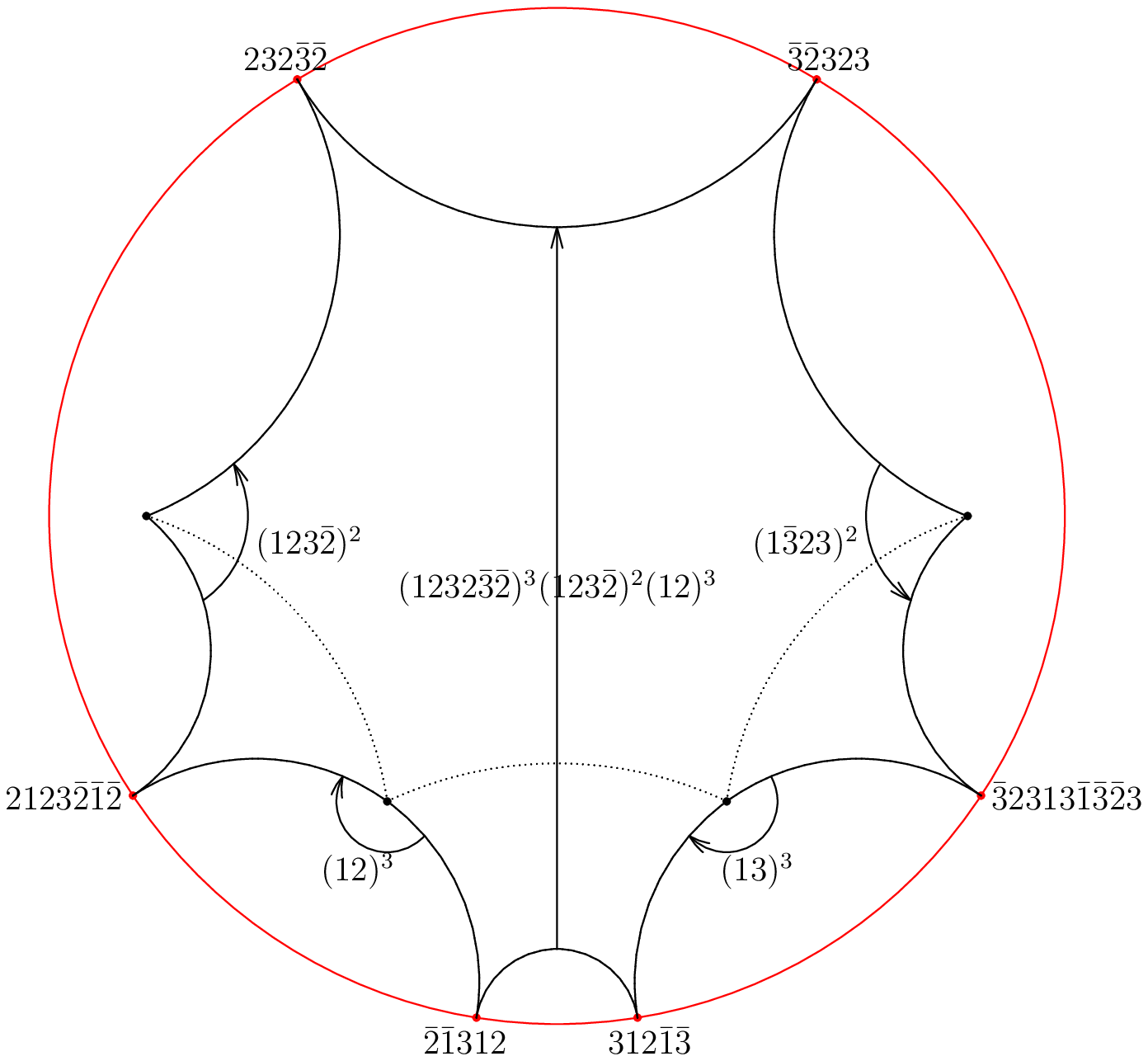}\\
    $p=6$
    \end{tabular}
  \caption{Pairing for the stabilizer of the mirror of $R_1$, for $\Sc(p,\sigma_1)$ groups.}
  \label{fig:stab-s1}
\end{figure}

\begin{figure}[htbp]
  \centering
  \begin{tabular}{c}
    \includegraphics[width=0.3\textwidth]{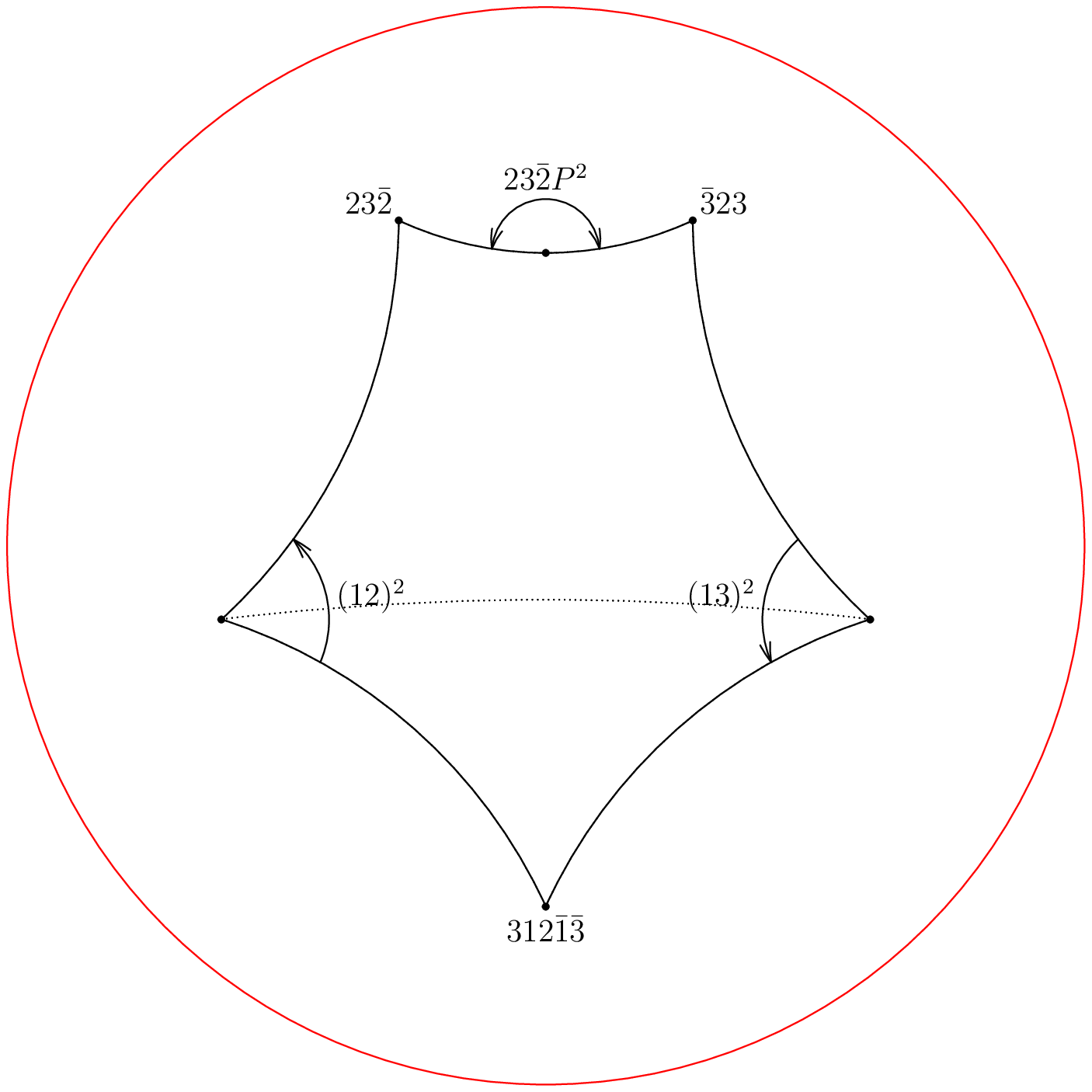}\\
    $p=3$
  \end{tabular}
  \begin{tabular}{c}
    \includegraphics[width=0.3\textwidth]{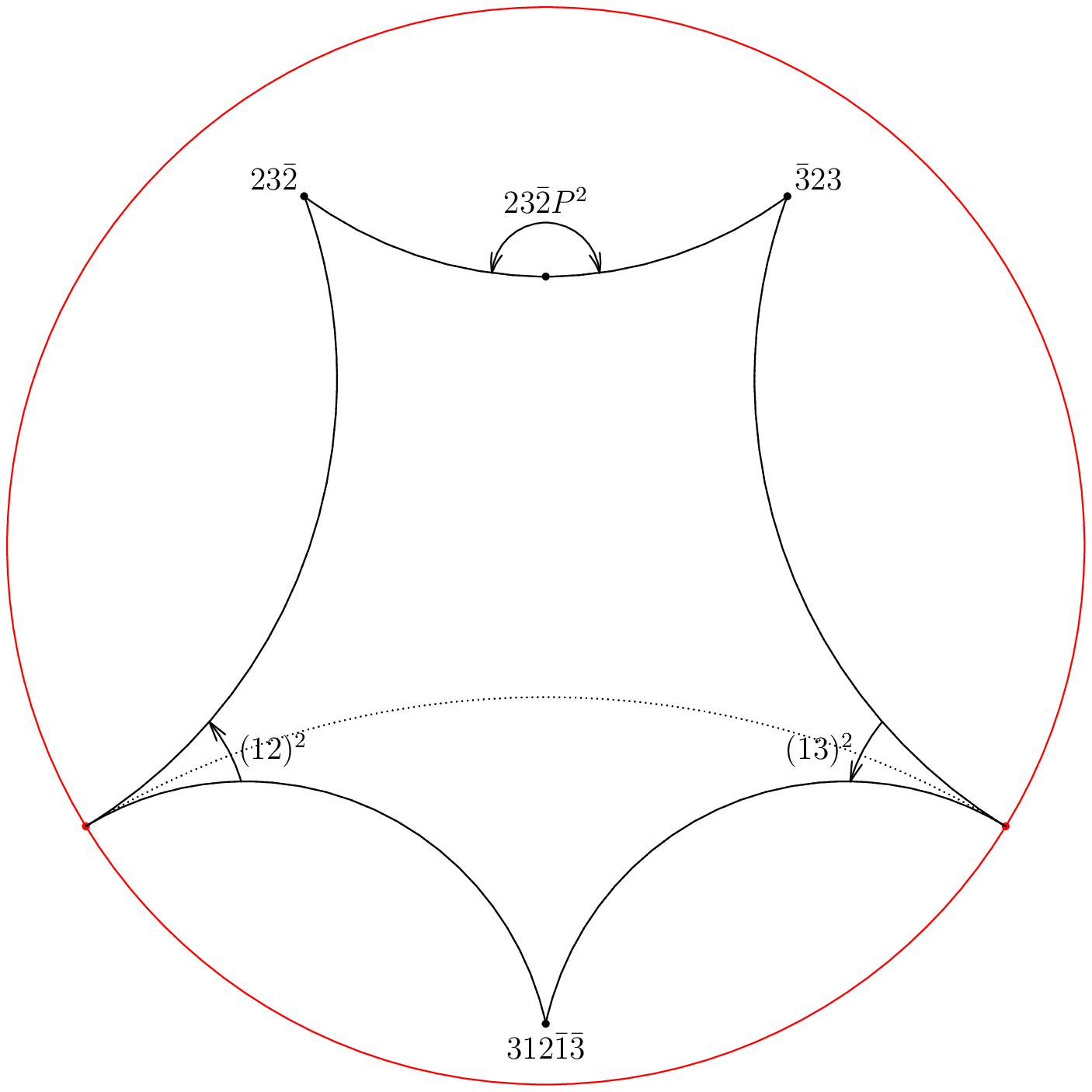}\\
    $p=4$
  \end{tabular}
  \begin{tabular}{c}
    \includegraphics[width=0.3\textwidth]{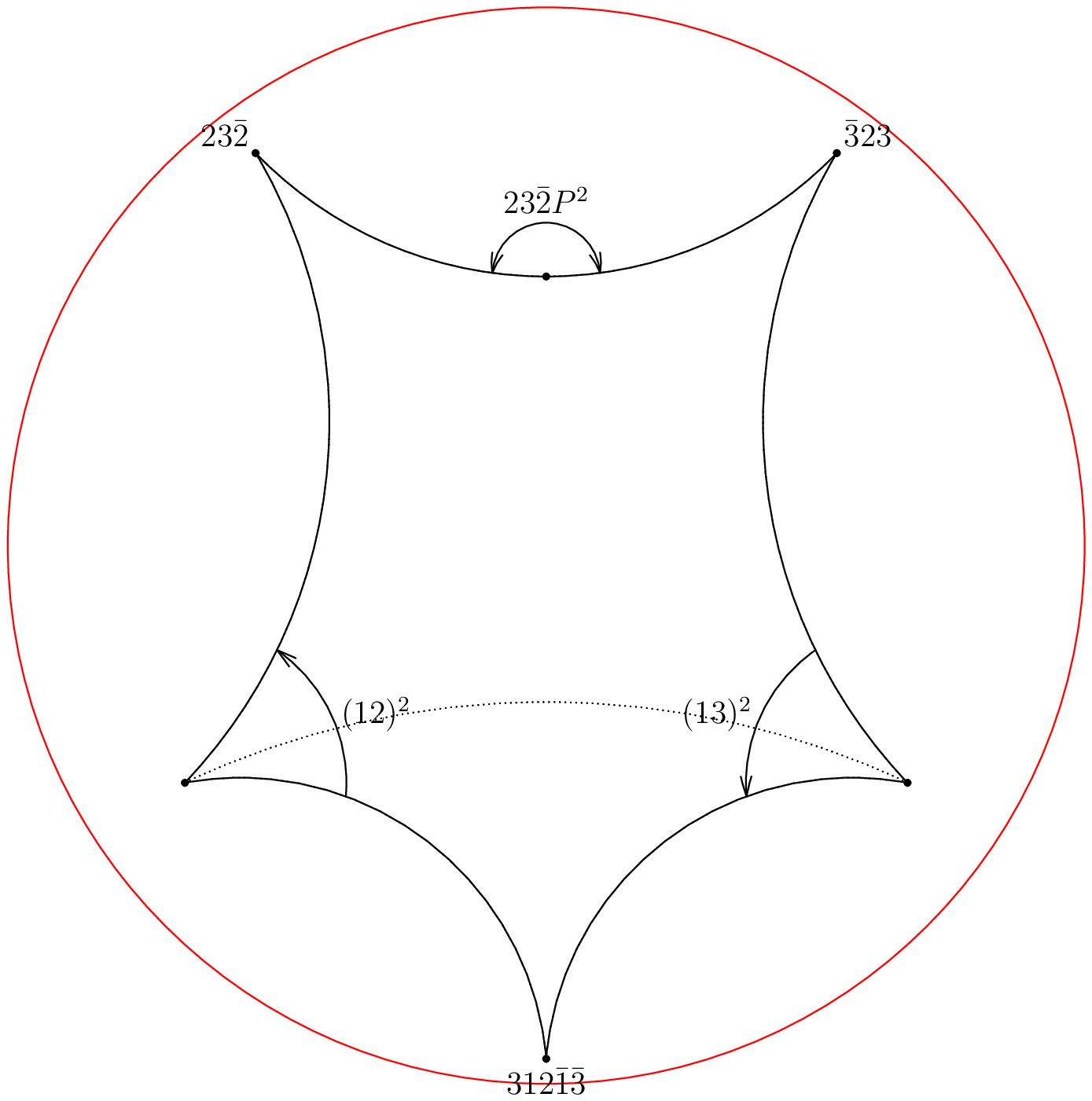}\\
    $p=5$
  \end{tabular}
  \\
  \begin{tabular}{c}
    \includegraphics[width=0.3\textwidth]{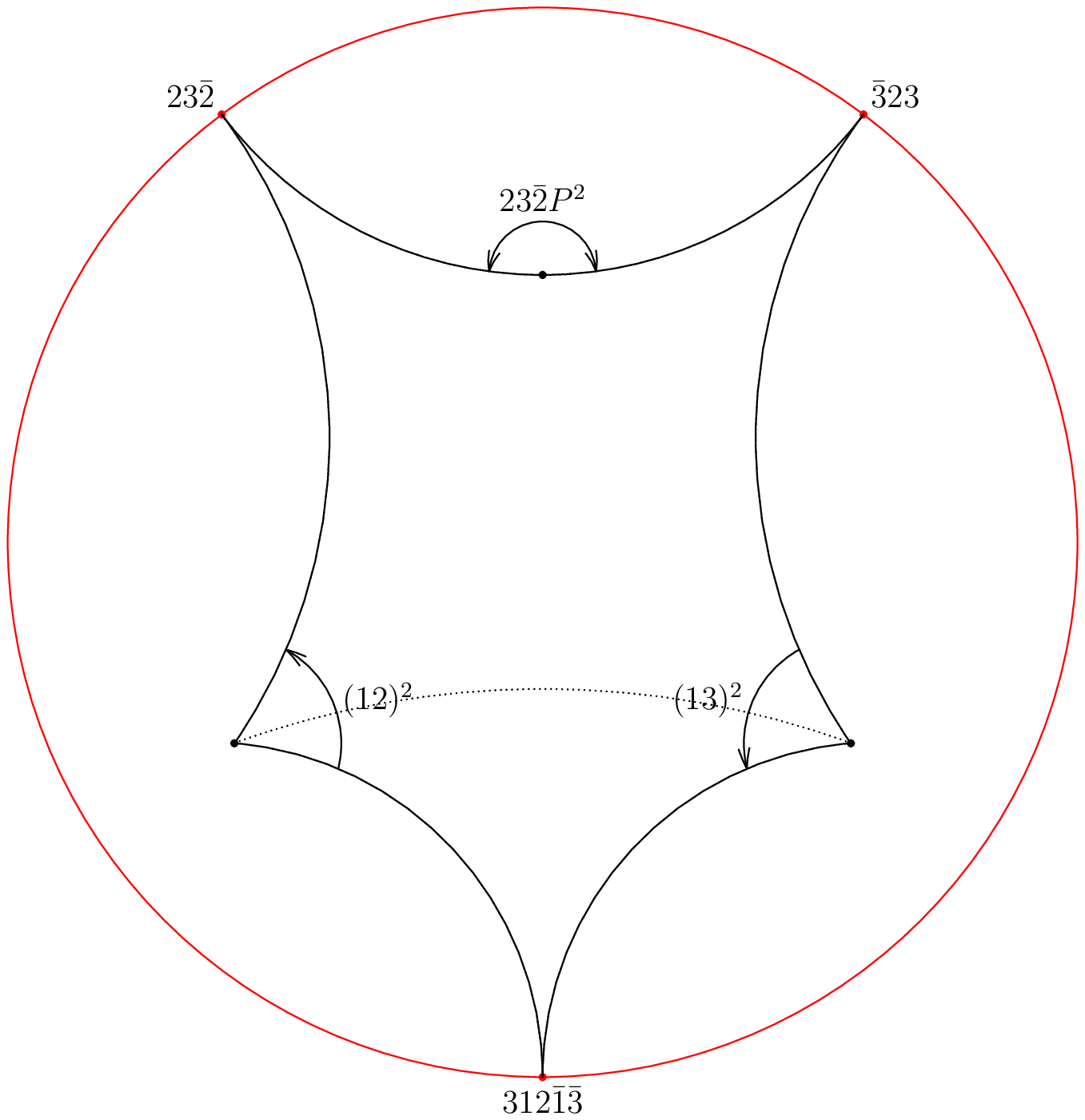}\\
    $p=6$
  \end{tabular}
  \begin{tabular}{c}
    \includegraphics[width=0.3\textwidth]{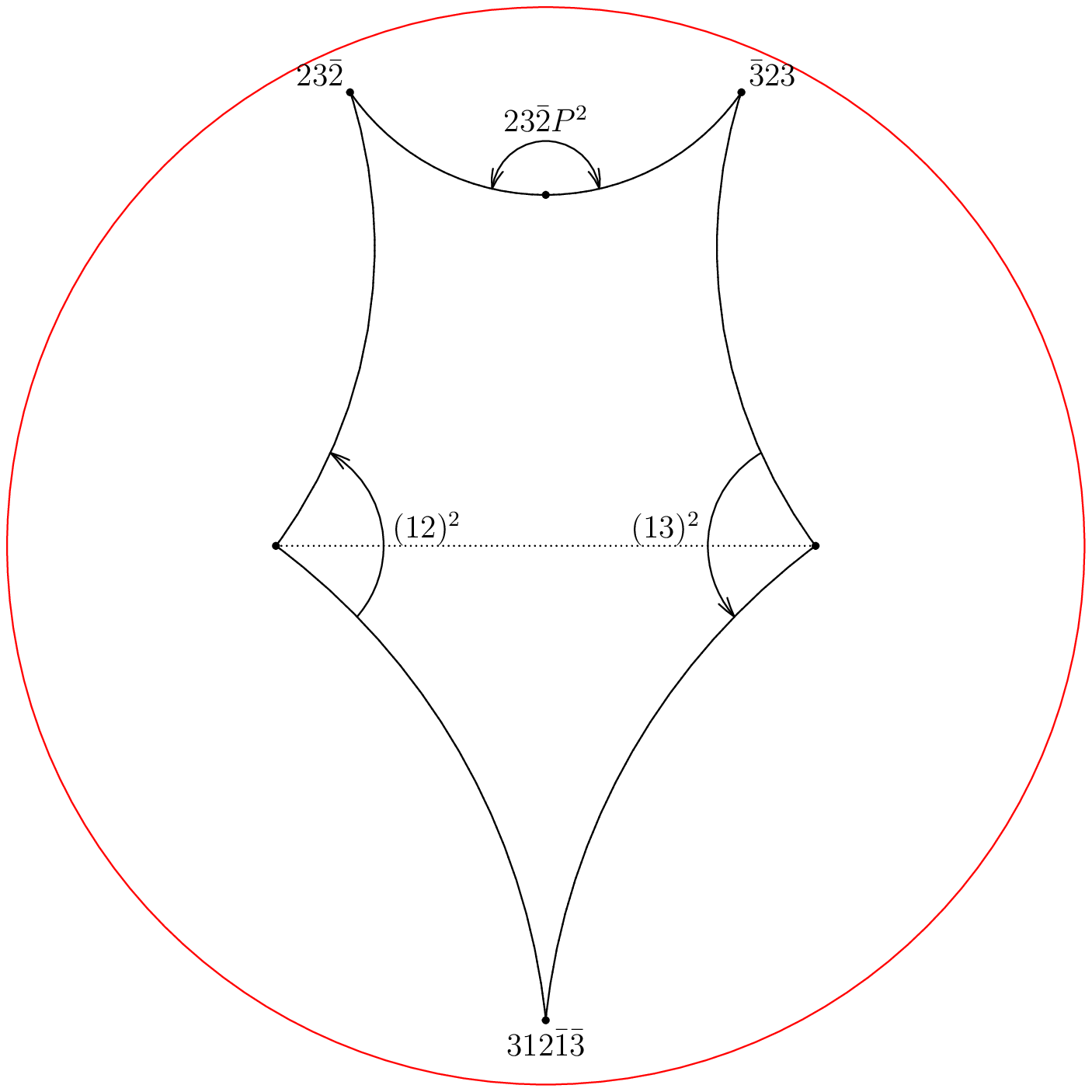}\\
    $p=8$
  \end{tabular}
  \begin{tabular}{c}
    \includegraphics[width=0.3\textwidth]{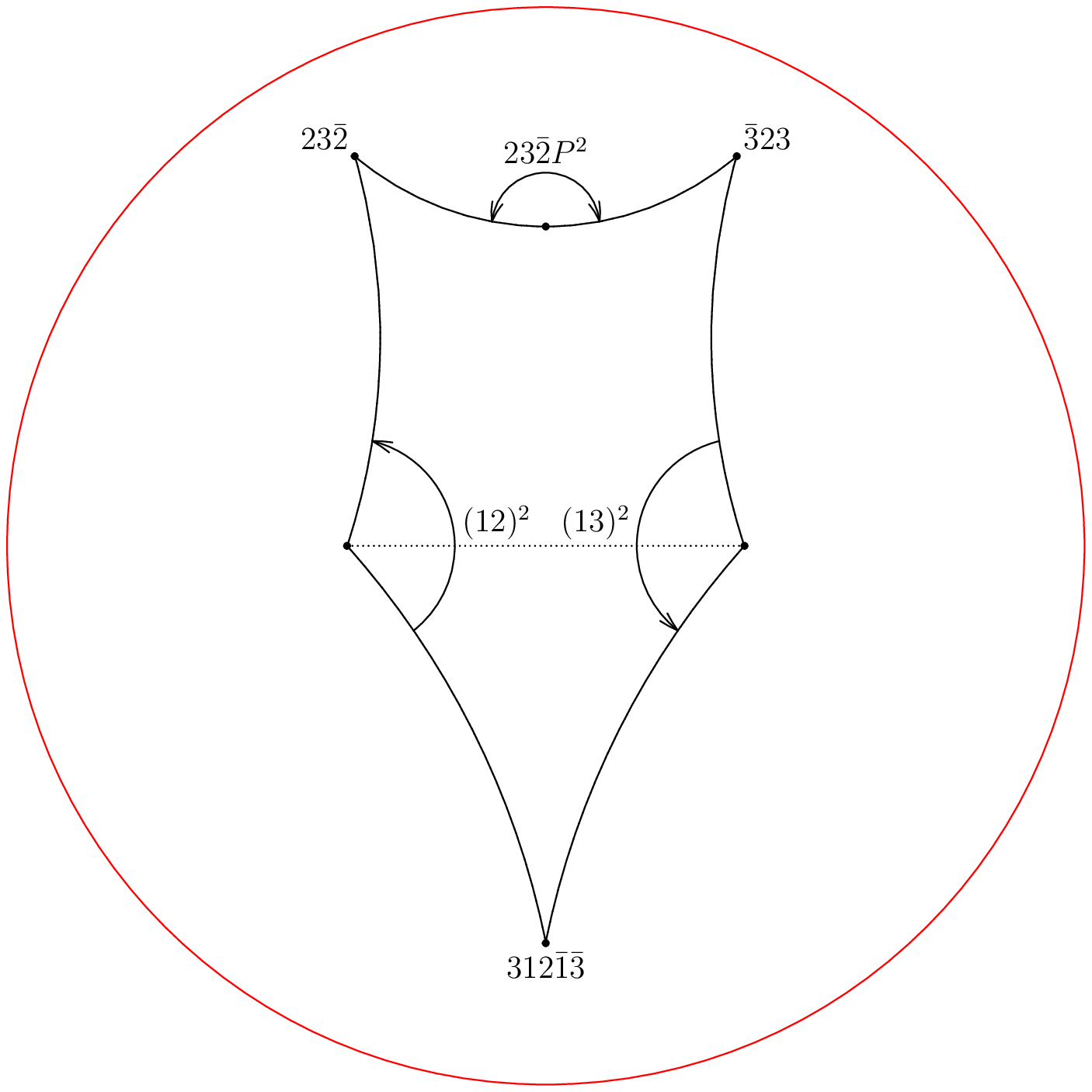}\\
    $p=12$
  \end{tabular}
  \caption{Pairing for the stabilizer of the mirror of $R_1$, for $\Sc(p,\bar\sigma_4)$ groups.}
  \label{fig:stab-s4c}
\end{figure}

\begin{figure}[htbp]
  \centering
  \begin{tabular}{c}
    \includegraphics[width=0.3\textwidth]{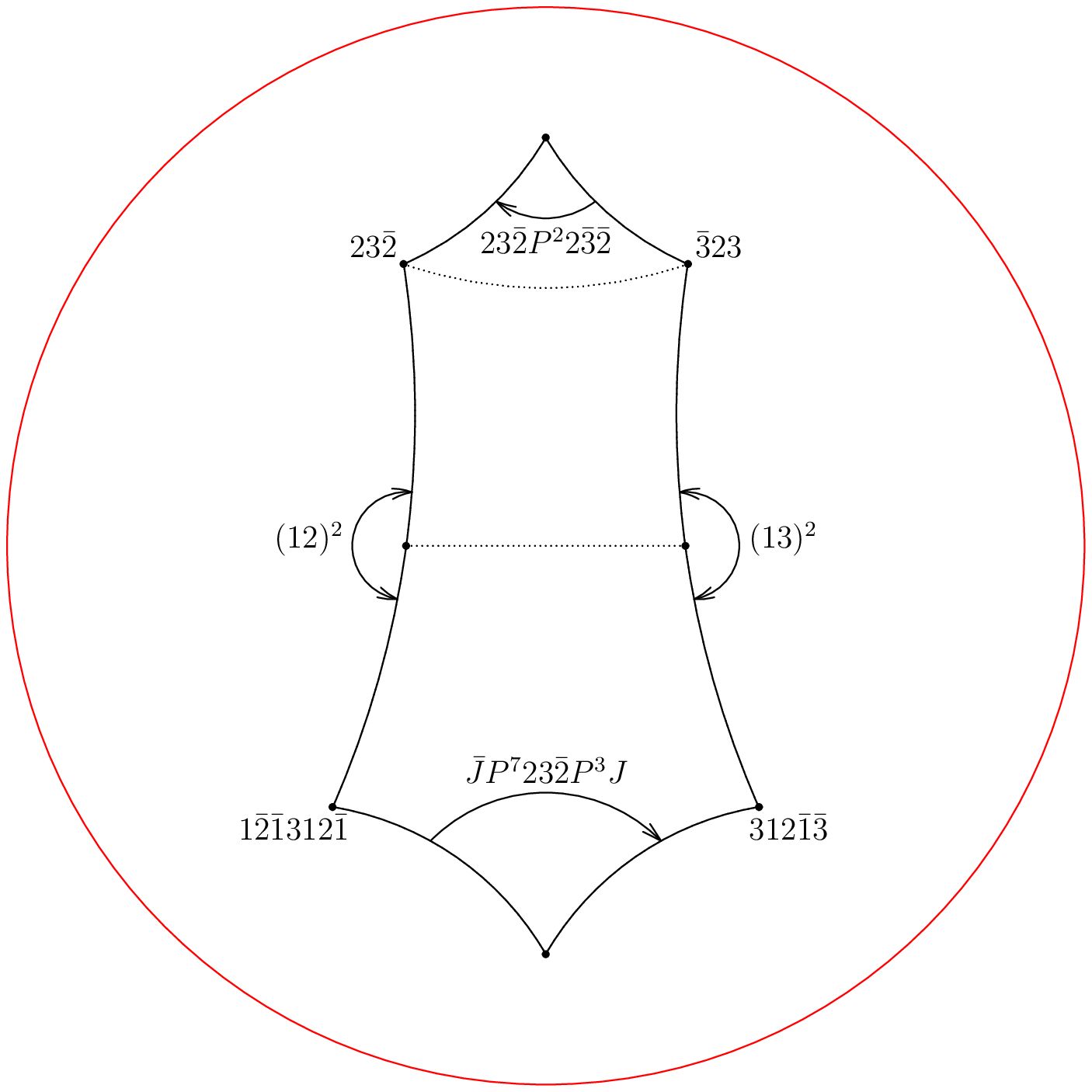}\\
    $p=2$
  \end{tabular}
  \begin{tabular}{c}
    \includegraphics[width=0.3\textwidth]{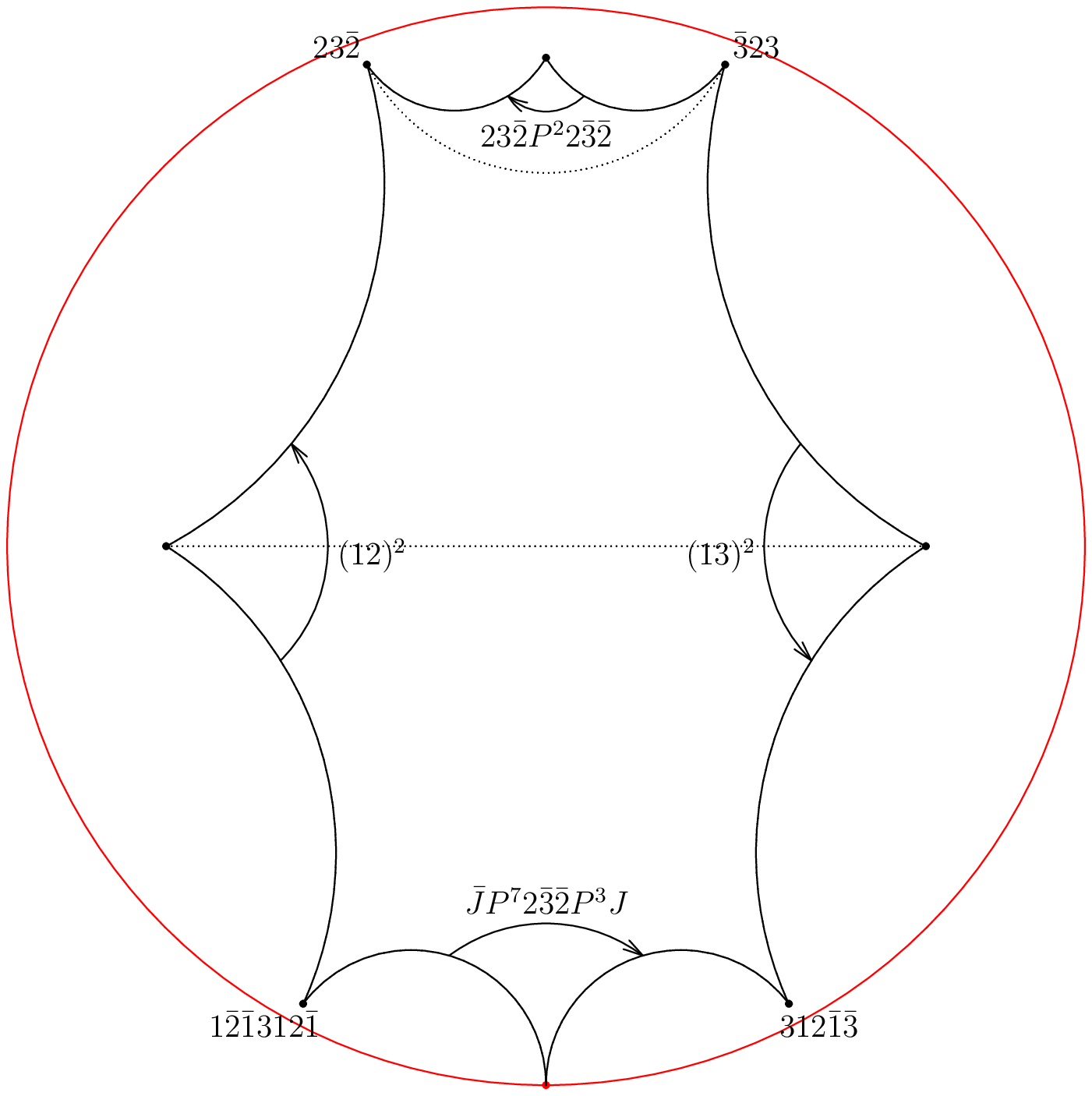}\\
    $p=3$
  \end{tabular}
  \begin{tabular}{c}
    \includegraphics[width=0.3\textwidth]{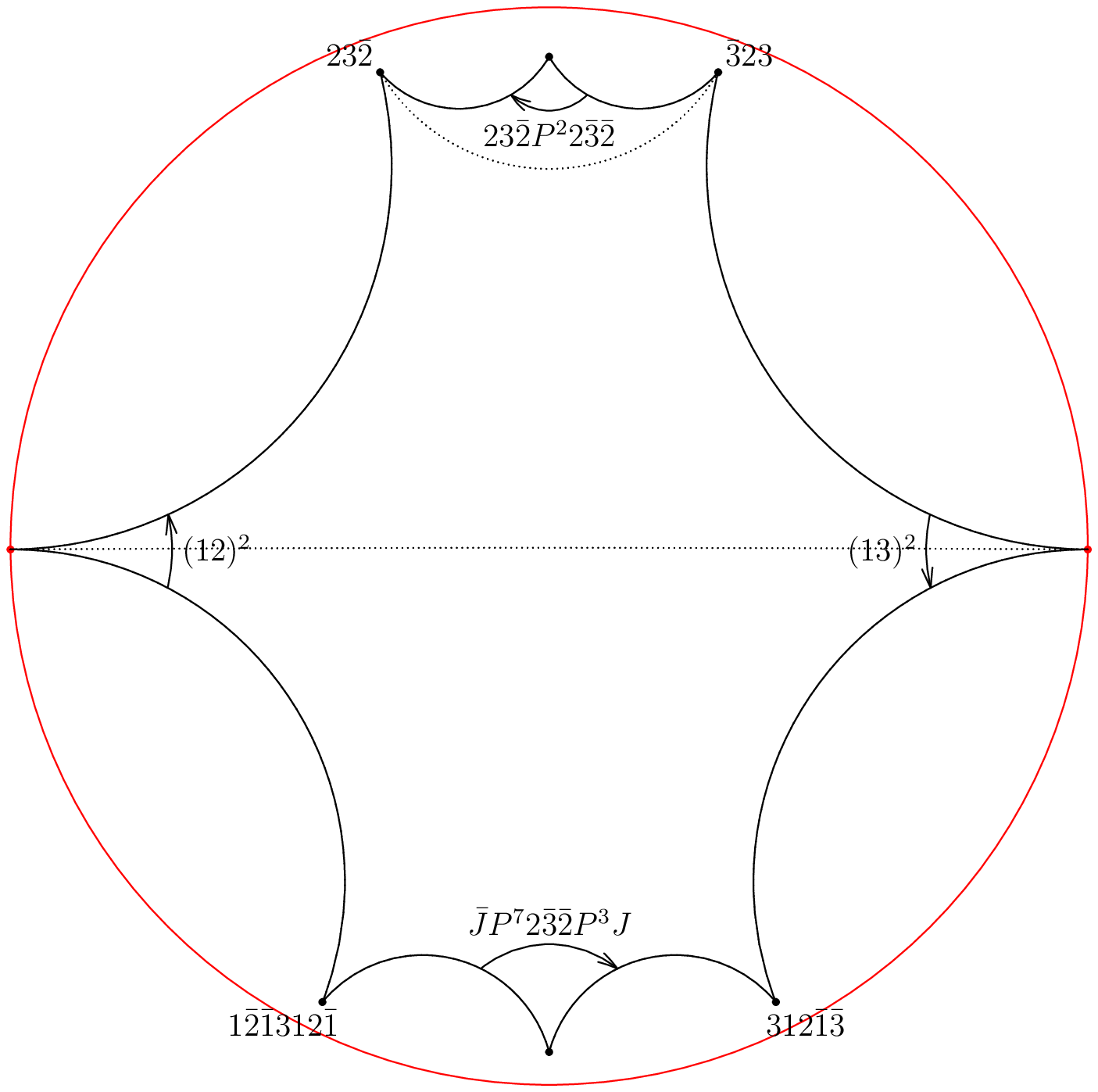}\\
    $p=4$
  \end{tabular}
  \caption{Pairing for the stabilizer of the mirror of $R_1$, for $\Sc(p,\sigma_5)$ groups.}
  \label{fig:stab-s5}
\end{figure}

\begin{figure}[htbp]
  \centering
  \begin{tabular}{c}
    \includegraphics[width=0.3\textwidth]{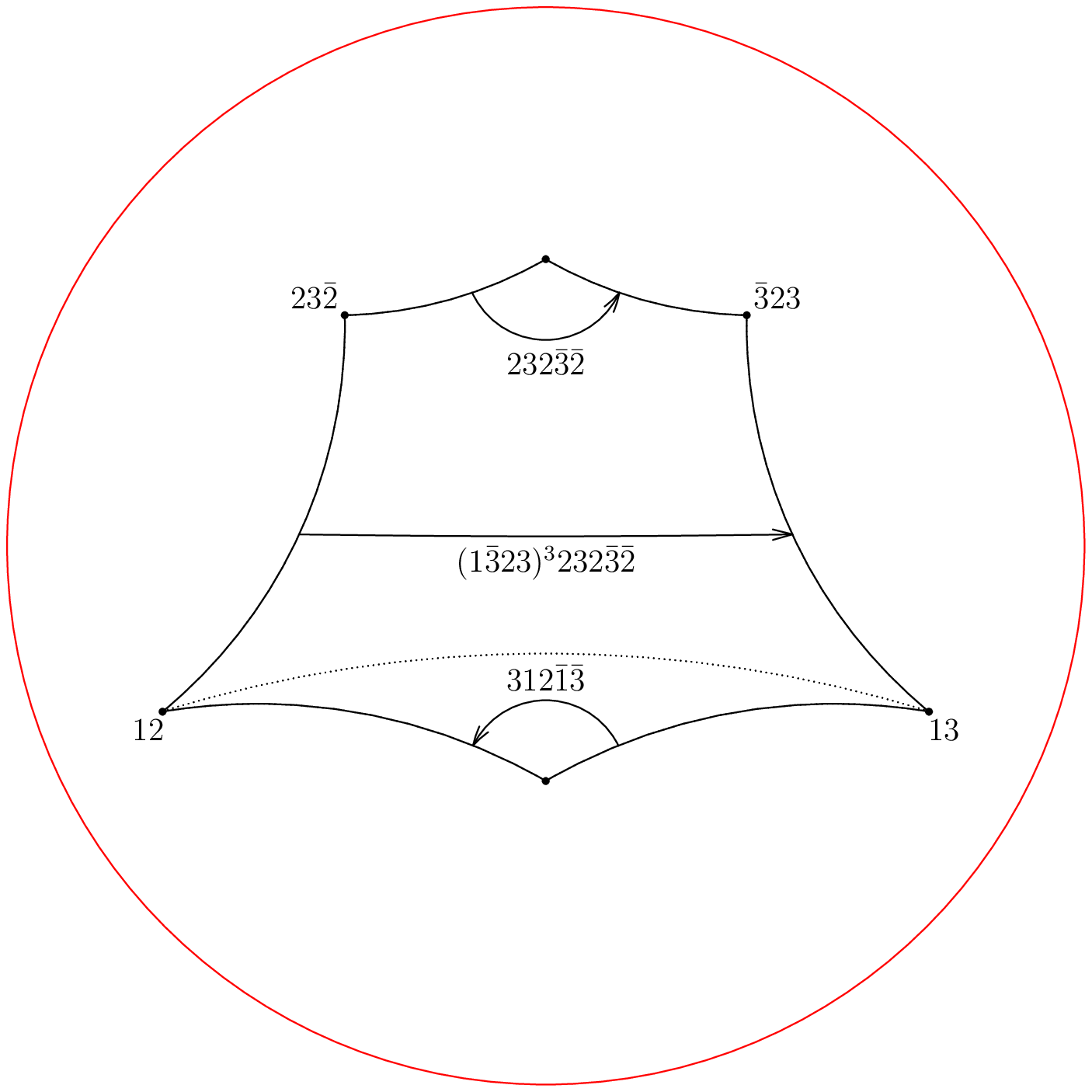}\\
    $p=3$
  \end{tabular}
  \begin{tabular}{c}
    \includegraphics[width=0.3\textwidth]{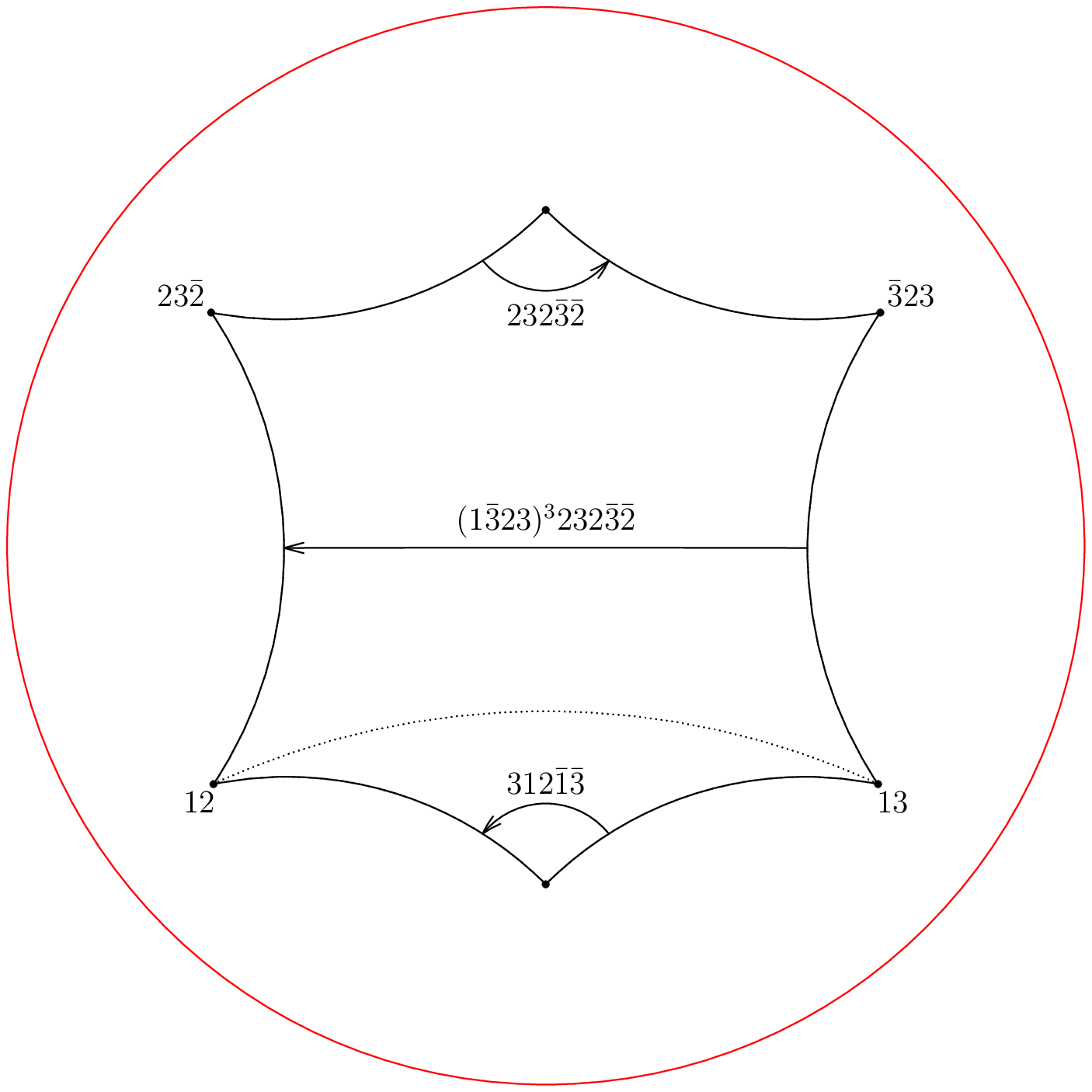}\\
    $p=4$
  \end{tabular}
  \\
  \begin{tabular}{c}
    \includegraphics[width=0.3\textwidth]{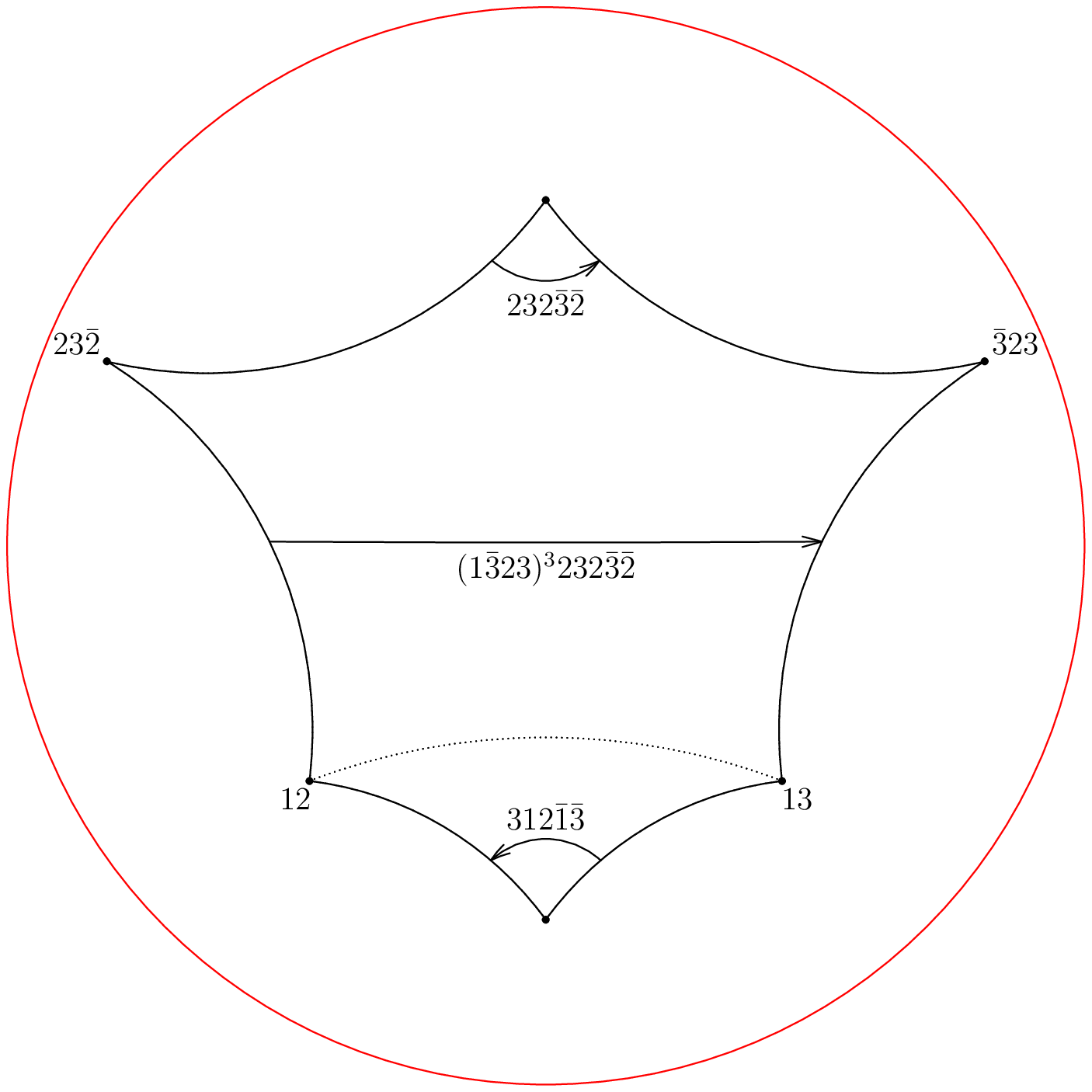}\\
    $p=5$
  \end{tabular}
  \begin{tabular}{c}
    \includegraphics[width=0.3\textwidth]{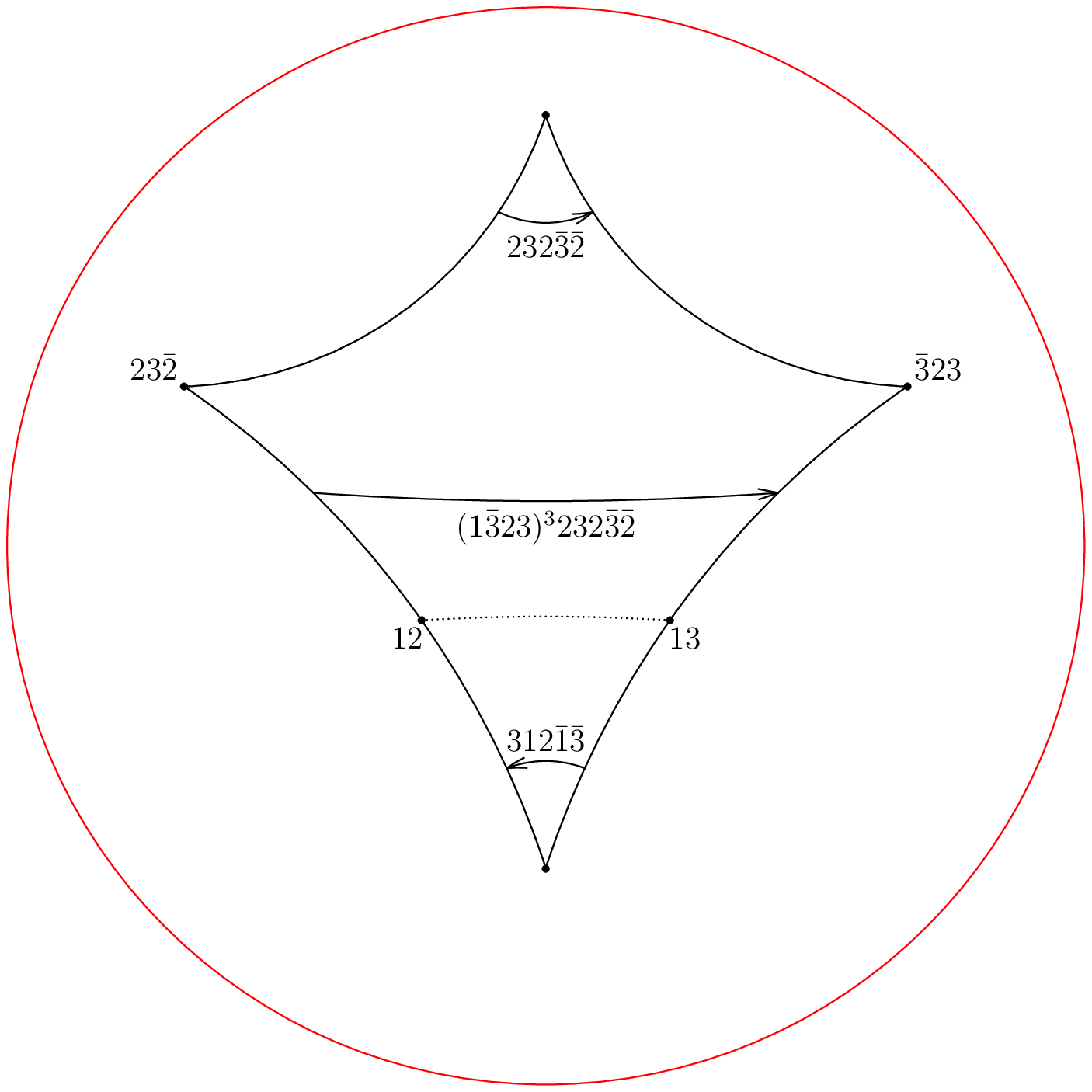}\\
    $p=10$
  \end{tabular}
  \caption{Pairing for the stabilizer of the mirror of $R_1$, for $\Sc(p,\sigma_{10})$ groups.}
  \label{fig:stab-s10}
\end{figure}

\begin{figure}[htbp]
  \centering
  \begin{tabular}{c}
    \includegraphics[width=0.3\textwidth]{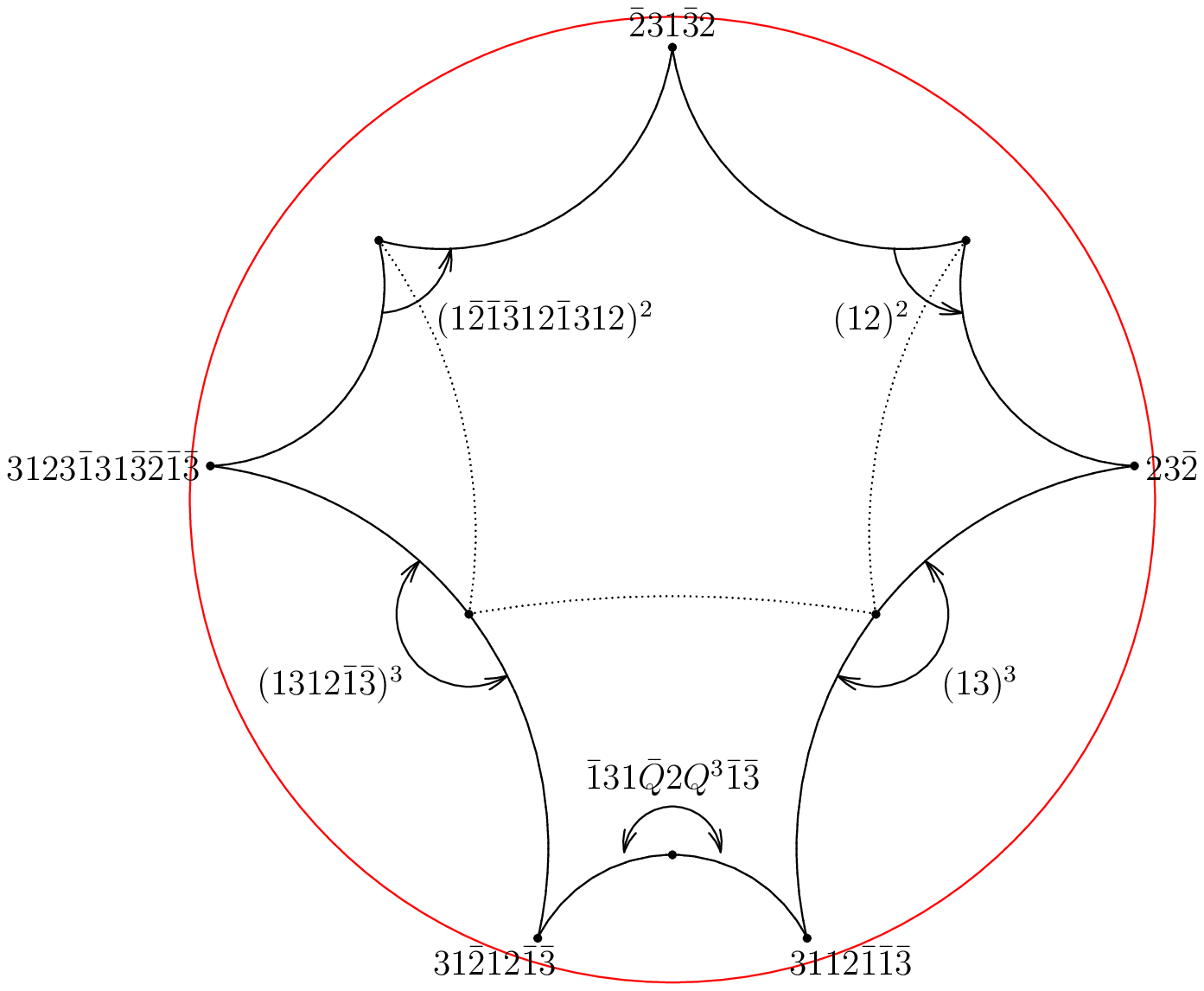}\\
    $p=3$
  \end{tabular}
  \begin{tabular}{c}
    \includegraphics[width=0.3\textwidth]{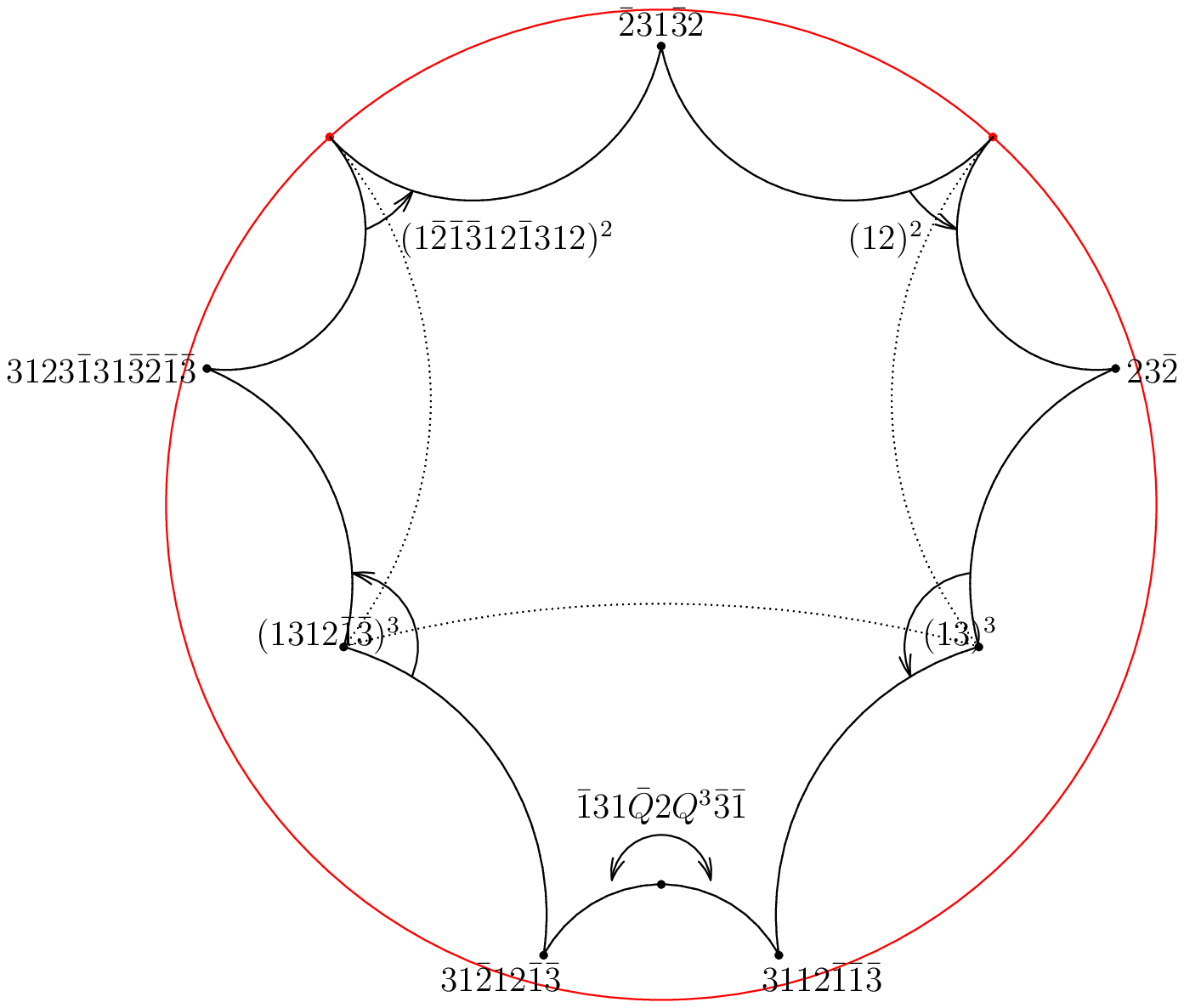}\\
    $p=4$
  \end{tabular}
  \begin{tabular}{c}
    \includegraphics[width=0.3\textwidth]{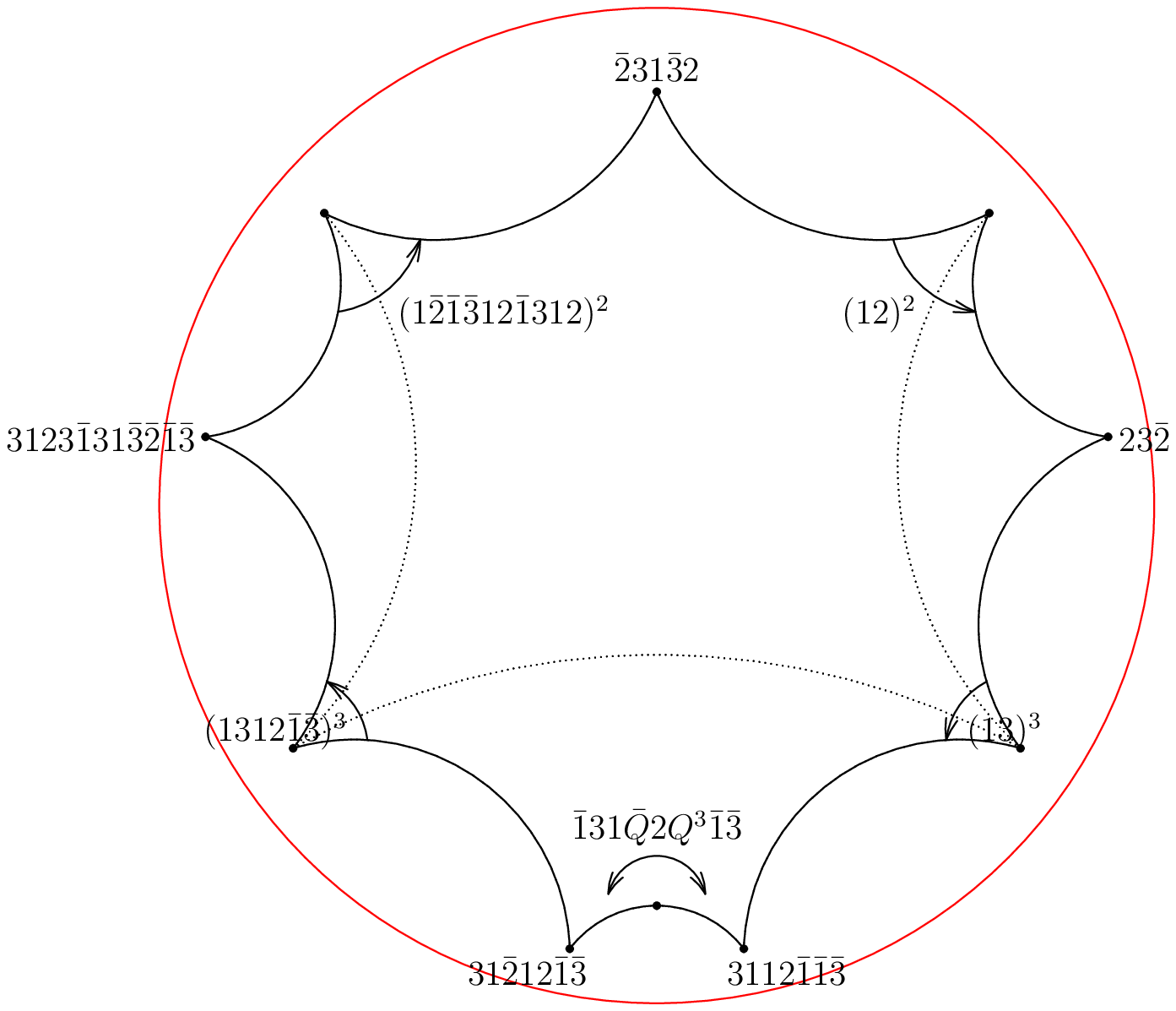}\\
    $p=5$
  \end{tabular}
  \caption{Pairing for the stabilizer of the mirror of $R_1$, for $\Tc(p,{\bf S_2})$ groups.}
  \label{fig:stab-S2}
\end{figure}

\begin{figure}[htbp]
  \centering
  \begin{tabular}{c}
    \includegraphics[width=0.3\textwidth]{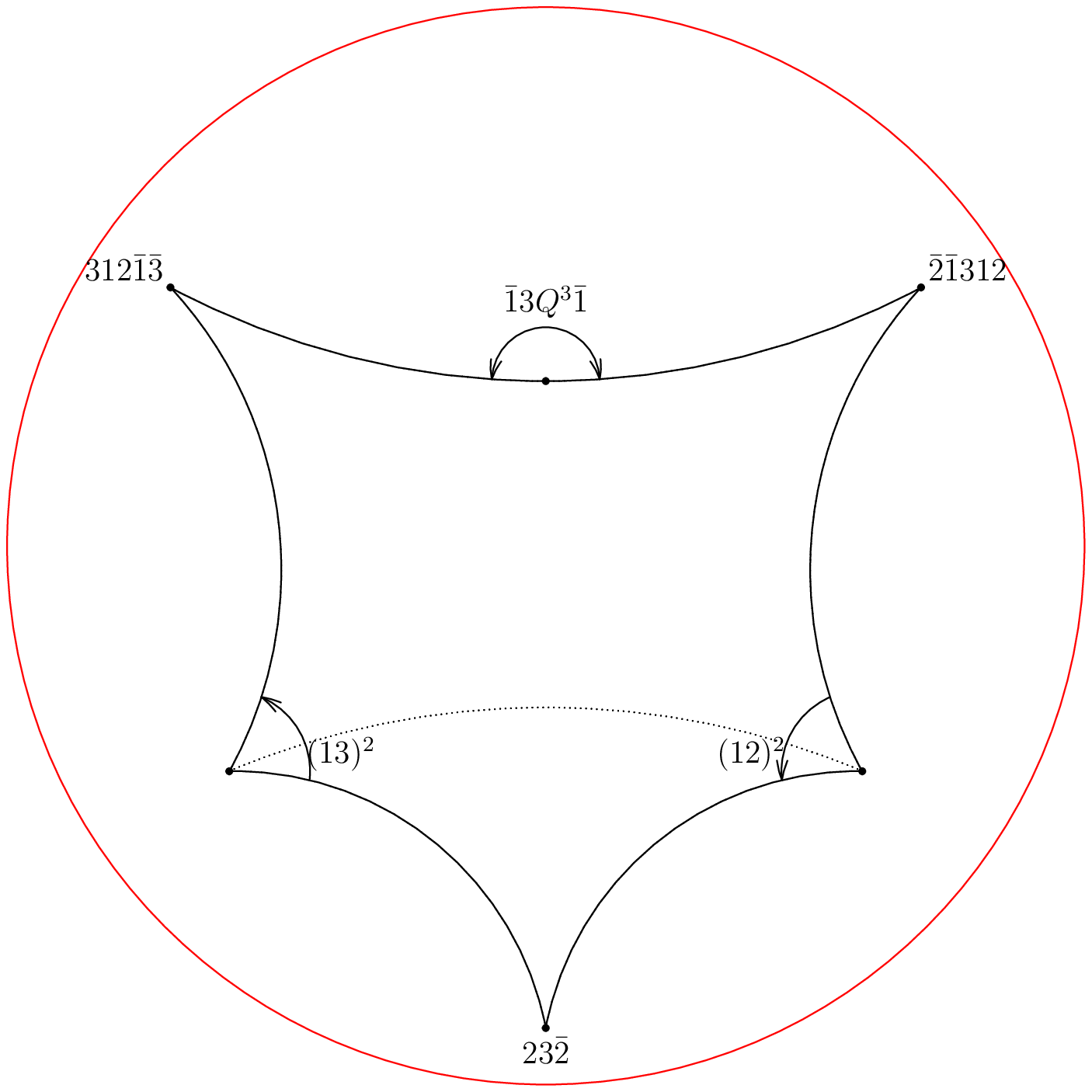}\\
    $p=3$
  \end{tabular}
  \begin{tabular}{c}
    \includegraphics[width=0.3\textwidth]{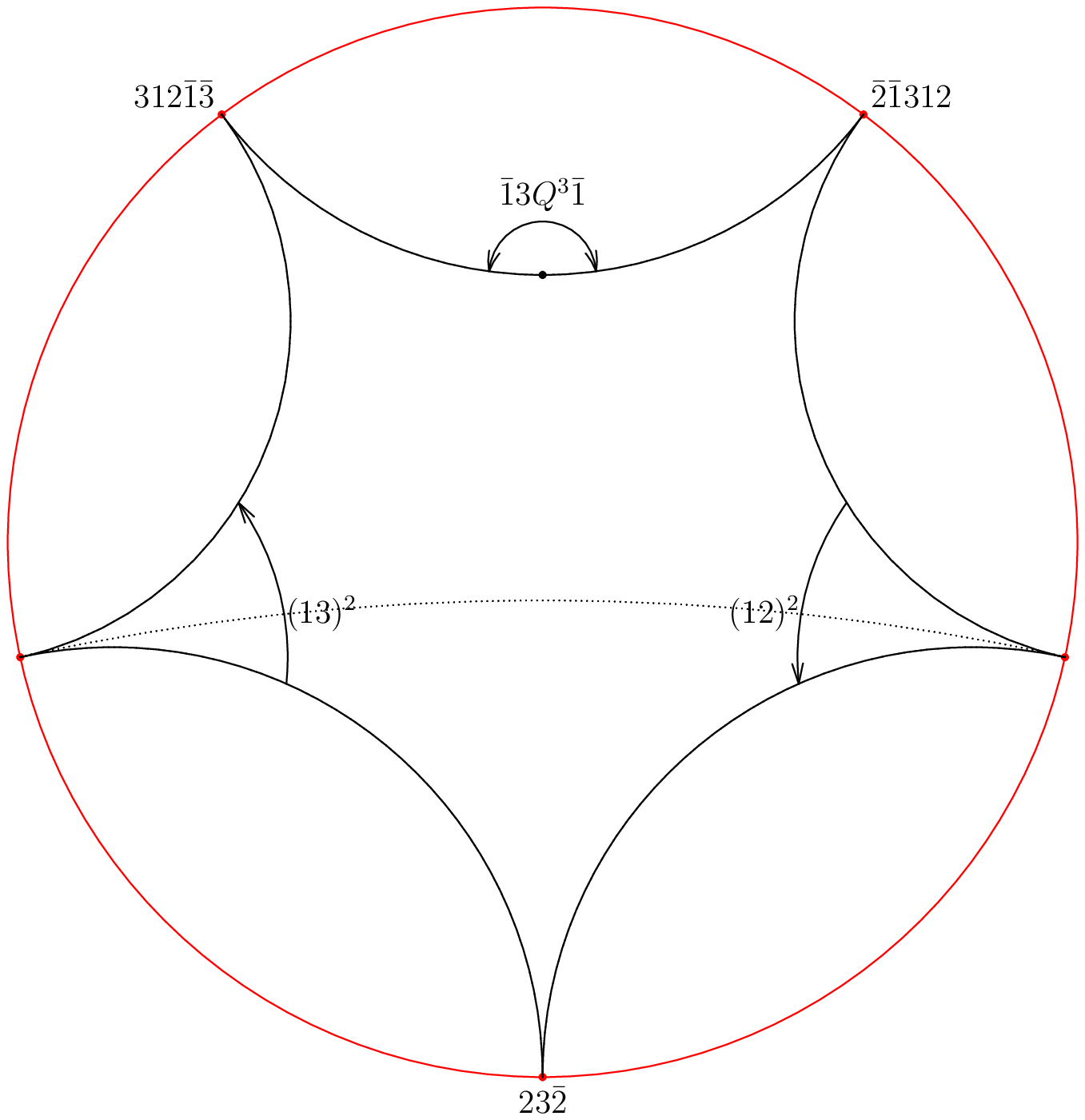}\\
    $p=4$
  \end{tabular}
  \begin{tabular}{c}
    \includegraphics[width=0.3\textwidth]{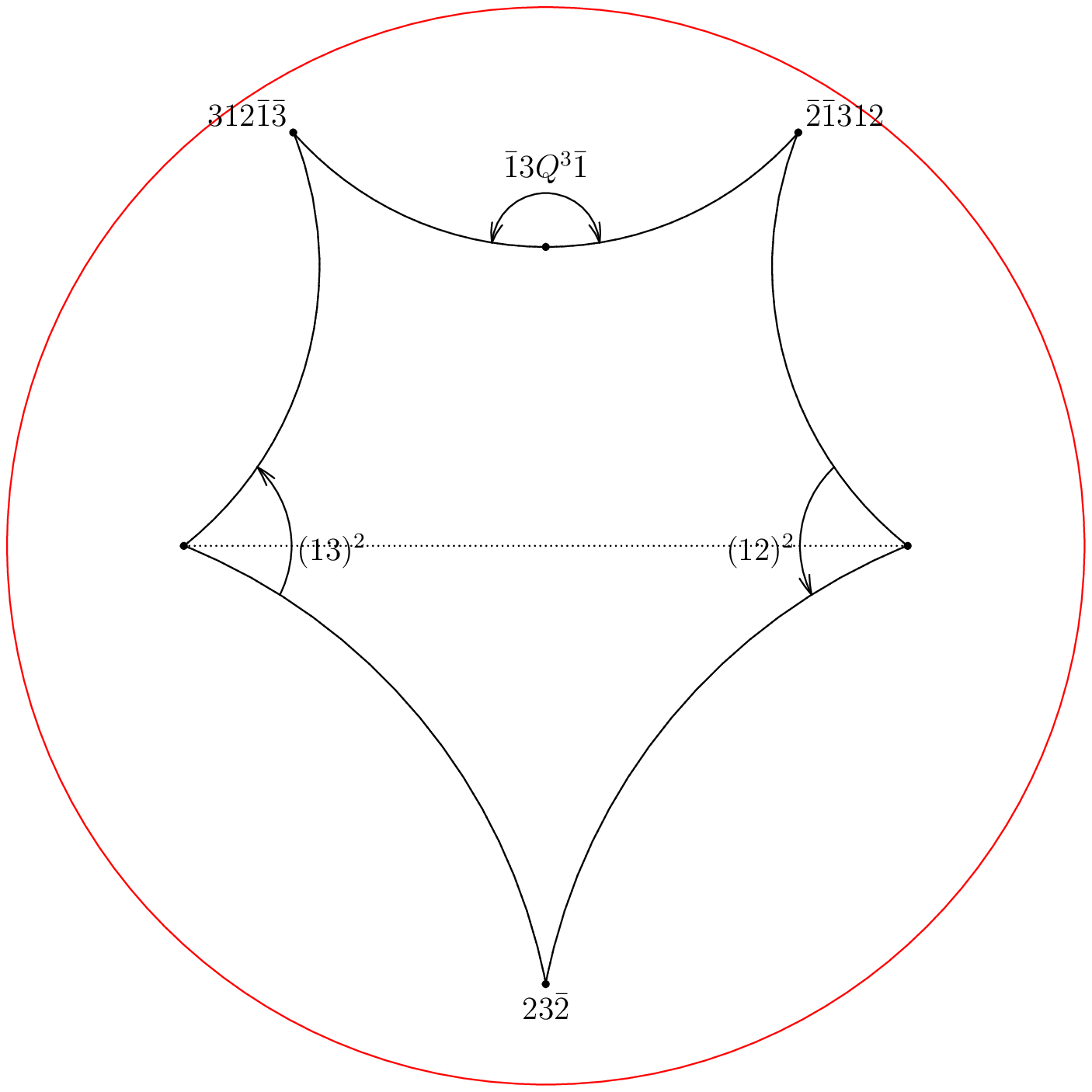}\\
    $p=6$
  \end{tabular}
  \caption{Pairing for the stabilizer of the mirror of $R_1$, for $\Tc(p,{\bf E_2})$ groups.}
  \label{fig:stab-R1-E2}
\end{figure}

\begin{figure}[htbp]
  \centering
  \begin{tabular}{c}
    \includegraphics[width=0.3\textwidth]{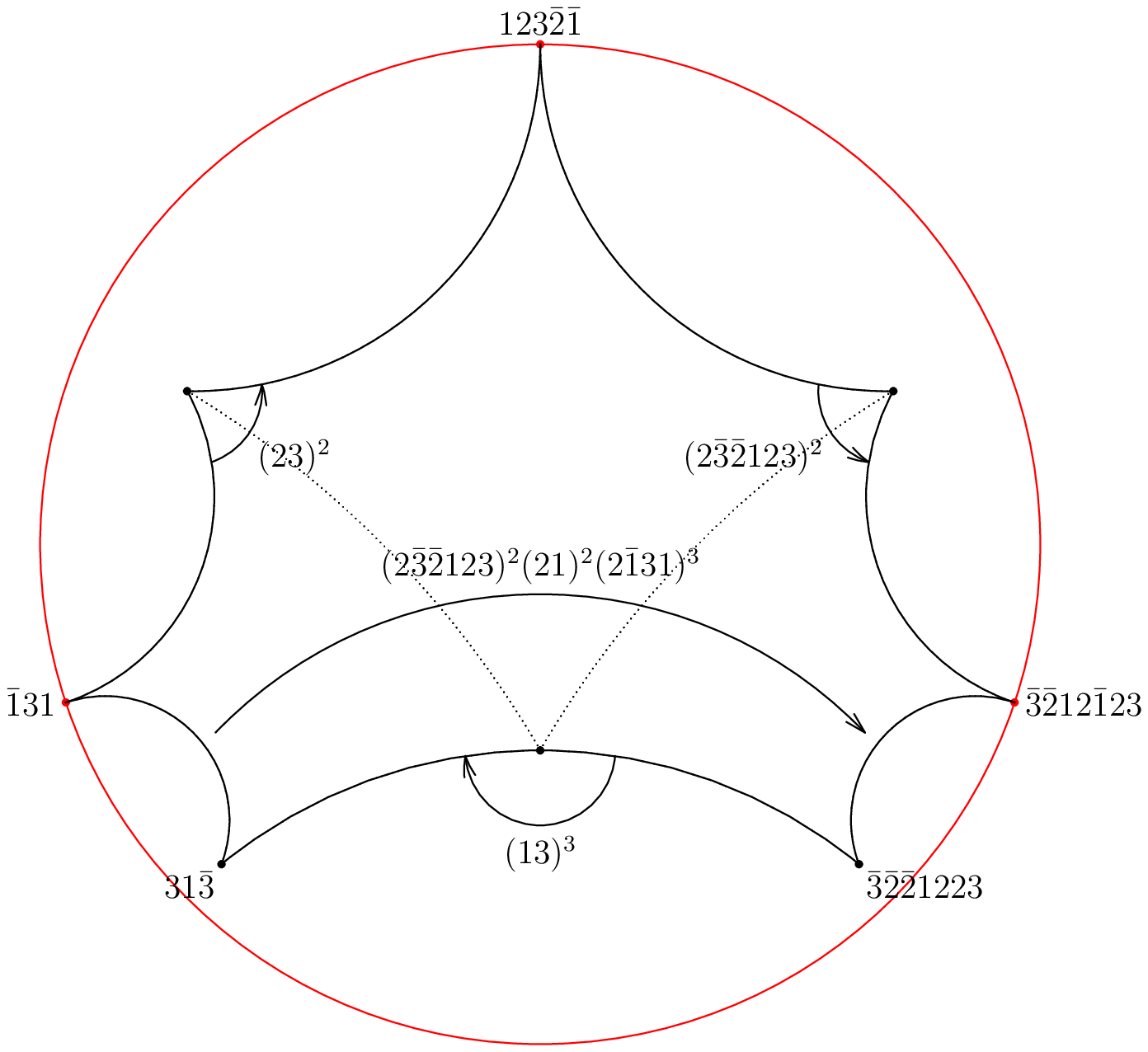}\\
    $p=3$
  \end{tabular}
  \begin{tabular}{c}
    \includegraphics[width=0.3\textwidth]{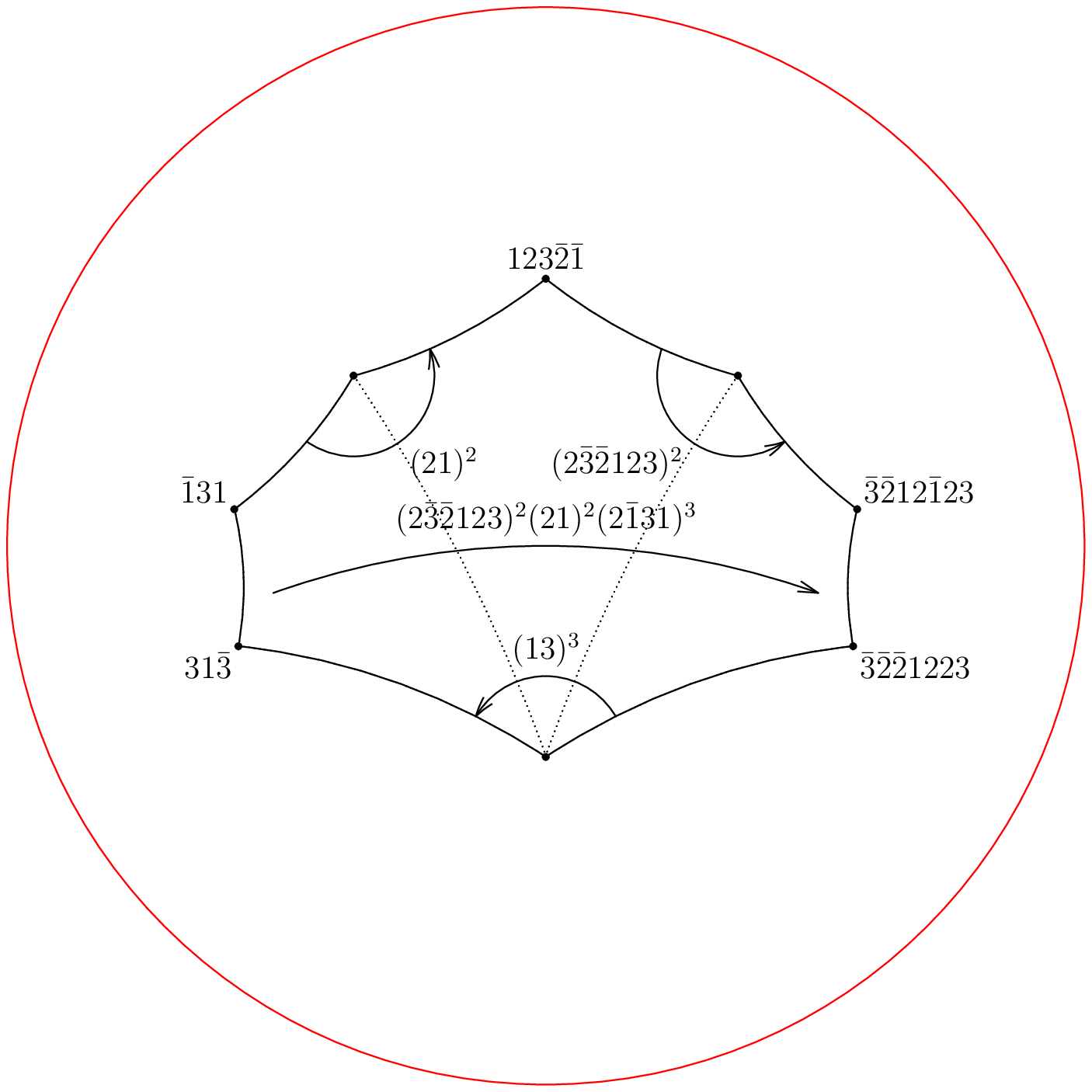}\\
    $p=4$
  \end{tabular}
  \begin{tabular}{c}
    \includegraphics[width=0.3\textwidth]{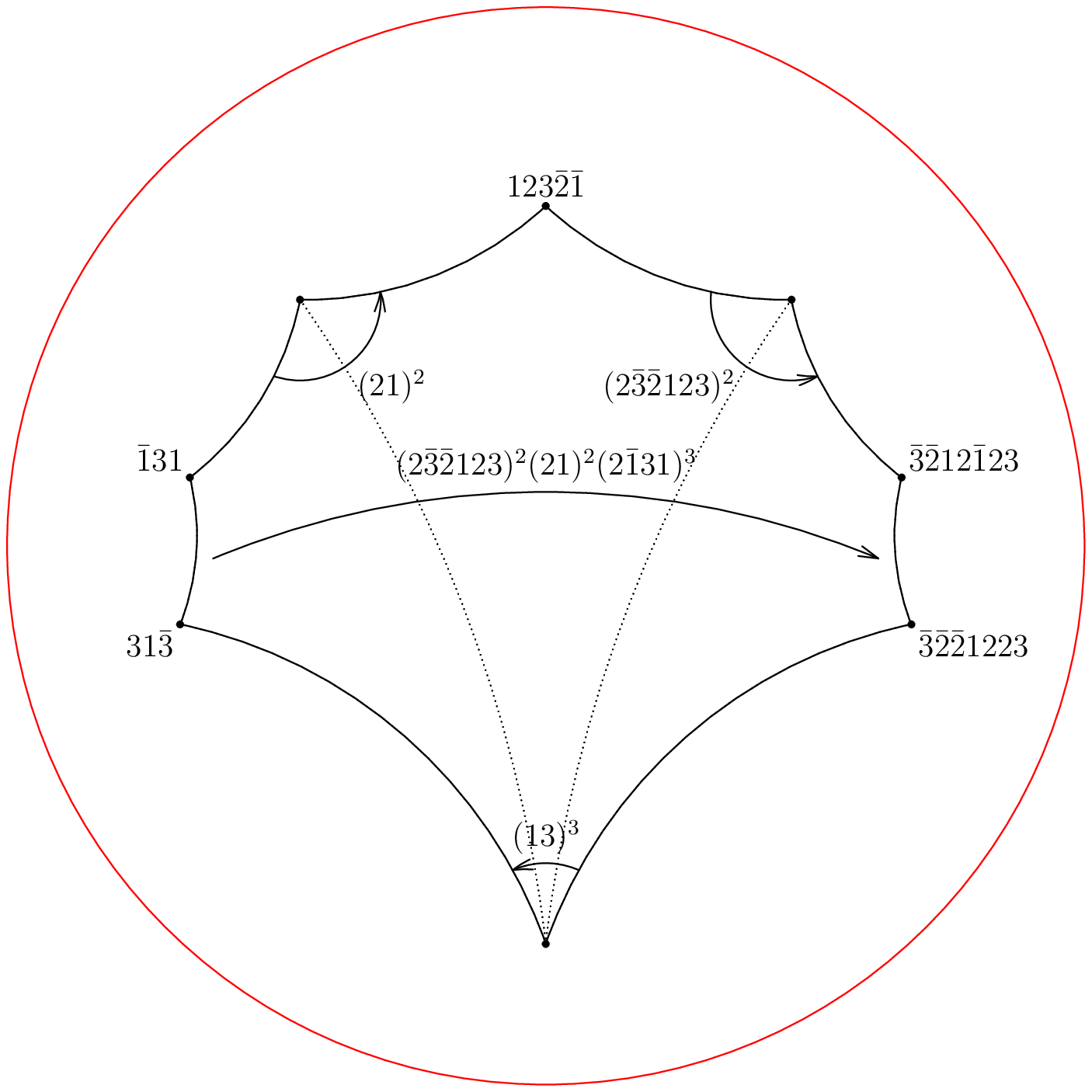}\\
    $p=6$
  \end{tabular}
  \caption{Pairing for the stabilizer of the mirror of $R_2$, for $\Tc(p,{\bf E_2})$ groups.}
  \label{fig:stab-R2-E2}
\end{figure}

\begin{figure}[htbp]
  \centering
  \begin{tabular}{c}
    \includegraphics[width=0.3\textwidth]{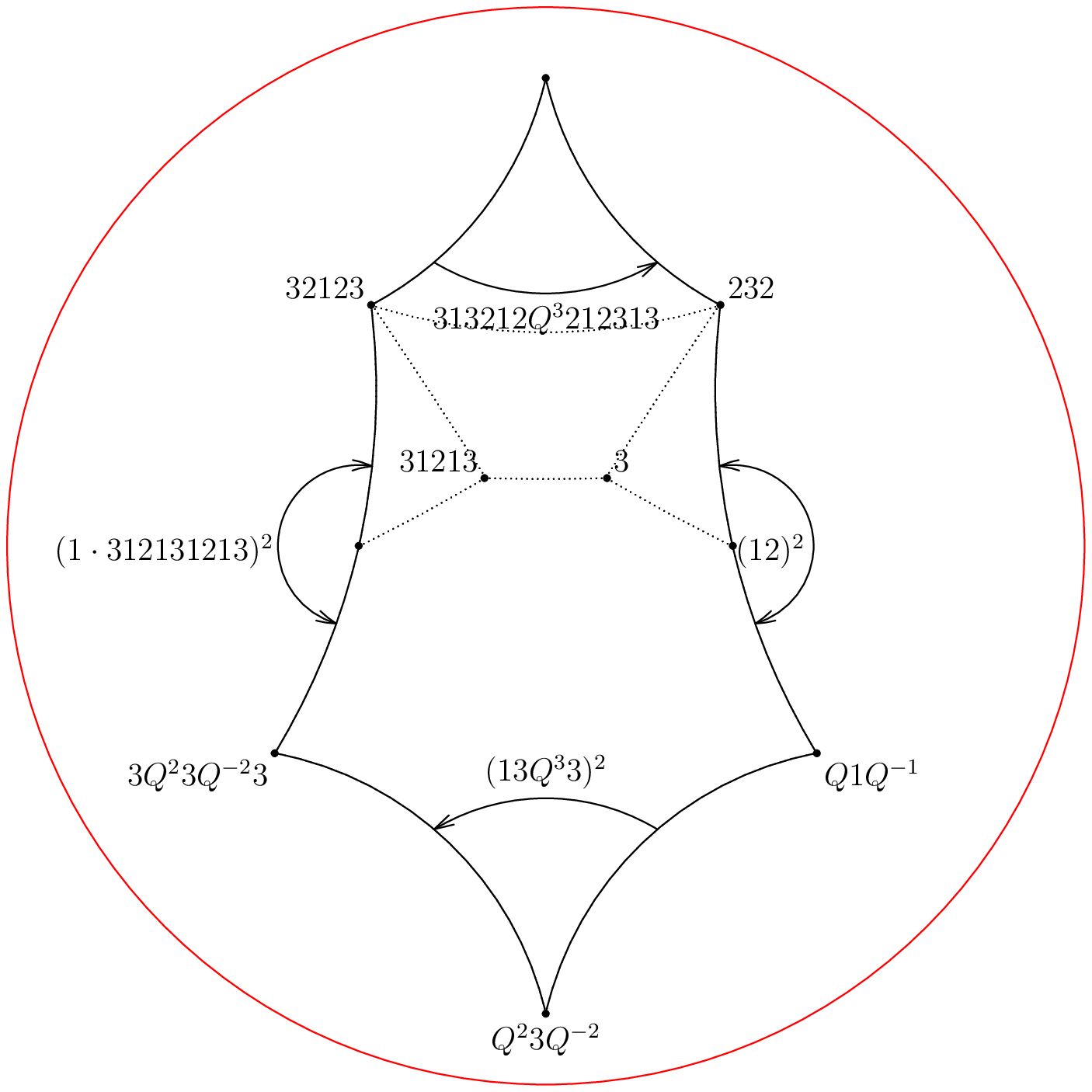}\\
  \end{tabular}
  \caption{Pairing for the stabilizer of the mirror of $R_1$, for $\Tc(2,{\bf H_1})$.}
  \label{fig:stab-H1}
\end{figure}

\begin{figure}[htbp]
  \centering
  \begin{tabular}{c}
    \includegraphics[width=0.3\textwidth]{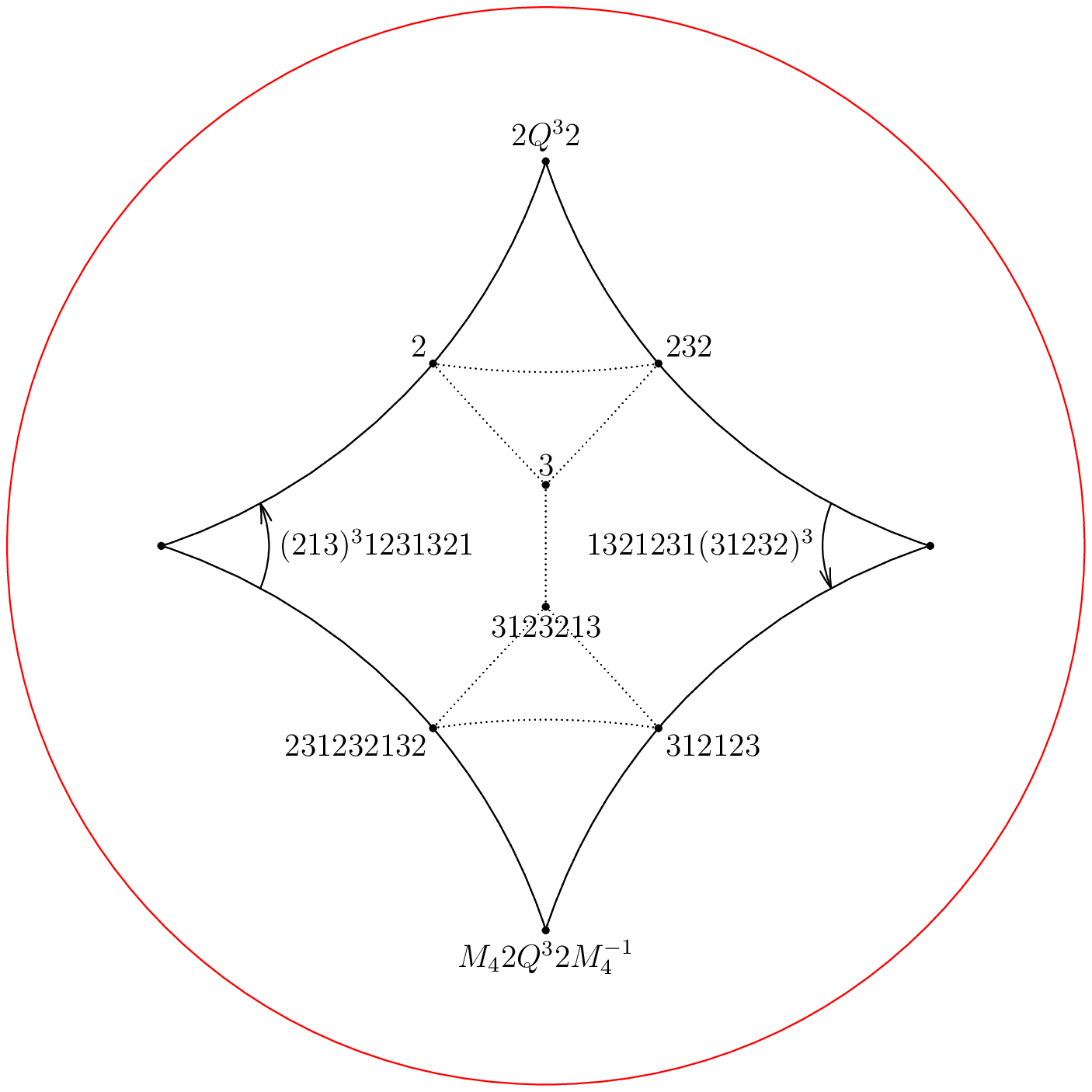}\\
    $p=3$
  \end{tabular}
  \begin{tabular}{c}
    \includegraphics[width=0.3\textwidth]{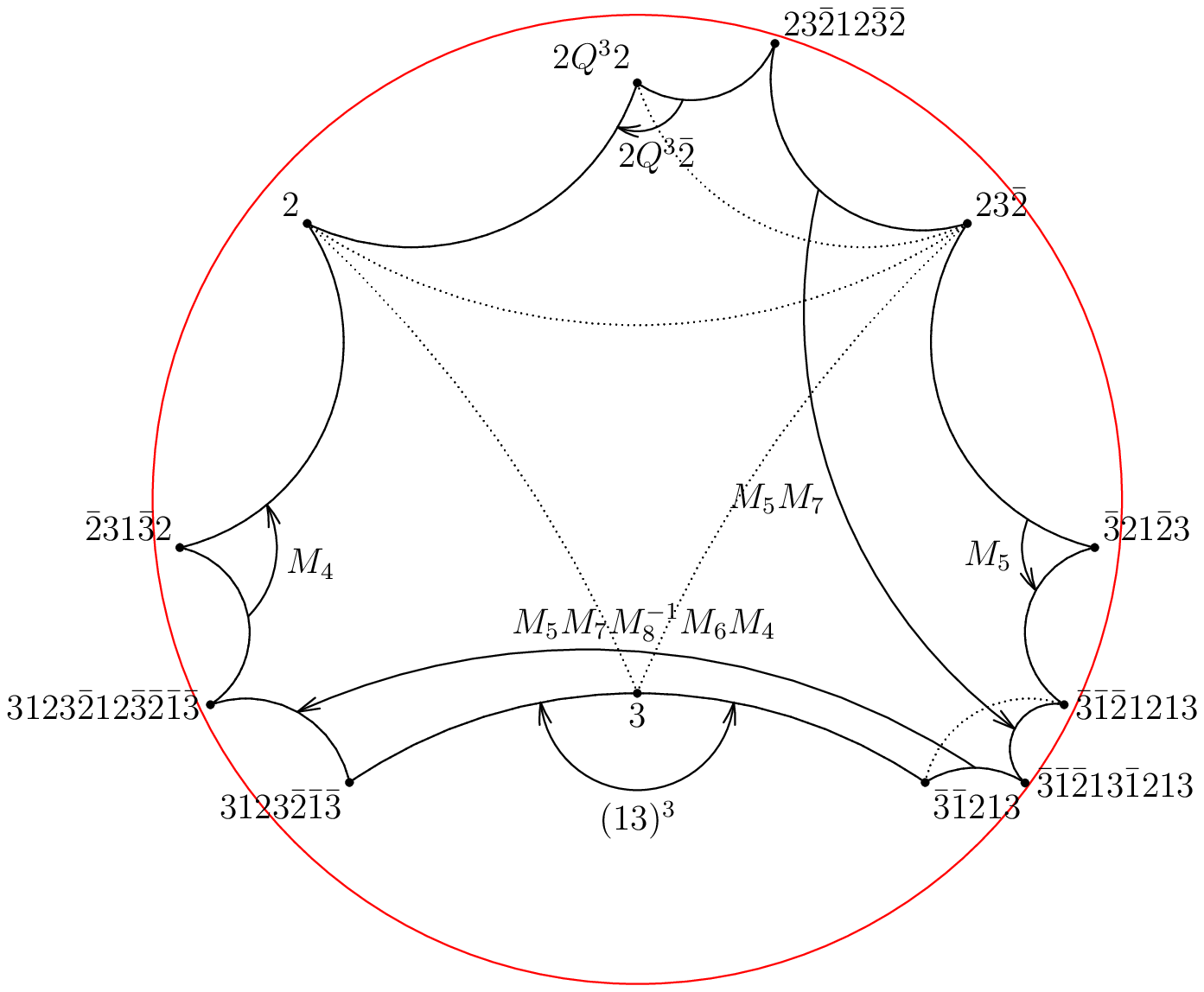}\\
    $p=4$
  \end{tabular}
  \begin{tabular}{c}
    \includegraphics[width=0.3\textwidth]{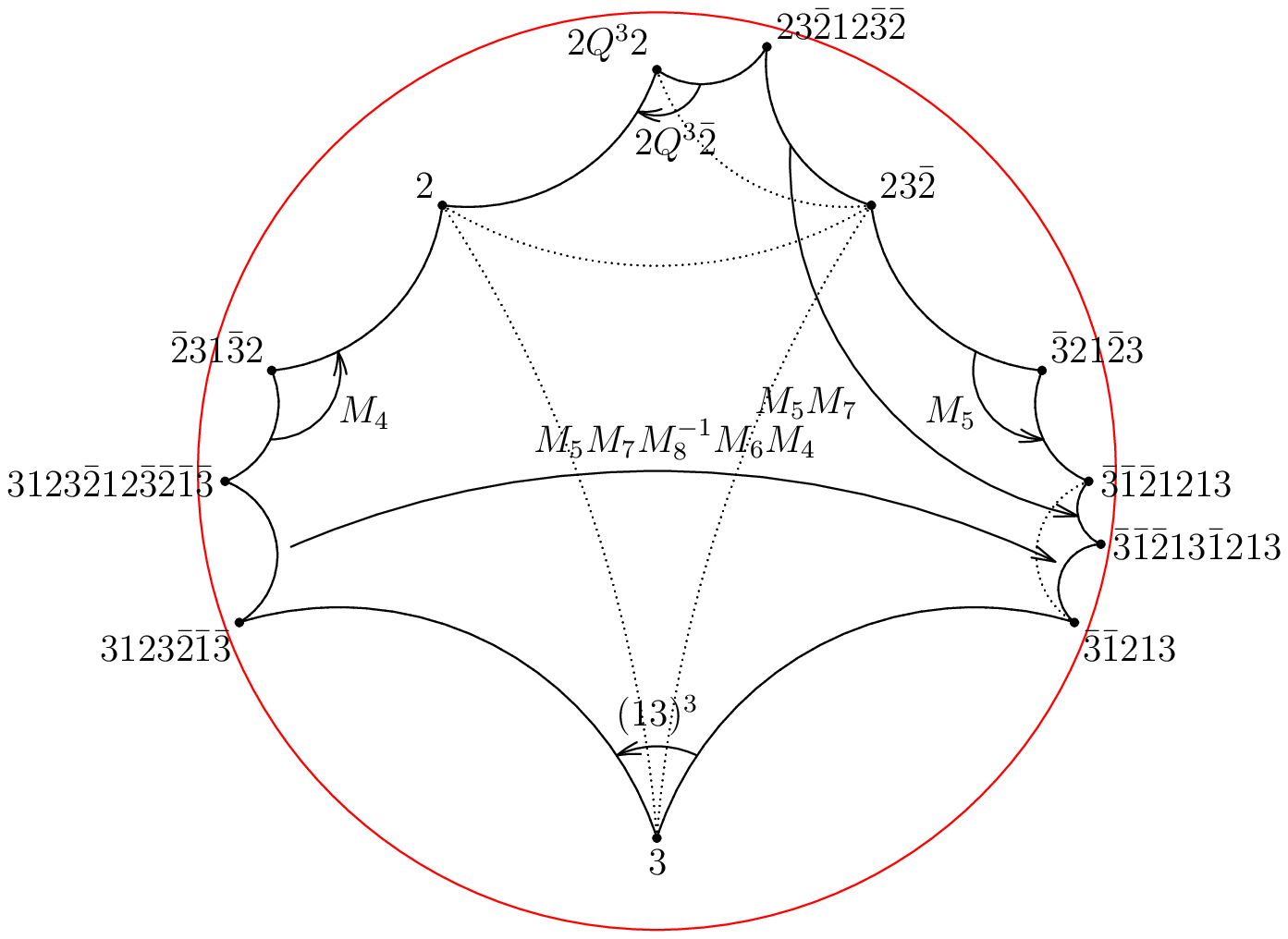}\\
    $p=6$
  \end{tabular}
  \caption{Pairing for the stabilizer of the mirror of $R_1$, for $\Tc(p,{\bf H_2})$ groups.}
  \label{fig:stab-H2}
\end{figure}

We do this only for the mirror of base facets, the other ones are much
easier (in fact, for top and central complex 2-facets, the jigsaw
puzzle has only two pieces) A fundamental domain (as well as
side-pairing transformations) for each of the groups in
Table~\ref{tab:list} is given in
Figures~\ref{fig:stab-s1}-~\ref{fig:stab-H2}.

In a given family $\Sc(p,\tau)$ with fixed $\tau$, the pictures are
very similar, note that the labels do not indicate whether or not the
polytope is truncated at that vertex. For example, for groups
$\Sc(p,\sigma_1)$, the vertex fixed by $(R_1R_2)^3$ is either
$e_1\boxtimes e_2$ (for $p=3$) or $e_1\boxtimes(e_1\boxtimes e_2)$
(for $p=4$ or $6$), see section~\ref{sec:sides}.

\section{Applying the Poincar\'e polyhedron theorem}

In order to determine the stabilizer of a mirror (in restriction to
the mirror), we use the Poincar\'e polyhedron theorem; in fact only
the 2-dimensional special case of that theorem is needed (for a nice
exposition of the special case, see~\cite{derham}; see
also~\cite{beardon}).

We work through a basic example in detail; the statement we prove in
the next few paragraphs is given in Proposition~\ref{prop:stab-s1}. In
the end of this section, we briefly comment on specific difficulties
that come about in treating other cases.

\begin{prop}\label{prop:stab-s1}
  The stabilizer of the mirror of $R_1$ in $\Sc(p,\sigma_1)$ is a
  central extension of a Fuchsian group of signature
  $(0;\frac{p}{p-3},\frac{p}{p-3},\frac{2p}{|p-4|},\frac{2p}{|p-4|},\frac{2p}{|p-6|},\frac{2p}{|p-6|})$,
  with center of order $p$ (generated by $R_1$).
\end{prop}

Recall that the statement about the signature means that the quotient
of the mirror is topologically a compact complex curve of genus 0
(i.e. a copy of $\CP^1$), with 6 cone points of multiplicity given by
$\frac{p}{p-3},\frac{p}{p-3},\frac{2p}{|p-4|},\frac{2p}{|p-4|},\frac{2p}{|p-6|},\frac{2p}{|p-6|})$
respectively. The fact that these numbers are (either $\infty$ or)
integers comes from the fact that the conditions of the Poincar\'e
polyhedron theorem for $\Sc(p,\sigma_1)$ are satisfied
(see~\cite{dpp2}). Note also that a cone point of $\infty$
multiplicity should actually be removed to get a punctured $\CP^1$.

\begin{pf} (of Proposition~\ref{prop:stab-s1})
  
  We invite the reader to have a look at Figure~\ref{fig:stab-s1}
  (page~\pageref{fig:stab-s1}). For each relevant value of $p$, the
  picture describes a geodesic decagon, with vertices that can easily
  be made explicit (see the discussion in section~\ref{sec:sides}, in
  particular for the issue of truncation).

  The arrows in the figure indicate an explicit side-pairing for that
  polygon; the quotient of the mirror of $R_1$ by the action of its
  stabilizer is then the quotient of the decagon under the equivalence
  relation generated by the identifications of the side-pairing.

  The generic points on sides of the decagon are identified with one
  and only one other point on a side, but vertices are usually
  identified with more than two vertices; to make these
  identifications explicit, one needs to track the cycles of vertices.

It is easy to see from the identifications in Figure~\ref{fig:stab-s1}
that the fundamental domain for the stabilizer of the mirror of $R_1$
gives a cell-decomposition of the image of the mirror in the quotient
with $6$ vertices, $5$ edges and a $1$ two-cell, so we get the the
topological Euler characteristic of the image of the mirror (possibly
compactified when some vertices are cusps) is
$$
6-5+1=2,
$$
so the image has genus $0$, and it is isomorphic to $\CP^1$. To figure
out the signature of the corresponding Fuchsian group, we need to
compute the order of the vertex stabilizers for all vertices of the
polygon.

For example, let us describe how to treat the vertex labelled
$\bar2\bar1312$. Recall that this has homogeneous coordinates given by
$e_1\boxtimes R_2^{-1}R_1^{-1}e_3$ (this is a null vector in the case
$p=6$). It cycle looks like the following:
\begin{equation}
  \label{eq:cycle_example}
  1,\bar2\bar31312 \xrightarrow{(12)^3}
  1,2123\bar2\bar1\bar2 \xrightarrow{(123\bar2)^2}
  1,232\bar3\bar2 \xrightarrow{((123\bar2\bar3\bar2)^3(123\bar2)^2(12)^3)^{-1}}
  1,\bar2\bar31312
\end{equation}
Hence the cycle transformation at that vertex is
$$
(12)^{-3}(123\bar2)^{-2}(1232\bar3\bar2)^{-3}(123\bar2)^2(12)^3,
$$
which is conjugate to
$$
(1232\bar3\bar2)^{-3},
$$
and the latter is a complex reflection in a point of order 2 (resp. 4)
if $p=3$ (resp. $4$), or a parabolic element if $p=6$. It acts on the
mirror as a rotation by angle $\pi$ (resp. $\pi/2$) if $p=3$
(resp. $p=4$). 

For this polygon, there is another cycle of the same nature (obtained
by the automorphism of the lattice induced by
$(1,2,3)\mapsto (\bar1,\bar3,\bar2)$), as well as four cycles consisting
of a single point, namely those for $1,2$, for $1,3$, for $1,23\bar2$
and for $1,\bar323$.

The maps $(12)^3$ and $(13)^3$ have the same order, they are parabolic
for $p=3$, and they have order $4$ (resp. $2$) for $p=4$ (resp. $p=6$).
The maps $(123\bar2)^2$ and $(1\bar323)^2$ have the same order, they
are parabolic for $p=4$, and have order $6$ for $p=4$ or $6$.

For $p=3$, we get a Fuchsian group of signature
$(0;2,2,6,6,\infty,\infty)$; more generally, the signature for other
values of $p$ (and also other group orbits of mirrors) are listed in
the table of page~\pageref{tab:s1}.
\end{pf}

From the above discussion, we can also get a presentation for the
stabilizer, which is a central extension (with center of order $p$) of
the corresponding subgroup of $PU(1,1)$. The resulting presentation
for the stabilizer of the mirror $R_1$ in $\Sc(p,\sigma_1)$ is given
in equation~\eqref{eq:pres-stab-s1}: { \small
\begin{equation}
  \label{eq:pres-stab-s1}
  \begin{array}{c}
    \langle\ z, a_{(12)^3},a_{(13)^3},a_{(123\bar2)^2},a_{(1\bar3\bar23)^2},a_{(1232\bar3\bar2)^3},a_{(1\bar3\bar2323)^3}\ |\ z \textrm{ central }, z^p,\hspace{2cm}\\
    a_{(12)^3}^{\frac{p}{p-3}},
    a_{(13)^3}^{\frac{p}{p-3}},
    a_{(123\bar2)^2}^{\frac{2p}{p-4}},
    a_{(1\bar323)^2}^{\frac{2p}{p-4}},
    a_{(1232\bar3\bar2)^2}^{\frac{2p}{p-6}},
    za_{(1\bar3\bar2323)^2}^{\frac{2p}{p-6}}\ \rangle
  \end{array}
\end{equation}
}
Note that the exponents that occur in the presentation are written as
fractions, and the denominator of that fraction vanishes for some
value of $p$; when this happens, the corresponding exponent should be
thought of as being infinite, and the corresponding relation should be
removed from the presentation. Note also that, when finite, the
fractions actually define integers, because the hypotheses of the
Poincar\'e polyhedron theorem hold (see~\cite{dpp2} for more details).

We now briefly sketch the diffulties that are to be overcome for other
families of groups. The general method is the following, where
$\Gamma$ refers to a given lattice triangle group
(see~\ref{tab:list}), and $m$ is the mirror of a given reflection
$a$. We list a number of ressources we have for finding maps that
define a well-defined side-pairing for the polygon. For the simplest
cases, ressource~(1) suffice (e.g. the example treated above).
\begin{enumerate}
\item Find the obvious complex reflections fixing each vertex of the
  polygon (for the mirror of $a$ and vertex labelled $b$, consider
  $(ab)^k$ for the smallest $k\in\IN^*$ such that $(ab)^k$ commutes
  with both $a$ and $b$, see Proposition~\ref{prop:center}). 
\item For each vertex $v$ of the polygon, find $Stab_\Gamma(v)$ by
  tracking vertex cycles in the fundamental polytope for $\Gamma$
  described in~\cite{dpp2}, and deduce $Stab_{Stab_\Gamma(m)}(v)$ by
  taking the subgroup of $Stab_\Gamma(v)$ stabilizing $m$.
\item Same as the previous ressource, replacing $v$ by an edge $e$ of
  the polygon. In some cases, $Stab_\Gamma(e)$ contains a flip,
  i.e. an element of order 2 exchanging the edges of $e$. In that
  case, subdivide the edge $e$ and add an extra vertex in its middle.
\item If the above still does not give a well-defined side-pairing, we
  consider sides stabilized by a complex reflection spower of $P$ or
  $Q$ as in Proposition~\ref{prop:ridge-injection} and replace the
  corresponding complex 2-facets by sectors.
\end{enumerate}
In the next few paragraphs, we sketch the specific difficulties of
each family of groups in Table~\ref{tab:list}.
\begin{itemize}
\item For the groups $\Sc(p,\sigma_1)$, ressource~(1) gives the
  side-pairing transformations $(12)^3$, $(13)^3$, $(123\bar2)^2$,
  $(1\bar323)^2$; the reflections $(1232\bar3\bar2)^3$,
  $(1\bar3\bar2323)^3$, $(1\bar2\bar1312)^3$, $(1312\bar1\bar3)^3$ are
  not used (directly) as side-pairing transformations, since for these
  complex reflections, the image of the polygon meets the original
  polygon only in a vertex. However, the side-pairing transformation
  $(1232\bar3\bar2)^4(123\bar2)^2(12)^3$, which is a product of
  well-chosen complex reflection as above, was found by guessing
  (based on a mix of explicit computation and drawing pictures).
\item For the groups $\Sc(p,\bar\sigma_4)$, ressource~(1) produces
  $(12)^2$ and $(13)^2$, $(123\bar2)^3$, $(1\bar323)^3$,
  $(1312\bar1\bar3)^3$. The first two give side-pairing
  transformations, while the last three do not (the image of the
  polygon intersects the polygon only in a point). Ressource~(3) gives
  the element $23\bar2P^2$, hence we subdivide the top edges in
  Figure~\ref{fig:stab-s4c}.
\item For the groups $\Sc(p,\sigma_5)$, ressource~(1) produces
  $(12)^2$ and $(13)^2$; ressource~(4) gives $23\bar2P^22\bar3\bar2$,
  while $J^{-1}P^723\bar2P^3J$ was obtained by combining (2) and (4).
\item For the groups $\Sc(p,\sigma_{10})$ we use ressource~(1), and
  combine the complex reflections to identify the vertical sides via
  $(1\bar323)^2232\bar3\bar2$.
\item For the groups $\Tc(p,{\bf S_2})$ we use~(1) and~(2) (which
  gives the sides-pairing transformation for the bottom edge).
\item For the mirror of $R_1$ for groups $\Tc(p,{\bf E_2})$, we
  use~(1) and~(2) (which gives the sides-pairing transformation for
  the top edge). For the mirror of $R_2$, we use~(1) and combine the
  corresponding complex reflections in an intricate manner.
\item For the group $\Tc(2,{\bf H_1})$ we use~(1) and~(2) and~(4).
\item For the groups $\Tc(p,{\bf H_2})$, which are by far the most
  complicated to handle, we use all ressources.  
\end{itemize}

\section{Summary of the results} \label{sec:results}

In this section, we gather information on the mirror stabilizers for
all conjugacy classes of complex reflections in lattice complex
hyperbolic triangle groups. For each group and each orbit of complex
reflection mirror, we list
\begin{itemize}
\item Generators for the stabilizer;
\item The order of the fix point stabilizer of the mirror (see the
  column headed ``Order fix pt stab'');
\item The signature of the corresponding Fuchsian subgroup, its
  orbifold Euler characteristic and its area;
\item The field $\IQ(\tr\Gamma^2)$ generated by
  $\textrm{tr}(\gamma^2)$ for $\gamma$ in the Fuchsian group (seen as
  a subgroup of $SU(1,1)$);
\item The arithmeticity of the Fuchsian group (A=arithmetic, NA =
  non-arithmetic).
\end{itemize}
Note that for the mirrors of $R_1$ (or $R_2$ for groups
$\Tc(p,{\bf E_2})$), the trace field $\IQ(\tr\Gamma^2)$ is the same as
the adjoint trace field of the ambient lattice, and the
(non-)arithmeticity is inherited from the ambient lattice.

This can easily be verified by a case by case analysis, for example by
using the method explained in~\cite{jiang-wang-yang}.

\renewcommand{\arraystretch}{1.2}

\begin{table}\label{tab:s1}
  \centering
  {\Large Groups $\Sc(p,\sigma_1)$}\\
\begin{tabular}{c|c|c}
  Reflection     & Values of order $p$ & Generators\\          
  \hline
  $1$            &  3,4,6              & $(12)^3,(13)^3,(123\bar2)^2,(1\bar3\,23)^2,(1232\bar3\,\bar2)^3, (1\bar3\bar2323)^3$ \\
  $(12)^3$       &   4,6               & $1,2$ \\
  $(123\bar2)^2$ &    6                & $1,23\bar2$
\end{tabular}
\ \\
Mirror of $R = R_1$\\
\begin{tabular}{c|c|c|c|c|c|c}
  $p$  & Order fix pt. stab & Signature                   & $\chi^{orb}$    &  Area      & $\IQ(\tr\Gamma^2)$ & A/NA\\
  \hline
  3    & 3                  & $(0;2,2,6,6,\infty,\infty)$ & $-8/3$          &  $16\pi/3$ & $\IQ(\sqrt{6})$   & NA\\
  4    & 4                  & $(0;4,4,4,4,\infty,\infty)$ & $-3$            &  $6\pi$    & $\IQ(\sqrt{2})$   & NA\\
  6    & 6                  & $(0;2,2,6,6,\infty,\infty)$ & $-8/3$          & $16\pi/3$  & $\IQ(\sqrt{6})$   & NA
\end{tabular}
\ \\
Mirror of $R = (R_1R_2)^3$\\
\begin{tabular}{c|c|c|c|c|c|c}
  $p$  & Order fix pt. stab & Signature     & $\chi^{orb}$    &  Area     & $\IQ(\tr\Gamma^2)$ & A/NA\\
  \hline
  4    &  4                 & $(0;3,4,4)$   & $-1/6$          &  $\pi/3$  & $\IQ$             & A\\
  6    &  2                 & $(0;3,6,6)$   & $-1/3$          &  $2\pi/3$ & $\IQ$             & A\\
\end{tabular}
\ \\
Mirror of $(R_1R_2R_3R_2^{-1})^2$\\
\begin{tabular}{c|c|c|c|c|c|c}
  $p$  & Order fix pt. stab & Signature     & $\chi^{orb}$    &  Area    & $\IQ(\tr\Gamma^2)$ & A/NA\\
  \hline
  6    &    6               & $(0;2,6,6)$   & $-1/6$          &  $\pi/3$ & $\IQ$             & A\\
\end{tabular}
\end{table}
\ \\

\begin{table}\label{tab:s4c}
\centering
  {\Large Groups $\Sc(p,\overline{\sigma}_4)$}\\
\begin{tabular}{c|c|c}
  Reflection     & Values of order $p$ & Generators\\          
  \hline
  $1$            &  3,4,5,6,8,12       & $(12)^2,(13)^2,(123\bar2)^3,(1\bar323)^3$,$23\bar2P^2$ \\
  $(12)^2$       &  5,6,8,12           & $1,2$ \\
  $(123\bar2)^2$ &  8,12               & $1,23\bar2$
\end{tabular}
\ \\
Mirror of $R = R_1$\\
\begin{tabular}{c|c|c|c|c|c|c}
  $p$  & Order fix pt. stab & Signature               & $\chi^{orb}$    &  Area         & $\IQ(\tr\Gamma^2)$ & A/NA\\
  \hline
  3    & 3                  & $(0;2,6,6,6)$           & $-1$            &  $2\pi$       & $\IQ(\sqrt{21})$  & A \\
  4    & 4                  & $(0;4,8,\infty,\infty)$ & $-13/8$         &  $13\pi/4$    & $\IQ(\sqrt{7})$   & NA\\
  5    & 5                  & $(0;10,10,10,10)$       & $-8/5$          &  $16\pi/5$    & $\IQ(\sqrt{\frac{5+\sqrt{5}}{14}})$  & NA\\
  6    & 6                  & $(0;6,6,12,\infty)$     & $-19/12$        &  $19\pi/6$    & $\IQ(\sqrt{21})$  & NA\\
  8    & 8                  & $(0;4,4,8,16)$          & $-21/16$        &  $21\pi/8$    & $\IQ(\sqrt{2},\sqrt{7})$  & NA\\
  12   & 12                 & $(0;3,3,4,24)$          & $-25/24$        &  $25\pi/12$   & $\IQ(\sqrt{3},\sqrt{7})$  & NA\\
\end{tabular}
\ \\
Mirror of $R = (R_1R_2)^2$\\
\begin{tabular}{c|c|c|c|c|c|c}
  $p$  & Order fix pt. stab & Signature     & $\chi^{orb}$    &  Area      & $\IQ(\tr\Gamma^2)$ & A/NA\\
  \hline
  5    &  10                & $(0;2,5,5)$   & $-1/10$         &  $\pi/5$   & $\IQ(\sqrt{5})$   & A\\
  6    &  6                 & $(0;2,6,6)$   & $-1/6$          &  $\pi/3$   & $\IQ$             & A\\
  8    &  4                 & $(0;2,8,8)$   & $-1/4$          &  $\pi/2$   & $\IQ(\sqrt{2})$   & A\\
  12   &  3                 & $(0;2,12,12)$ & $-1/3$          &  $2\pi/3$  & $\IQ(\sqrt{3})$   & A\\
\end{tabular}
\ \\
Mirror of $(R_1R_2R_3R_2^{-1})^3$\\
\begin{tabular}{c|c|c|c|c|c|c}
  $p$  & Order fix pt. stab & Signature     & $\chi^{orb}$    &  Area       & $\IQ(\tr\Gamma^2)$ & A/NA\\
  \hline
  8    &    8               & $(0;2,3,8)$   & $-1/24$         &  $\pi/12$   & $\IQ(\sqrt{2})$   & A\\
  12   &    4               & $(0;2,3,12)$  & $-1/12$         &  $\pi/6$    & $\IQ(\sqrt{3})$   & A\\
\end{tabular}
\end{table}
\ \\

\begin{table}\label{tab:s5}
\centering
  {\Large Groups $\Sc(p,\sigma_5)$}\\
\begin{tabular}{c|c|c}
  Reflection     & Values of order $p$ & Generators\\          
  \hline
  $1$            &  2,3,4              & $(12)^2,(13)^2,23\bar2P^52\bar3\bar2, J^{-1}P^723\bar2P^3J$ \\
  $(123\bar2)^5$ &  4                  & $1,23\bar2$\\
  $P^5$          &  2,3,4              & $P^7, 2\bar3\bar2123\bar2$
\end{tabular}
\ \\
Mirror of $R = R_1$\\
\begin{tabular}{c|c|c|c|c|c|c}
  $p$  & Order fix pt. stab & Signature                  & $\chi^{orb}$    &  Area       & $\IQ(\tr\Gamma^2)$        & A/NA\\
  \hline
  2    & 2                  & $(0;2,2,6,6)$              & $-2/3$          &  $4\pi/3$   & $\IQ(\sqrt{5})$          & A\\
  3    & 3                  & $(0;6,6,6,6,\infty)$       & $-7/3$          &  $14\pi/3$  & $\IQ(\sqrt{5})$          & NA\\
  4    & 4                  & $(0;4,6,12,\infty,\infty)$ & $-5/2$          &  $5\pi$     & $\IQ(\sqrt{3},\sqrt{5})$ & NA\\
\end{tabular}
\ \\
Mirror of $(R_1R_2R_3R_2^{-1})^5$\\
\begin{tabular}{c|c|c|c|c|c|c}
  $p$  & Order fix pt. stab & Signature     & $\chi^{orb}$    &  Area      & $\IQ(\tr\Gamma^2)$  & A/NA\\
  \hline
  4    &    4               & $(0;2,4,5)$   & $-1/20$         &  $\pi/10$  & $\IQ(\sqrt{5})$    & A\\
\end{tabular}
\ \\
Mirror of $R = P^5$\\
\begin{tabular}{c|c|c|c|c|c|c}
  $p$  & Order fix pt. stab & Signature         & $\chi^{orb}$    &  Area        & $\IQ(\tr\Gamma^2)$  & A/NA\\
  \hline 
  2    &  6                 & $(0;2,5,6)$       & $-2/15$         &  $4\pi/15$   & $\IQ(\sqrt{5})$    & A\\
  3    &  6                 & $(0;3,5,\infty)$  & $-7/15$         &  $14\pi/15$  & $\IQ(\sqrt{5})$    & NA\\
  4    &  6                 & $(0;4,5,12)$      & $-7/15$         &  $14\pi/15$  & $\IQ(\sqrt{3},\sqrt{5})$ & NA\\
\end{tabular}
\end{table}
\ \\

\begin{table}\label{tab:s10}
\centering
  {\Large Groups $\Sc(p,\sigma_{10})$}\\
\begin{tabular}{c|c|c}
  Reflection     & Values of order $p$ & Generators\\          
  \hline
  $1$            &  3,4,5,10           & $312\bar1\bar3, 232\bar3\bar2, (123\bar2)^3, (1\bar323)^3$ \\
  $(12)^5$       &  4,5,10             & $1,2$\\
  $(123\bar2)^3$ &  10                 & $1,23\bar2$\\
\end{tabular}
\ \\
Mirror of $R = R_1$\\
\begin{tabular}{c|c|c|c|c|c|c}
  $p$  & Order fix pt. stab & Signature           & $\chi^{orb}$  &  Area      & $\IQ(\tr\Gamma^2)$  & A/NA \\
  \hline
  3    & 3                  & $(0;2,3,3,6)$       & $-2/3$        &  $4\pi/3$  & $\IQ(\sqrt{5})$ & A\\
  4    & 4                  & $(0;4,4,4,4)$       & $-1$          &  $2\pi$    & $\IQ(\sqrt{5})$ & A\\
  5    & 5                  & $(0;2,5,5,10)$      & $-1$          &  $2\pi$    & $\IQ(\sqrt{5})$ & A\\
  10   & 10                 & $(0;5,10,10)$       & $-3/5$        &  $6\pi/5$  & $\IQ(\sqrt{5})$ & A\\
\end{tabular}
\ \\
Mirror of $R = (R_1R_2)^5$\\
\begin{tabular}{c|c|c|c|c|c|c}
  $p$  & Order fix pt. stab & Signature         & $\chi^{orb}$    &  Area     & $\IQ(\tr\Gamma^2)$  & A/NA\\
  \hline 
  4    &  4                 & $(0;2,4,5)$       & $-1/20$         &  $\pi/10$ & $\IQ(\sqrt{5})$ & A\\
  5    &  2                 & $(0;2,5,5)$       & $-1/10$         &  $\pi/5$  & $\IQ(\sqrt{5})$ & A\\
  10   &  1                 & $(0;2,5,10)$      & $-1/5$          &  $2\pi/5$ & $\IQ(\sqrt{5})$ & A\\
\end{tabular}
\ \\
Mirror of $(R_1R_2R_3R_2^{-1})^3$\\
\begin{tabular}{c|c|c|c|c|c|c}
  $p$  & Order fix pt. stab & Signature     & $\chi^{orb}$    &  Area       & $\IQ(\tr\Gamma^2)$  & A/NA\\
  \hline
  10   &    5               & $(0;2,3,10)$  & $-1/15$         &  $2\pi/15$  & $\IQ(\sqrt{5})$    & A\\
\end{tabular}
\end{table}
\ \\

\begin{table}\label{tab:S2}
\centering
  {\Large Groups $\Tc(p,{\bf S_2})$}\\
\begin{tabular}{c|c|c}
  Reflection     & Values of order $p$ & Generators\\          
  \hline
  $1$            &  3,4,5              & $(131\bar212\bar1\bar3)^5,(13112\bar1\bar1\bar3)^5,(1312\bar1\bar3)^3,(13)^3,$ \\
                 &                     &  $\quad (13123\bar131\bar3\bar2\bar1\bar3)^5,(123\bar2)^5,(1\bar2\bar1\bar312\bar1312)^2,(12)^2$\\

  $(12)^2$       &  5                  & $1,2$\\
  $(123\bar2)^5$ &  4,5                & $1,23\bar2$\\
\end{tabular}
\ \\
Mirror of $R = R_1$\\
\begin{tabular}{c|c|c|c|c|c|c}
  $p$  & Order fix pt. stab & Signature                   & $\chi^{orb}$    &  Area        & $\IQ(\tr\Gamma^2)$           & A/NA\\
  \hline
  3    & 3                  & $(0;2,2,2,6,6,6)$           & $-2$            &  $4\pi$      & $\IQ(\sqrt{5})$             & A\\
  4    & 4                  & $(0;2,4,4,4,\infty,\infty)$ & $-11/4$         &  $11\pi/2$   & $\IQ(\sqrt{3},\sqrt{5})$    & NA\\
  5    & 5                  & $(0;2,2,10,10,10,10)$       & $-13/5$         &  $26\pi/5$   & $\IQ(\cos\frac{2\pi}{15})$  & NA\\
\end{tabular}
\ \\
Mirror of $R = (R_1R_2)^2$\\
\begin{tabular}{c|c|c|c|c|c|c}
  $p$  & Order fix pt. stab & Signature         & $\chi^{orb}$    &  Area     & $\IQ(\tr\Gamma^2)$           & A/NA\\
  \hline 
  5    &  10                & $(0;2,5,5)$       & $-1/10$         &  $\pi/5$  & $\IQ(\sqrt{5})$             & A\\
\end{tabular}
\ \\
Mirror of $(R_1R_2R_3R_2^{-1})^5$\\
\begin{tabular}{c|c|c|c|c|c|c}
  $p$  & Order fix pt. stab & Signature     & $\chi^{orb}$    &  Area      & $\IQ(\tr\Gamma^2)$           & A/NA\\
  \hline
  4    &    4               & $(0;2,4,5)$   & $-1/20$         &  $\pi/10$  & $\IQ(\sqrt{5})$             & A\\
  5    &    2               & $(0;2,5,5)$   & $-1/10$         &  $\pi/5$   & $\IQ(\sqrt{5})$             & A\\
\end{tabular}
\end{table}
\ \\

\begin{table}\label{tab:E2}
\centering
  {\Large Groups $\Tc(p,{\bf E_2})$}\\
\begin{tabular}{c|c|c}
  Reflection     & Values of order $p$ & Generators\\          
  \hline
  $1$            &  3,4,6              & $(13)^2,(12)^2,(1312\bar1\bar3)^2,(1\bar2\bar1312)^2$ \\
  $2$            &  3,4,6              & $(21)^2,(2\bar3\bar2123)^2,(2\bar131)^3,(23)^3$\\
  $(12)^2$       &  6                  & $1,2$\\
  $(13)^2$       &  6                  & $1,3$\\
  $(123\bar2)^2$ &  6                  & $1,23\bar2$\\
  $(312\bar1)^3$ &  4,6                & $3,12\bar1$\\
  $Q^3$          &  3,4,6              & $Q,Q^{-1}2\bar3\bar2123\bar2$
\end{tabular}
\ \\
Mirror of $R = R_1$\\
\begin{tabular}{c|c|c|c|c|c|c}
  $p$  & Order fix pt. stab & Signature                    & $\chi^{orb}$    &  Area    & $\IQ(\tr\Gamma^2)$           & A/NA\\
  \hline
  3    & 3                  & $(0;2,6,6,6)$                & $-1$            &  $2\pi$  & $\IQ$                       & A\\
  4    & 4                  & $(0;2,\infty,\infty,\infty)$ & $-3/2$          &  $3\pi$  & $\IQ(\sqrt{3})$             & NA\\
  6    & 6                  & $(0;2,6,6,6)$                & $-1$            &  $2\pi$  & $\IQ$                       & A\\
\end{tabular}
\ \\
Mirror of $R = R_2$\\
\begin{tabular}{c|c|c|c|c|c|c}
  $p$  & Order fix pt. stab & Signature                      & $\chi^{orb}$    &  Area     & $\IQ(\tr\Gamma^2)$           & A/NA\\
  \hline
  3    & 3                  & $(0;2,6,6,6,\infty)$           & $-2$            &  $4\pi$   & $\IQ$                       & A\\
  4    & 4                  & $(0;4,4,\infty,\infty,\infty)$ & $-5/2$          &  $5\pi$   & $\IQ(\sqrt{3})$             & NA\\
  6    & 6                  & $(0;2,6,6,6,\infty)$           & $-2$            &  $4\pi$   & $\IQ$                       & A\\
\end{tabular}
\ \\
Mirror of $R = (R_1R_2)^2$\\
\begin{tabular}{c|c|c|c|c|c|c}
  $p$  & Order fix pt. stab & Signature         & $\chi^{orb}$    &  Area     & $\IQ(\tr\Gamma^2)$           & A/NA\\
  \hline 
  6    &  4                 & $(0;2,6,6)$       & $-1/6$          &  $\pi/3$  & $\IQ$                       & A\\
\end{tabular}
\ \\
\ \\
Mirror of $R = (R_1R_3)^2$\\
\begin{tabular}{c|c|c|c|c|c|c}
  $p$  & Order fix pt. stab & Signature         & $\chi^{orb}$    &  Area     & $\IQ(\tr\Gamma^2)$           & A/NA\\
  \hline 
  6    &  4                 & $(0;2,6,6)$       & $-1/6$          &  $\pi/3$  & $\IQ$                       & A\\
\end{tabular}
\ \\
Mirror of $(R_1R_2R_3R_2^{-1})^2$\\
\begin{tabular}{c|c|c|c|c|c|c}
  $p$  & Order fix pt. stab & Signature     & $\chi^{orb}$    &  Area     & $\IQ(\tr\Gamma^2)$           & A/NA\\
  \hline
  6    &    6               & $(0;2,6,6)$   & $-1/6$          &  $\pi/3$  & $\IQ$                       & A\\
\end{tabular}
\ \\
Mirror of $(R_3R_1R_2R_1^{-1})^3$\\
\begin{tabular}{c|c|c|c|c|c|c}
  $p$  & Order fix pt. stab & Signature     & $\chi^{orb}$   &  Area      & $\IQ(\tr\Gamma^2)$   & A/NA\\
  \hline
  4    &    4               & $(0;3,4,4)$   & $-1/6$         &  $\pi/3$   & $\IQ$               & A\\
  6    &    2               & $(0;3,6,6)$   & $-1/3$         &  $2\pi/3$  & $\IQ$               & A\\
\end{tabular}
\ \\
Mirror of $Q^3$\\
\begin{tabular}{c|c|c|c|c|c|c}
  $p$  & Order fix pt. stab & Signature         & $\chi^{orb}$    &  Area      & $\IQ(\tr\Gamma^2)$           & A/NA\\
  \hline
  3    &     2              & $(0;3,3,\infty)$  & $-1/3$          &  $2\pi/3$  & $\IQ$             & A\\
  4    &     2              & $(0;3,4,12)$      & $-1/3$          &  $2\pi/3$  & $\IQ(\sqrt{3})$   & A\\
  6    &     2              & $(0;3,6,6)$       & $-1/3$          &  $2\pi/3$  & $\IQ$             & A\\
\end{tabular}
\end{table}
\ \\

\begin{table}\label{tab:H1}
\centering
  {\Large Groups $\Tc(p,{\bf H_1})$}\\
\begin{tabular}{c|c|c}
  Reflection     & Values of order $p$ & Generators\\          
  \hline
  $1$            &  2                  & $(12)^2,(1312131213)^2,313212Q^3212313,(13Q^33)^2$ \\
  $Q^3$          &  2                  & $Q,2Q^32$\\
\end{tabular}
\ \\
Mirror of $R = R_1$\\
\begin{tabular}{c|c|c|c|c|c|c}
  $p$  & Order fix pt. stab & Signature                   & $\chi^{orb}$    &  Area      & $\IQ(\tr\Gamma^2)$           & A/NA\\
  \hline
  2    & 2                  & $(0;2,2,14,14)$             & $-6/7$         &  $12\pi/7$  & $\IQ(\cos\frac{2\pi}{7})$   & A\\
\end{tabular}
\ \\
Mirror of $Q^3$\\
\begin{tabular}{c|c|c|c|c|c|c}
  $p$  & Order fix pt. stab & Signature     & $\chi^{orb}$    &  Area       & $\IQ(\tr\Gamma^2)$           & A/NA\\
  \hline
  2    &    14              & $(0;2,3,14)$  & $-2/21$         &  $4\pi/21$  & $\IQ(\cos\frac{2\pi}{7})$   & A\\
\end{tabular}
\end{table}
\ \\

\begin{table}\label{tab:H2}
\centering
  {\Large Groups $\Tc(p,{\bf H_2})$}\\
\begin{tabular}{c|c|c}
  Reflection             & Values of order $p$ & Generators\\          
  \hline
  $1$                    &  2,3,5              & $(1\bar231\bar32)^5,(1\bar321\bar23)^5,(13)^3,(123\bar2)^5,(12)^5,2Q^3\bar2$ \\  
  $(12)^5$               &  5                  & $1,23\bar2$\\
  $(123\bar2)^5$         &  5                  & $1,23\bar2$\\
  $(3121\bar2\bar1)^5$ &  3,5                & $3,121\bar2\bar1$\\
  $Q^3$                  &  2                  & $Q^{-2}\bar2\bar12, 123\bar212\bar3\bar2\bar1$\\
\end{tabular}
\ \\
Mirror of $R = R_1$\\
\begin{tabular}{c|c|c|c|c|c|c}
  $p$  & Order fix pt. stab & Signature                   & $\chi^{orb}$    &  Area        & $\IQ(\tr\Gamma^2)$           & A/NA\\
  \hline
  2    & 2                  & $(0;5,10,10)$               & $-8/5$          &  $16\pi/5$   & $\IQ(\sqrt{5})$             & A\\
  3    & 3                  & $(0;2,2,5,6,6,15,15)$       & $-10/3$         &  $20\pi/3$   & $\IQ(\cos\frac{2\pi}{15})$  & NA\\
  5    & 5                  & $(0;2,2,5,5,5,10,10)$       & $-16/5$         &  $32\pi/5$   & $\IQ(\sqrt{3})$             & A\\
\end{tabular}
\ \\
Mirror of $(R_1R_2)^5$\\
\begin{tabular}{c|c|c|c|c|c|c}
  $p$  & Order fix pt. stab & Signature     & $\chi^{orb}$    &  Area      & $\IQ(\tr\Gamma^2)$ & A/NA\\
  \hline
  5    &    2               & $(0;2,5,5)$    & $-1/10$        &  $\pi/5$   & $\IQ(\sqrt{5})$   & A\\
\end{tabular}
\ \\
Mirror of $(R_1R_2R_3R_2^{-1})^5$\\
\begin{tabular}{c|c|c|c|c|c|c}
  $p$  & Order fix pt. stab & Signature     & $\chi^{orb}$    &  Area     & $\IQ(\tr\Gamma^2)$ & A/NA\\
  \hline
  5    &    2               & $(0;2,5,5)$   & $-1/10$         &  $\pi/5$  & $\IQ(\sqrt{5})$   & A\\
\end{tabular}
\ \\
Mirror of $(R_3R_1R_2R_1R_2^{-1}R_1^{-1})^5$\\
\begin{tabular}{c|c|c|c|c|c|c}
  $p$  & Order fix pt. stab & Signature     & $\chi^{orb}$    &  Area       & $\IQ(\tr\Gamma^2)$ & A/NA\\
  \hline
  3    &    3               & $(0;3,3,5)$   & $-2/15$         &  $4\pi/15$  & $\IQ(\sqrt{5})$   & A\\
  5    &    1               & $(0;5,5,5)$   & $-2/5$          &  $4\pi/5$   & $\IQ(\sqrt{5})$   & A\\
\end{tabular}
\ \\
Mirror of $Q^3$\\
\begin{tabular}{c|c|c|c|c|c|c}
  $p$  & Order fix pt. stab & Signature     & $\chi^{orb}$    &  Area        & $\IQ(\tr\Gamma^2)$           & A/NA\\
  \hline
  2    &    14              & $(0;2,3,10)$  & $-1/15$         &  $2\pi/15$   & $\IQ(\sqrt{5})$             & A\\
  3    &    14              & $(0;3,3,15)$  & $-4/15$         &  $8\pi/15$   & $\IQ(\cos\frac{2\pi}{15})$  & A\\
  5    &    14              & $(0;3,5,5)$   & $-4/15$         &  $8\pi/15$   & $\IQ(\sqrt{5})$             & A\\
\end{tabular}
\end{table}

\section{Hybrids}\label{sec:hybrids}

Using the results of section~\ref{sec:results}, we can obtain a
description of some lattices in the list as hybrids (in the sense
of~\cite{wells},~\cite{falbel-pasquinelli}). By this, we mean that we
find two orthogonal $\IC$-Fuchsian subgroups $\Sigma_1$ and $\Sigma_2$
of $\Gamma$ such that $\Sigma_1\cup\Sigma_2$ generates $\Gamma$.  Here
``orthogonal'' means that the complex lines preserved by the groups
are orthogonal.

In this paper, we do not address the issue of arithmeticity or
commensurability of the Fuchsian subgroups of complex hyperbolic
lattice triangle groups, but we do exhibit several pairs of orthogonal
$\IC$-Fuchsian subgroup that generate the ambient lattice.

\begin{thm}\label{thm:hybrids-s4c}
  For $p=8$ and $12$, the group $\Sc(p,\bar\sigma_4)$ is generated by
  the union of two mirror stabilizers; one can take the following
  pairs of reflections:
  \begin{itemize}
  \item $R_1$ and $(1 23\bar2)^3$,
  \item $R_1$ and $(1 \bar323)^3$,
  \item $R_1$ and $(1 312\bar1\bar3)^3$.
  \end{itemize}
\end{thm}

\begin{pf}
  Given the table on page~\pageref{tab:s4c}, the first statement is that
  $\Sc(p,\bar\sigma_4)$ is generated by $1$, $23\bar2$, $(12)^3$,
  $(13)^3$, $(1\bar323)^2$, $23\bar2(1J)^2$. This can be checked using
  computational group theory software (GAP or Magma), using the
  command to compute the index of a subgroup.

  The second statement is similar.

  The third one is slightly different, but similar as well; we check
  that $1$, $312\bar1\bar3$, $(12)^2$, $(13)^2$, $(123\bar2)^2$,
  $(1\bar323)^2$ and $23\bar2(1J)^2$ generate $\Sc(p,\bar\sigma_4)$. Once
  again, we check this using Magma.

  Note that these generation statements are actually true for
  $p=3,4,5,6$ as well, but in those cases, the elements
  $(123\bar2)^3$, $(1\bar323)^3$, $(1312\bar1\bar3)^3$ are complex
  reflections in points rather than lines.
\end{pf}

With a similar proof, we get the following generation results.
\begin{thm}\label{thm:hybrids-s5}
  The group $\Sc(2,\sigma_5)$ is generated by the union of stabilizers
  of the mirrors of $P^5$ and $2\bar3\bar2123\bar2$. For $p=4$, it is
  also generated by union of the stabilizers of the following pairs of
  reflections:
  \begin{itemize}
  \item $R_1$ and $(123\bar2)^5$;
  \item $R_1$ and $(1\bar323)^5$;
  \item $R_1$ and $(1\bar2\bar1312)^5$;
  \item $R_1$ and $(1312\bar1\bar3)^5$;
  \end{itemize}
\end{thm}

\begin{thm}\label{thm:hybrids-s10}
  For $p=4$, $5$ and $10$, the group $\Sc(p,\sigma_{10})$ generated by
  the union of the stabilizers of the mirrors of
  \begin{itemize}
  \item $R_1$ and $(12)^5$;
  \item $R_1$ and $(13)^5$.
  \end{itemize}
  For $p=10$, it is also generated by the stabilizers of
  \begin{itemize}
  \item $R_1$ and $(123\bar2)^3$;
  \item $R_1$ and $(1\bar323)^3$.
  \end{itemize}
\end{thm}

\begin{thm}
  For $p=4$ or $5$, the group $\Tc(p,{\bf S_2})$ is generated by the
  union of the stabilizers of $R_1$ and $(123\bar2)^5$.
\end{thm}

\end{document}